\documentclass[11pt]{amsart}

\usepackage{amsmath,amssymb,amsthm,mathtools}
\usepackage[margin=1.15in]{geometry}
\usepackage{microtype}
\usepackage{booktabs,array}
\usepackage[shortlabels]{enumitem}
\usepackage{xcolor}
\usepackage[colorlinks=true,linkcolor=blue!55!black,citecolor=green!45!black,urlcolor=blue!55!black]{hyperref}

\hypersetup{pdftitle={Functorial height and Orlov's rationality conjecture},pdfauthor={Xun Lin and Shizhuo Zhang}}
\newtheorem{theorem}{Theorem}[section]
\newtheorem{proposition}[theorem]{Proposition}
\newtheorem{lemma}[theorem]{Lemma}
\newtheorem{corollary}[theorem]{Corollary}

\theoremstyle{definition}
\newtheorem{definition}[theorem]{Definition}
\newtheorem{example}[theorem]{Example}
\theoremstyle{remark}
\newtheorem{remark}[theorem]{Remark}
\numberwithin{equation}{section}

\newcommand{\C}{\mathbb{C}}
\newcommand{\Z}{\mathbb{Z}}
\newcommand{\PP}{\mathbb{P}}
\newcommand{\cO}{\mathcal{O}}
\newcommand{\cA}{\mathcal{A}}
\newcommand{\cB}{\mathcal{B}}
\newcommand{\cC}{\mathcal{C}}
\newcommand{\cD}{\mathcal{D}}
\newcommand{\cE}{\mathcal{E}}
\newcommand{\cH}{\mathcal{H}}
\newcommand{\cP}{\mathcal{P}}
\newcommand{\cQ}{\mathcal{Q}}
\newcommand{\cX}{\mathcal{X}}
\newcommand{\bE}{\mathbf{E}}
\newcommand{\bT}{\mathbf{T}}
\newcommand{\Db}{\mathrm{D}^{\mathrm{b}}}
\newcommand{\Perf}{\mathrm{Perf}}
\newcommand{\id}{\mathrm{id}}
\newcommand{\ph}{\mathrm{ph}}
\newcommand{\Kuz}{\mathrm{Kuz}}
\newcommand{\op}{\mathrm{op}}
\DeclareMathOperator{\Hom}{Hom}
\DeclareMathOperator{\RHom}{RHom}
\DeclareMathOperator{\Ext}{Ext}
\DeclareMathOperator{\End}{End}
\DeclareMathOperator{\Tr}{Tr}
\DeclareMathOperator{\tr}{tr}
\DeclareMathOperator{\Tot}{Tot}
\DeclareMathOperator{\Cone}{Cone}
\DeclareMathOperator{\coker}{coker}
\DeclareMathOperator{\HH}{HH}
\DeclareMathOperator{\NHH}{NHH}
\DeclareMathOperator{\Pic}{Pic}
\DeclareMathOperator{\Aut}{Aut}
\DeclareMathOperator{\Bl}{Bl}
\DeclareMathOperator{\Conf}{Conf}
\DeclareMathOperator{\Supp}{Supp}
\DeclareMathOperator{\thick}{thick}
\DeclareMathOperator{\rk}{rk}
\DeclareMathOperator{\codim}{codim}
\DeclareMathOperator{\gr}{gr}

\begin{document}

\title{Functorial height and Orlov's rationality conjecture}

\author{Xun Lin}
\address{School of Science and Engineering, The Chinese University of Hong Kong, Shenzhen, Shenzhen 518172, P.~R.~China}
\email{linxun@cuhk.edu.cn}

\author{Shizhuo Zhang}
\address{School of Mathematics, Sun Yat-sen University, Guangzhou 510275, P.~R.~China}
\email{zhangshzh28@mail.sysu.edu.cn}

\date{November 17, 2026}

\subjclass[2020]{Primary 14F08; Secondary 14J26, 16E40, 18G80}
\keywords{Exceptional collection, functorial height, normal Hochschild cohomology, rational surface, Fourier--Mukai kernel, phantom category}

\begin{abstract}
An exceptional collection in the derived category of a smooth projective variety may fail to generate the whole category, and this failure cannot be detected by K-theory or by Hochschild homology. Kuznetsov's height detects it through Hochschild cohomology, that is, through natural transformations from the identity functor to its shifts. We replace the identity functor by an arbitrary endofunctor $F$ and define the height of a collection with coefficient $F$. It controls, in a range of degrees, the restriction of natural transformations $\mathrm{id}\to F$ to the orthogonal complement of the collection, and it is computed from the collection by a finite complex. For collections of line bundles on surfaces whose Hom divisors have nonnegative canonical degree, the tensor product with a skyscraper sheaf has height exactly two, and the projections of skyscraper sheaves form a two-dimensional family of point-like objects in the orthogonal complement. We apply these tools to Orlov's rationality conjecture. A smooth projective complex surface is rational if its derived category has a full exceptional collection in which all members except possibly one are line bundles, or one that consists of at most three completely orthogonal blocks of arbitrary objects; we also obtain rationality in several cases with four blocks or with a hereditary component. If $S\times\mathbb P^1$ has a full exceptional collection of the form $(M_i\boxtimes\mathcal O(a_i))$, then the $M_i$ with relative degrees of a fixed parity form a full exceptional collection on $S$; consequently $S$ is rational whenever $S\times\mathbb P^1$ has a full exceptional collection of line bundles. Finally, we show that the general case of the conjecture for surfaces is equivalent to the surjectivity and injectivity of two maps in Hochschild cohomology.
\end{abstract}

\maketitle

\tableofcontents

\section{Introduction}\label{sec:intro}

\subsection{Exceptional collections and the fullness problem}

A full exceptional collection is the simplest possible description of the derived category of a variety: it presents $\Db(X)$ as built from finitely many copies of the derived category of a point, glued by finite-dimensional data. Projective spaces, quadrics, Grassmannians and rational surfaces have such collections, and a variety that has one behaves, from the point of view of its derived category, like a finite-dimensional algebra.

Two questions arise at once. The first is \emph{recognition}: given an exceptional collection, how does one decide whether it is full? The obvious necessary condition, that its classes form a basis of the Grothendieck group, is far from sufficient. Phantom categories, admissible subcategories with vanishing Grothendieck group and Hochschild homology, occur as orthogonals of exceptional collections on surfaces of general type \cite{BGKS15,GO13}, and Krah has shown that even a blowup of the projective plane in ten general points carries an exceptional collection of maximal length whose orthogonal is a nonzero phantom \cite{Kra24}. So the Grothendieck group and Hochschild homology cannot detect fullness.

The second question is \emph{geometric}: what does the existence of a full exceptional collection say about the variety? The surface form of Orlov's rationality conjecture predicts that a smooth projective complex surface with a full exceptional collection is rational; see \cite[Introduction]{Via17}. For collections of line bundles on surfaces the conjecture has a long history, through the constructions of Hille and Perling on rational surfaces \cite{HP11}, the theorem of Brown and Shipman for full strong collections \cite{BS17}, results for small Picard rank \cite{Zha19}, and the classification of cyclic strong collections \cite{EXZ21}; it was recently proved for all full exceptional collections of line bundles on surfaces in \cite{LZ26}.

The sharpest general tool for the first question is Kuznetsov's \emph{height} \cite{Kuz15}. For an exceptional collection generating $\cB\subset\Db(X)$, with orthogonal $\cA$, normal Hochschild cohomology measures the difference between the Hochschild cohomology of $X$ and that of $\cA$; its first nonzero degree is the height, and it controls the range in which restriction $\HH^\bullet(X)\to\HH^\bullet(\cA)$ is an isomorphism. A collection of positive height cannot be full. Height is a non-additive invariant, and it sees quasiphantoms and phantoms that additive invariants miss; when the height is zero, however, its value alone gives no information (Appendix~\ref{app:ordinary}).

Hochschild cohomology is the space of natural transformations from the identity functor to its shifts. The starting point of this paper is that there is no reason to use only the identity. Any perfect endofunctor $F$ of $\Db(X)$, such as tensoring with a skyscraper sheaf, a line bundle, a sheaf on a curve, or the graph of an automorphism, gives a complex $\RHom(\id,F)$ of natural transformations, and each of them can be restricted to the orthogonal of a collection. Different coefficients probe the orthogonal in different ways, and several coefficients together retain the compositions between them. We develop this idea systematically and apply it to Orlov's conjecture.

\subsection{Functorial height}

Let $\cD$ be a smooth proper dg category over a field $k$ with a semiorthogonal decomposition
\[
  \cD=\langle\cA,\cB\rangle,\qquad i\colon\cA\hookrightarrow\cD,\qquad j\colon\cB\hookrightarrow\cD .
\]
Write $S_\cD$ for its Serre functor and $Q_\cB=jj^!$ for the projector onto $\cB$. For a perfect endofunctor $F$, define the \emph{ambient} and \emph{normal} coefficient complexes
\[
  C_F(\cD)=\RHom(\id_\cD,F),\qquad N_F(\cB,\cD)=\Tr_\cD\bigl(Q_\cB S_\cD^{-1}F\bigr),
\]
where $\Tr$ is the categorical trace. Their initial degrees are $\alpha_F=\inf C_F$ and $h_F=\inf N_F$, with $\inf 0=+\infty$. We call $h_F$ the \emph{functorial height} of $\cB$ with coefficient $F$; for $F=\id$ it is Kuznetsov's height.

\begin{theorem}[Functorial localization]\label{thm:intro-localization}
There is a natural distinguished triangle
\begin{equation}\label{eq:intro-triangle}
  N_F(\cB,\cD)\longrightarrow C_F(\cD)\xrightarrow{\ \mathrm{res}_F\ }\RHom\bigl(\id_\cA,i^!Fi\bigr).
\end{equation}
The restriction map is an isomorphism in degrees $t\le h_F-2$ and injective in degree $h_F-1$. In particular, if $C_F(\cD)\ne0$ and $h_F>\alpha_F$, then $\cA\ne0$.
\end{theorem}

This is the coefficient version of Kuznetsov's localization theorem \cite[Theorems~3.3 and~4.2]{Kuz15}. When $\cB$ is generated by an exceptional collection, $N_F$ is computed by a finite \emph{normal bar complex} (Proposition~\ref{prop:bar}), whose first page gives an elementary lower bound, the \emph{pseudoheight}, while its differentials can give stronger vanishing.

Coefficients can see more than the identity. In Appendix~\ref{app:fivefold} we construct a rational fivefold with an exceptional collection of maximal length whose orthogonal is a universal phantom and whose ordinary height is zero; a point-type coefficient on the same collection has $h_F=2>0=\alpha_F$ (Proposition~\ref{prop:two-heights}). In this example, and in the example of a Barlow surface times $\PP^1$ (Example~\ref{ex:barlow}), the nonfullness is known from the construction; what the functorial height adds is a certificate computed from the collection alone.

\subsection{Geometric and simultaneous coefficients}

For $\cD=\Perf(X)$ with $X$ smooth projective, perfect endofunctors are Fourier--Mukai kernels $K\in\Perf(X\times X)$, and \eqref{eq:intro-triangle} becomes a triangle of Hom complexes on $X\times X$. The simplest kernels are external products: for perfect objects $T$ and $Y$ and $K_{T,Y}=T^\vee\boxtimes S_X(Y)$, the triangle is
\begin{equation}\label{eq:intro-external}
  \RHom_X(T,Q_\cB Y)\longrightarrow\RHom_X(T,Y)\longrightarrow\RHom_\cA(i^*T,i^*Y).
\end{equation}
So the height of $K_{T,Y}$ measures the range in which morphisms from $T$ to $Y$ are preserved by projection to the orthogonal. External coefficients together detect every nonzero orthogonal (Corollary~\ref{cor:detect}); what makes the method effective is a height estimate that can be computed without knowing the orthogonal.

For several probes $T_1,\dots,T_m$ the kernels $K_{ab}=T_a^\vee\boxtimes S_X(T_b)$ give a matrix of heights $(h_{ab})$. Their comparison maps all come from the single functor $i^*$, so they are compatible with composition: with $T=\bigoplus_aT_a$ they form a map of dg algebras
\[
  \RHom_X(T,T)\longrightarrow\RHom_\cA(i^*T,i^*T),\qquad\text{with fibre }\ \RHom_X(T,Q_\cB T).
\]
More generally, for arbitrary kernels $K_a$ with sum $K_\Sigma$ and $\cH=\RHom(K_\Sigma,K_\Sigma)$, the three coefficient complexes are dg modules over $\cH$, and the normal complex of any kernel built from the $K_a$ is a derived tensor product over $\cH$ (Proposition~\ref{prop:coefficient-algebra}). This gives cancellation criteria: a morphism of kernels can kill low normal cohomology without killing the corresponding ambient cohomology (Proposition~\ref{prop:cancellation}).

\subsection{Line bundles on surfaces}

For exceptional collections of line bundles on surfaces the heights of point coefficients can be computed exactly. Let $\bE=(L_1,\dots,L_N)$ be such a collection on a smooth projective complex surface $X$, and assume that the zero divisors of all nonzero maps $L_i\to L_j$, $i<j$, have components of nonnegative canonical degree; this holds, for example, whenever $K_X$ is nef. The \emph{Hom graph} of $\bE$ has the vertices $1,\dots,N$ and an edge $\{i,j\}$ whenever $\Hom(L_i,L_j)\ne0$; let $c$ be the number of its connected components.

\begin{theorem}[Theorems~\ref{thm:point-height} and~\ref{thm:point-like}]\label{thm:intro-point-like}
Let $\Db(X)=\langle\cA,\bE\rangle$ with a nonempty line-bundle collection as above. For $x,y$ off the Hom divisors the height of $K_{\cO_x,\cO_y}$ is two. Write $m_{xy}=\dim H^2N_{K_{\cO_x,\cO_y}}$ and $W_x=i^*\cO_x$. Then $m_{xy}\ge c$ and
\[
  \End(W_x)=\C,\qquad T_xX\hookrightarrow\Ext^1(W_x,W_x),
  \qquad\dim\Ext^1(W_x,W_x)=1+m_{xx}.
\]
For $x\ne y$, $\Ext^{\le0}(W_x,W_y)=0$ and $\dim\Ext^1(W_x,W_y)=m_{xy}$.
\end{theorem}
The lower bound comes from the length-zero part of the normal bar filtration. Positive bar lengths can contribute in the same degree, as Example~\ref{ex:point-multiplicity} shows. For a hereditary orthogonal the Euler matrix gives an exact formula (Corollary~\ref{cor:hereditary-surface}).

\subsection{Orlov's rationality conjecture}

Our main applications concern Orlov's conjecture for surfaces. A \emph{block} is an exceptional collection whose distinct members are completely orthogonal, that is, have vanishing derived Hom in both directions. In the following theorem the blocks are consecutive, and morphisms between different blocks may be nonzero in any degree.

\begin{theorem}[Theorems~\ref{thm:any-three-blocks} and~\ref{thm:rational}]\label{thm:intro-rational}
Let $X$ be a smooth connected projective complex surface. Suppose that $\Db(X)$ has a full exceptional collection satisfying one of the following conditions:
\begin{enumerate}[(i)]
  \item it consists of at most three blocks of arbitrary objects;
  \item all its members except possibly one are shifted line bundles.
\end{enumerate}
Then $X$ is rational.
\end{theorem}

The two parts are proved by different methods. For (i), a nonrational surface with a full exceptional collection has $\HH^1(X)=0$. A composable word of nonidentity morphisms through three blocks has at most two letters, so the minimal model of the collection is formal. Its Euler derivation is then inner, and shifting the individual objects makes the collection strong. A strong collection in three blocks has $\HH^4=0$, whereas $\HH^4(X)$ is dual to $H^0(X,2K_X)$, so Castelnuovo's criterion gives rationality. For (ii), the point coefficient rules out an orthogonal generated by one exceptional object, under the canonical nonnegativity condition; moving the blowup centres and contracting exceptional curves reduces the general case to this one.

The same ideas give further cases. The block argument extends to four blocks if the collection is formal, or if the forward Ext spaces between adjacent blocks can be shifted into degree zero (Theorem~\ref{thm:four-pure-adjacent}), and to some nonformal collections with a prescribed higher product (Theorem~\ref{thm:fixed-arity}). A different argument, based on Kuznetsov's gluing sequence for Hochschild cohomology, applies when one component is hereditary: a component equivalent to the derived category of a finite acyclic quiver, of arbitrary representation type, together with one block of arbitrary exceptional objects forces rationality (Theorem~\ref{thm:gluing-hereditary}). We also treat pairs of arbitrary objects next to a collection of line bundles (Section~\ref{subsec:pair-profile}), exceptional generators with square-zero radical (Theorem~\ref{thm:radical-square}), and components with few isomorphism classes of objects (Theorem~\ref{thm:discrete-projections}). For an arbitrary full exceptional collection, the gluing sequence shows that rationality is equivalent to the surjectivity of one map and the injectivity of another in Hochschild cohomology (Proposition~\ref{prop:last-generator}).

\subsection{Products with the projective line}

\begin{theorem}[Theorem~\ref{thm:parity-descent}]\label{thm:intro-product}
Let $S$ be a smooth projective complex variety. If $(M_i\boxtimes\cO_{\PP^1}(a_i))$ is a full exceptional collection on $S\times\PP^1$, then the $M_i$ with even $a_i$, in their original order, form a full exceptional collection on $S$, and so do the $M_i$ with odd $a_i$. In particular, if $S$ is a connected surface and $S\times\PP^1$ has a full exceptional collection of line bundles, then $S$ and $S\times\PP^1$ are rational.
\end{theorem}

The proof realizes the spaces of sections of $\cO_{\PP^1}(n)$ as representations of the Sweedler Hopf algebra. For positive odd $n$ these representations are projective, whereas the trivial representation has infinite projective dimension, and this distinction persists through the evaluation triangles of a full exceptional collection.

\subsection{Relation to previous work}

The functorial localization theorem and the normal bar complex are the coefficient versions of the corresponding results of Kuznetsov \cite{Kuz15}. The proof of part~(ii) of Theorem~\ref{thm:intro-rational} relies on \cite{LZ26} at the following points: the structure of the algebra of morphisms between the line bundles \cite[Lemma~3.2 and Proposition~3.8]{LZ26}, used in Proposition~\ref{prop:LZ}; the filtration by positive-degree arrows \cite[Lemmas~3.9--3.15]{LZ26}, which we adapt to a two-sided bar complex in Lemma~\ref{lem:connectivity}; the family of blowups and the choice of general centres \cite[Lemmas~3.20--3.22]{LZ26}; the analysis of left-orthogonal divisors on a blowup \cite[Lemmas~3.24--3.26 and Proposition~3.28]{LZ26}; and the mutations of line bundles in the contraction step \cite[Lemmas~3.29--3.30 and Proposition~3.31]{LZ26}. We also use Elagin's criterion for left-orthogonal divisors \cite{Ela16}, Hu's deformation and openness results for exceptional collections \cite{Hu19}, Kuznetsov's gluing sequence for Hochschild cohomology \cite{Kuz09}, the purity theorem of Bridgeland and Stern \cite{BSt10}, Pirozhkov's results on admissible subcategories supported on curves \cite{PirD,PirC}, and an action of the Sweedler Hopf algebra on a polynomial ring due to Walton and Witherspoon \cite{WalWit14}. The new ingredients are the functorial height and its point calculation, the exclusion of an orthogonal generated by one exceptional object, the treatment of arbitrary members under deformation and contraction, the Euler-derivation argument for blocks, the use of hereditary gluing data, and the parity descent from $S\times\PP^1$.

\subsection{Organization}

Sections~\ref{sec:localization}--\ref{sec:special} develop functorial localization, the normal bar complex, Fourier--Mukai and simultaneous coefficients, and a collection of special coefficients. Section~\ref{sec:point-height} computes the heights of point coefficients for line bundles on surfaces, and Section~\ref{sec:rationality} proves part~(ii) of Theorem~\ref{thm:intro-rational} and treats pairs of arbitrary objects. Sections~\ref{sec:three-object-blocks} and~\ref{sec:hereditary-new} treat orthogonal blocks and hereditary components. Section~\ref{sec:products} proves the descent from a product with a projective line, and Section~\ref{sec:remaining-glue} formulates the general case in terms of two gluing maps. The appendices compare ordinary and functorial height.

\subsection{Conventions}

Sections~\ref{sec:localization}--\ref{sec:simultaneous} work over an arbitrary field $k$, except where a variety over $\C$ is specified; in Section~\ref{sec:special} characteristic zero is assumed where indicated; from Section~\ref{sec:point-height} on the base field is $\C$. All functors are enhanced and all gradings are cohomological: $H^t(C[m])=H^{t+m}(C)$. Functors compose from right to left. We write $\HH^t$ for Hochschild cohomology and $\HH_t$ for Hochschild homology. An exceptional object satisfies $\RHom(E,E)\simeq k$, and an exceptional collection $(E_1,\dots,E_N)$ also satisfies $\RHom(E_j,E_i)=0$ for $j>i$. Fullness means generation as a thick subcategory. We use $\cA=\cB^\perp=\{D:\RHom(B,D)=0\text{ for all }B\in\cB\}$, so that $\cD=\langle\cA,\cB\rangle$.

\subsection*{Use of AI tools}
The main ideas of this article are due to the authors. The notion of functorial height and its basic properties were obtained by the authors before AI tools became available. The details of the arguments were obtained with the assistance of AI tools, and possible applications, including Orlov's rationality conjecture, the existence of phantom categories, and categorical Torelli theorems, were also found with AI assistance. The authors are responsible for the correctness of the article.

\subsection*{Acknowledgements}
Shizhuo Zhang is supported by the Sun Yat-sen University Starting Grant 34000-12256019. Xun Lin is supported by the Chinese University of Hong Kong(Shenzhen) Start Up Grant.

\section{Normal localization and functorial height}\label{sec:localization}

\subsection{Enhanced coefficients and the trace}

Let $k$ be a field. A dg category is \emph{proper} if all its morphism complexes are perfect over $k$, and \emph{smooth} if its diagonal bimodule is perfect. We work with small, pretriangulated, idempotent-complete smooth proper dg categories. A \emph{perfect endofunctor} is a quasi-functor represented by a perfect bimodule, and derived morphisms of endofunctors are computed in the dg category of bimodules.

For a perfect endofunctor $G$, its categorical trace is the derived tensor product of its bimodule with the diagonal. Tensor--Hom duality for the perfect diagonal gives
\begin{equation}\label{eq:trace}
  \Tr_\cD(G)\simeq\RHom(S_\cD^{-1},G)\simeq\RHom(\id_\cD,S_\cD G),
\end{equation}
where the inverse dualizing bimodule represents $S_\cD^{-1}$ \cite[Theorems~4.2--4.5]{Shk07}, and the second equivalence applies $S_\cD$ to both arguments. We also use cyclicity,
\begin{equation}\label{eq:cyclicity}
  \Tr_\cD(UV)\simeq\Tr_\cC(VU)
\end{equation}
for perfect functors $U\colon\cC\to\cD$ and $V\colon\cD\to\cC$ \cite{PS13,HSS17}. These equivalences are compatible with natural transformations, so they identify the maps in a triangle of traces as well as its objects.

\subsection{The localization triangle}

Fix an admissible semiorthogonal decomposition $\cD=\langle\cA,\cB\rangle$, with inclusions $i,j$ and adjoints $i^*\dashv i\dashv i^!$, $j^*\dashv j\dashv j^!$. The projection triangle is
\begin{equation}\label{eq:projection}
  Q_\cB:=jj^!\longrightarrow\id_\cD\longrightarrow P_\cA:=ii^*,
\end{equation}
whose arrows are the counit of $j\dashv j^!$ and the unit of $i^*\dashv i$.

\begin{definition}\label{def:complexes}
For a perfect endofunctor $F$, the \emph{ambient}, \emph{normal} and \emph{restricted} coefficient complexes are
\begin{equation}\label{eq:complexes}
\begin{gathered}
  C_F(\cD)=\RHom(\id_\cD,F),\qquad N_F(\cB,\cD)=\Tr_\cD\bigl(Q_\cB S_\cD^{-1}F\bigr),\\
  R_F(\cA\subset\cD)=\RHom(\id_\cA,i^!Fi).
\end{gathered}
\end{equation}
\end{definition}

\begin{theorem}\label{thm:localization}
The complexes \eqref{eq:complexes} fit into a natural distinguished triangle
\[
  N_F(\cB,\cD)\longrightarrow C_F(\cD)\xrightarrow{\ \mathrm{res}_F\ }R_F(\cA\subset\cD).
\]
On cohomology, the restriction of $\eta\colon\id_\cD\to F[t]$ is $i^!\eta i\colon\id_\cA\to i^!Fi[t]$. For $\cD=\Perf(X)$ and $F=\id$ this is \cite[Theorem~3.3]{Kuz15}.
\end{theorem}

\begin{proof}
Compose \eqref{eq:projection} on the right with $S_\cD^{-1}F$ and take the trace; exactness of the derived tensor product gives a distinguished triangle. Its first term is the normal complex by definition, and its second term is the ambient complex by \eqref{eq:trace}. For the third term, cyclicity gives
\[
  \Tr_\cD(ii^*S_\cD^{-1}F)\simeq\Tr_\cA(i^*S_\cD^{-1}Fi)\simeq\RHom(\id_\cA,S_\cA i^*S_\cD^{-1}Fi).
\]
There is a natural equivalence $S_\cA i^*S_\cD^{-1}\simeq i^!$: for $A\in\cA$ and $D\in\cD$, Serre duality and adjunction give
\[
  \RHom_\cA(A,S_\cA i^*S_\cD^{-1}D)\simeq\RHom_\cA(i^*S_\cD^{-1}D,A)^\vee\simeq\RHom_\cD(S_\cD^{-1}D,iA)^\vee\simeq\RHom_\cD(iA,D),
\]
and the Yoneda lemma identifies this functor with the right adjoint $i^!$. Under these adjunctions, the unit $\id\to ii^*$ of \eqref{eq:projection} sends a transformation $\eta$ to $i^!\eta i$, which proves the description of the restriction map.
\end{proof}

\subsection{Initial degrees and restriction}

For a perfect $k$-complex $C$, put $\inf C=\inf\{t:H^t(C)\ne0\}$, with $\inf0=+\infty$.

\begin{definition}\label{def:height}
The \emph{$F$-height} of $\cB\subset\cD$ is $h_F(\cB\subset\cD)=\inf N_F(\cB,\cD)$. Put $\alpha_F(\cD)=\inf C_F(\cD)$. If $C_F(\cD)\ne0$, the \emph{excess height} is
\[
  \varepsilon_F(\cB\subset\cD)=h_F(\cB\subset\cD)-\alpha_F(\cD).
\]
For an exceptional collection $\bE$ generating $\cB$ we also write $h_F(\bE)$.
\end{definition}

The next proposition is the coefficient analogue of \cite[Theorem~4.2]{Kuz15}; together with Theorem~\ref{thm:localization} it proves Theorem~\ref{thm:intro-localization}.

\begin{proposition}\label{prop:range}
Put $h=h_F(\cB\subset\cD)$. If $h$ is finite, restriction is an isomorphism in degrees $t\le h-2$ and injective in degree $h-1$; if $h=+\infty$, it is an isomorphism in every degree. In either case, every nonzero class in $H^tC_F(\cD)$ with $t<h$ has nonzero restriction; in particular, if such a class exists then $\cA\ne0$.
\end{proposition}

\begin{proof}
In the exact sequence $H^tN_F\to H^tC_F\to H^tR_F\to H^{t+1}N_F$, both outer terms vanish for $t\le h-2$, and the left term vanishes whenever $t<h$. If $\cA=0$ then $R_F=0$, so no nonzero ambient class can have nonzero restriction.
\end{proof}

\begin{proposition}\label{prop:formal}
The normal and ambient complexes satisfy
\begin{gather*}
  \cB=\cD\ \Longrightarrow\ N_F\simeq C_F,\ \ h_F=\alpha_F;\qquad h_{F[m]}=h_F-m,\quad\alpha_{F[m]}=\alpha_F-m;\\
  h_{F\oplus G}=\min(h_F,h_G),\qquad\alpha_{F\oplus G}=\min(\alpha_F,\alpha_G).
\end{gather*}
The excess is invariant under shifts of the coefficient. Morita equivalences transport the three coefficient complexes and their triangle, and a distinguished triangle of coefficients induces distinguished triangles of all three complexes.
\end{proposition}

\begin{proof}
If $\cB=\cD$, its projector is the identity, and \eqref{eq:trace} applies. The three constructions are exact in $F$, which gives the assertions on shifts, sums and triangles. A Morita equivalence transports the diagonal, the inverse dualizing bimodule and the adjunction maps, hence the coefficient triangle.
\end{proof}

In particular, height is invariant under mutations of the exceptional generators of a fixed subcategory, extending \cite[Corollary~4.3]{Kuz15}, and positive excess obstructs fullness, as positive ordinary height does in \cite[Proposition~6.1]{Kuz15}.

\section{The finite normal bar complex}\label{sec:bar}

Let $\bE=(E_1,\dots,E_N)$ be an exceptional collection in $\cD$, and put $\cB=\langle\bE\rangle$. For a coefficient $F$, abbreviate $H_F=S_\cD^{-1}F$. Its \emph{normal bimodule} on the exceptional objects has entries
\begin{equation}\label{eq:normal-bimodule}
  M_F(E_i,E_j)=\RHom_\cD\bigl(E_j,H_F(E_i)\bigr).
\end{equation}
For $a\colon E_i\to E_{i'}$, $c\colon E_{j'}\to E_j$ and $\xi\colon E_j\to H_F(E_i)$, the two actions are
\begin{equation}\label{eq:two-actions}
  a\cdot\xi=H_F(a)\circ\xi,\qquad \xi\cdot c=\xi\circ c .
\end{equation}
Both actions enter the computation of the height; the graded coefficient spaces alone do not determine it.

\subsection{The directed dg model}

Choose representatives of the $E_i$ in a dg enhancement. Keep the complexes $\RHom(E_i,E_j)$ for $i<j$, replace each diagonal entry by $k\,\id_{E_i}$, and set the reverse entries equal to zero. The inclusion into the full dg subcategory is a quasi-equivalence, since exceptionality identifies the diagonal cohomology with $k$ and semiorthogonality makes the reverse complexes acyclic. The algebra $\cE$ of this directed dg category decomposes as
\begin{equation}\label{eq:directed}
  \cE=S_0\oplus J,\qquad S_0=\bigoplus_{i=1}^Nke_i,\qquad J=\bigoplus_{i<j}e_j\cE e_i,
\end{equation}
with the matrix convention $e_j\cE e_i=\RHom(E_i,E_j)$. The coefficient bimodule $T_F$ is placed in the opposite matrix positions, $e_iT_Fe_j=M_F(E_i,E_j)$. For an $S_0$-bimodule $W$, write $\tr_{S_0}W=\bigoplus_ie_iWe_i$.

\begin{proposition}\label{prop:bar}
There is a finite normal bar model
\begin{equation}\label{eq:bar}
  N_F(\cB,\cD)\simeq\Tot\Bigl(\bigoplus_{p=0}^{N-1}\tr_{S_0}\bigl(T_F\otimes_{S_0}J^{\otimes_{S_0}p}\bigr)[p]\Bigr).
\end{equation}
Its bar-length filtration gives a convergent spectral sequence
\begin{equation}\label{eq:bar-ss}
\begin{aligned}
  E_1^{-p,q}(F)=\bigoplus_{\substack{i_0<\dots<i_p\\ a_0+\dots+a_p=q}}&\Ext^{a_0}(E_{i_0},E_{i_1})\otimes_k\dots\otimes_k\Ext^{a_{p-1}}(E_{i_{p-1}},E_{i_p})\\[-0.8em]
  &\otimes_k\Ext^{a_p}\bigl(E_{i_p},H_F(E_{i_0})\bigr)\ \Longrightarrow\ H^{q-p}N_F(\cB,\cD),
\end{aligned}
\end{equation}
where for $p=0$ the internal tensor product is $k$. The differentials have bidegree $d_r\colon E_r^{-p,q}\to E_r^{-p+r,q-r+1}$, and the first differential is given by Yoneda products and the two actions \eqref{eq:two-actions}. For $F=\id$ this is the reduced normal bar complex of \cite[Proposition~3.7]{Kuz15}.
\end{proposition}

\begin{proof}
Cyclicity identifies $\Tr_\cD(jj^!H_F)\simeq\Tr_\cB(j^!H_Fj)$. In the chosen generators, the bimodule of the last functor is $T_F$, and its trace is the derived tensor product with the diagonal of $\cE$. The normalized relative bar resolution of the diagonal over the separable algebra $S_0$ has terms $\cE\otimes_{S_0}J^{\otimes_{S_0}p}\otimes_{S_0}\cE$, which are homotopically projective over $\cE^e$, and the bar augmentation is a quasi-isomorphism (see the explicit contraction in Lemma~\ref{lem:bar-resolution}). A nonzero tensor of $p$ factors of $J$ traverses $p$ strictly ordered arrows, so $J^{\otimes_{S_0}p}=0$ for $p\ge N$. Taking the relative trace after tensoring the resolution with $T_F$ proves \eqref{eq:bar}. Taking internal cohomology gives \eqref{eq:bar-ss} by the K\"unneth formula over $k$; the shift $[p]$ subtracts $p$ from the internal degree $q$, and the filtration is finite, so the spectral sequence converges. On internal cohomology the bar faces become Yoneda composition and the endpoint actions. The construction uses the dg algebra $\cE$ itself; in a minimal $A_\infty$ model, the higher products of the algebra and of the coefficient module contribute to the later differentials.
\end{proof}

\subsection{The differential}

We use unshifted tensor factors and shift the whole length-$p$ summand by $[p]$. First let $A$ be a graded algebra and $T$ a graded bimodule with zero differential, and write $C_p=\tr_{S_0}(T\otimes_{S_0}J^{\otimes_{S_0}p})$. For
\[
  z=x\otimes a_1\otimes\dots\otimes a_p,\qquad |x|=t,\quad|a_\nu|=d_\nu,
\]
the internal degree is $t+\sum_\nu d_\nu$ and the total degree is this number minus $p$. For $p\ge1$ define
\begin{equation}\label{eq:faces}
\begin{aligned}
  b^{(p)}_0(z)&=xa_1\otimes a_2\otimes\dots\otimes a_p,\\
  b^{(p)}_\nu(z)&=x\otimes\dots\otimes(a_\nu a_{\nu+1})\otimes\dots\otimes a_p\qquad(0<\nu<p),\\
  b^{(p)}_p(z)&=(-1)^{d_p(t+d_1+\dots+d_{p-1})}(a_px)\otimes a_1\otimes\dots\otimes a_{p-1},
\end{aligned}
\end{equation}
and
\begin{equation}\label{eq:b}
  b|_{C_p[p]}=\sum_{\nu=0}^p(-1)^\nu b^{(p)}_\nu,\qquad b|_{C_0}=0 .
\end{equation}
Here $xa_1$ is the right action on the coefficient and $a_px$ is the left action; the sign in the last face moves the final arrow past all preceding factors. In length one this is the graded commutator $b(x\otimes a)=xa-(-1)^{|a||x|}ax$. Associativity and the bimodule identities give $b^{(p-1)}_ib^{(p)}_j=b^{(p-1)}_{j-1}b^{(p)}_i$ for $0\le i<j\le p$, so the terms of $b^2$ cancel in pairs. Each face preserves the internal degree and decreases the length by one, so $b$ has total degree one.

For a path $i_0<\dots<i_p$, a tensor in matrix order has $x\in e_{i_0}Te_{i_p}$ and $a_\nu\in e_{i_{p-\nu+1}}Je_{i_{p-\nu}}$. Consequently $xa_1\in e_{i_0}Te_{i_{p-1}}$ and $a_px\in e_{i_1}Te_{i_p}$, and for the normal coefficient
\[
  xa_1=x\circ a_1,\qquad a_px=H_F(a_p)\circ x .
\]
The forward order in \eqref{eq:bar-ss} is obtained by the Koszul symmetry of the tensor factors.

For a dg algebra and a dg bimodule, add the tensor differential
\[
  d_{\mathrm{int}}z=d_Tx\otimes a_1\otimes\dots\otimes a_p+\sum_{\nu=1}^p(-1)^{t+d_1+\dots+d_{\nu-1}}x\otimes\dots\otimes d_Aa_\nu\otimes\dots\otimes a_p .
\]
With the whole-summand shift convention the total differential is
\begin{equation}\label{eq:dtot}
  d_{\mathrm{tot}}|_{C_p[p]}=(-1)^pd_{\mathrm{int}}+b .
\end{equation}
The faces commute with the unshifted internal differential, so the signs $(-1)^p$ and $(-1)^{p-1}$ cancel the mixed terms of $d_{\mathrm{tot}}^2$. After taking internal cohomology, \eqref{eq:b} is the differential $d_1$.

\subsection{Pseudoheight and cancellation}

Put $e(E,G)=\inf\RHom(E,G)$. A path $P=(i_0<\dots<i_p)$ is \emph{internally active} if each of its internal morphism complexes $\RHom(E_{i_s},E_{i_{s+1}})$ has nonzero cohomology. Define
\begin{equation}\label{eq:kappa}
  \kappa(P)=p-\sum_{s=0}^{p-1}e(E_{i_s},E_{i_{s+1}}),\qquad \tau_F(P)=e\bigl(E_{i_p},H_F(E_{i_0})\bigr),
\end{equation}
where the empty sum is zero; paths with acyclic terminal coefficient are omitted.

\begin{definition}\label{def:pseudoheight}
The \emph{$F$-pseudoheight}, generalizing \cite[Definition~4.4]{Kuz15}, is
\begin{equation}\label{eq:pseudoheight}
  \ph_F(\bE)=\min_P\{\tau_F(P)-\kappa(P)\},
\end{equation}
where the minimum runs over internally active paths with nonacyclic terminal coefficient, and is $+\infty$ if there is no such path.
\end{definition}

\begin{proposition}\label{prop:pseudoheight}
One has $h_F(\bE)\ge\ph_F(\bE)$. If $C_F\ne0$ and $\ph_F(\bE)>\alpha_F$, then $\bE$ is not full.
\end{proposition}

\begin{proof}
The lowest total degree of the summand of \eqref{eq:bar-ss} indexed by $P$ is $\sum_se(E_{i_s},E_{i_{s+1}})+\tau_F(P)-p=\tau_F(P)-\kappa(P)$. So the first page vanishes in total degrees below $\ph_F(\bE)$, and so does the abutment. Apply Proposition~\ref{prop:range}.
\end{proof}

More generally, if every $E_2^{a,b}$ with $a+b<t$ vanishes, then $h_F\ge t$, even if the first page is nonzero in those degrees; this form of the estimate is used in Section~\ref{sec:point-height}. Independent shifts $E_i\mapsto E_i[s_i]$ preserve the pseudoheight, since the total change along a path is $\sum_r(s_{i_r}-s_{i_{r+1}})+(s_{i_p}-s_{i_0})=0$.

\section{Fourier--Mukai kernels and exceptional approximation}\label{sec:kernels}

Let $X$ be smooth projective of dimension $d$. Derived Morita theory identifies perfect enhanced endofunctors of $\Perf(X)$ with kernels in $\Perf(X\times X)$ \cite[Theorem~8.15]{Toe07}; our convention is
\[
  F_K(E)=Rp_{2*}\bigl(Lp_1^*E\otimes^LK\bigr).
\]
We write $C_K,N_K,R_K$ for the coefficient complexes of $F_K$; the functor $F_K$ need not preserve the components of the decomposition.

\subsection{Kernel localization}

\begin{proposition}\label{prop:kernel-localization}
Let $Q^R$ and $P^R$ be the kernels of the right adjoints of $Q_\cB$ and $P_\cA$. Then
\[
  C_K\simeq\RHom_{X\times X}(\cO_\Delta,K),\qquad N_K\simeq\RHom_{X\times X}(Q^R,K),\qquad R_K\simeq\RHom_{X\times X}(P^R,K),
\]
and the localization triangle is obtained by applying $\RHom(-,K)$ to $P^R\to\cO_\Delta\to Q^R$.
\end{proposition}

\begin{proof}
The kernel triangle is obtained by taking right adjoints in \eqref{eq:projection}. Since $P_\cA^R=ii^!$, adjunction gives
\[
  \RHom(ii^!,F_K)\simeq\RHom(i,F_Ki)\simeq\RHom(\id_\cA,i^!F_Ki);
\]
explicitly, if $\eta\colon i\to F_Ki$ and $\epsilon\colon ii^!\to\id$ is the counit, the first adjunction sends $\eta$ to the composite $ii^!\xrightarrow{\eta i^!}F_Kii^!\xrightarrow{F_K\epsilon}F_K$. For $Q=Q_\cB$, double right adjunction and Serre duality yield $Q^{RR}\simeq S_XQS_X^{-1}$, so
\[
  \RHom(Q^R,F_K)\simeq\RHom(\id,Q^{RR}F_K)\simeq\Tr_X(QS_X^{-1}F_K).
\]
Derived Morita theory identifies these complexes of transformations with Hom complexes of kernels, compatibly with the maps.
\end{proof}

Convolution of kernels is defined by $L\star K=Rp_{13*}(Lp_{12}^*K\otimes^LLp_{23}^*L)$, so that $F_{L\star K}=F_LF_K$.

\subsection{External products}

For external kernels, cyclicity reduces localization to the projection triangle of an object.

\begin{theorem}\label{thm:external}
For $T,Y\in\Perf(X)$, put $K_{T,Y}=T^\vee\boxtimes S_X(Y)$. There are natural identifications
\begin{equation}\label{eq:external}
\begin{gathered}
  F_{K_{T,Y}}(E)\simeq\RHom_X(T,E)\otimes_kS_X(Y),\qquad C_{K_{T,Y}}\simeq\RHom_X(T,Y),\\
  N_{K_{T,Y}}\simeq\RHom_X(T,Q_\cB Y),\qquad R_{K_{T,Y}}\simeq\RHom_\cA(i^*T,i^*Y).
\end{gathered}
\end{equation}
Under these identifications the localization triangle is \eqref{eq:intro-external}, and its second arrow is the map on morphisms induced by $i^*$.
\end{theorem}

\begin{proof}
The formula for the functor follows from the projection formula. Put $U(V)=V\otimes_kY$ and $V_T(E)=\RHom_X(T,E)$. Cyclicity identifies the trace of $UV_T$ with that of $V_TU$ on $\Perf(k)$, namely $\RHom_X(T,Y)$. Applying this to $S_X^{-1}F_{K_{T,Y}}$ and to $Q_\cB S_X^{-1}F_{K_{T,Y}}$ gives the identities for $C$ and $N$. Now apply $\RHom_X(T,-)$ to the projection triangle $Q_\cB Y\to Y\to ii^*Y$: its last term is $\RHom_X(T,ii^*Y)\simeq\RHom_\cA(i^*T,i^*Y)$, and the adjunction unit sends a morphism $T\to Y$ to its image under $i^*$. All these constructions use the same projector triangle and evaluation maps, so they agree with the restriction map of Theorem~\ref{thm:localization}.
\end{proof}

For an exceptional collection, the terminal coefficient of the bar complex becomes
\begin{equation}\label{eq:external-terminal}
  M_{K_{T,Y}}(E_i,E_j)=\RHom_X(T,E_i)\otimes_k\RHom_X(E_j,Y),
\end{equation}
on which the first endpoint acts by postcomposition and the second by precomposition.

\begin{corollary}\label{cor:detect}
For $\Perf(X)=\langle\cA,\cB\rangle$, the following are equivalent:
\begin{enumerate}[(i)]
  \item $\cA\ne0$;
  \item some external-product coefficient has nonzero ambient complex and zero normal complex;
  \item some perfect kernel has positive excess height.
\end{enumerate}
\end{corollary}

\begin{proof}
Choose $0\ne T\in\cA=\cB^\perp$. Then $j^!T=0$, so $Q_\cB T=0$, and Theorem~\ref{thm:external} gives $N_{K_{T,T}}=0$ and $C_{K_{T,T}}=\RHom_X(T,T)\ne0$. This proves (i)$\Rightarrow$(ii); the normal height in (ii) is $+\infty$, which gives (iii); and (iii) implies (i) by Proposition~\ref{prop:range}.
\end{proof}

To apply the corollary one needs a coefficient whose normal complex can be computed from the collection alone, without knowing $\cA$; Section~\ref{sec:special} provides such computations.

\section{Simultaneous functorial heights}\label{sec:simultaneous}

Several coefficients carry more information than each of them separately. Mixed external coefficients retain the morphisms of one residual dg category, and morphisms between kernels give cones with stronger normal vanishing. We develop both constructions in the enhanced setting of Sections~\ref{sec:localization} and~\ref{sec:kernels}, and end with a constraint on hereditary orthogonal components.

\subsection{A matrix of coefficients and one residual functor}

Fix a semiorthogonal decomposition $\Db(X)=\langle\cA,\cB\rangle$, with $Q=jj^!$, $P=ii^*$ and $Q\to\id\to P$. Choose perfect objects $T_1,\dots,T_m$, put $W_a=i^*T_a$, and consider all the ordered external coefficients $K_{ab}=T_a^\vee\boxtimes S_X(T_b)$. Write
\[
  N_{ab}=\RHom_X(T_a,QT_b),\qquad C_{ab}=\RHom_X(T_a,T_b),\qquad R_{ab}=\RHom_\cA(W_a,W_b)
\]
for their normal, ambient and residual complexes, with triangles
\begin{equation}\label{eq:matrix-triangle}
  N_{ab}\longrightarrow C_{ab}\longrightarrow R_{ab},
\end{equation}
and put $h_{ab}=\inf N_{ab}$.

\begin{definition}\label{def:profile}
The \emph{multiple-height profile} of the ordered probe system $\bT=(T_1,\dots,T_m)$ relative to $\cB\subset\Db(X)$ is the matrix $h_\cB(\bT)=(h_{ab})_{1\le a,b\le m}$. The \emph{simultaneous coefficient data} consist of the triangles \eqref{eq:matrix-triangle} together with the composition maps induced by $i^*$.
\end{definition}

For shifted probes $\widetilde T_a=T_a[t_a]$ one has $\widetilde N_{ab}\simeq N_{ab}[t_b-t_a]$ and $\widetilde h_{ab}=h_{ab}+t_a-t_b$, so the diagonal heights and the sums of heights around directed cycles are invariant.

\begin{theorem}[Simultaneous coefficient comparison]\label{thm:simultaneous}
The maps $C_{ab}\to R_{ab}$ are the morphism maps of the single enhanced functor $i^*$; hence they are compatible with identities and with all Yoneda compositions. If $h_{ab}\ge d_{ab}$, then
\[
  \Ext^q_X(T_a,T_b)\xrightarrow{\ \sim\ }\Ext^q_\cA(W_a,W_b)\quad(q\le d_{ab}-2),\qquad
  \Ext^{d_{ab}-1}_X(T_a,T_b)\hookrightarrow\Ext^{d_{ab}-1}_\cA(W_a,W_b).
\]
For $T=\bigoplus_aT_a$ and $W=i^*T$, these comparisons assemble into a map of dg algebras
\[
  E=\RHom_X(T,T)\longrightarrow E_\cA=\RHom_\cA(W,W),
\]
whose homotopy fibre, as an $E$-bimodule, is $N=\RHom_X(T,QT)=\bigoplus_{a,b}N_{ab}$.
\end{theorem}

\begin{proof}
Under the adjunction $\RHom_X(T_a,ii^*T_b)\simeq\RHom_\cA(i^*T_a,i^*T_b)$, the map induced by $T_b\to ii^*T_b$ is the morphism map of $i^*$. Applying $\RHom_X(T_a,-)$ to the projection triangle gives \eqref{eq:matrix-triangle}, and naturality gives compatibility with compositions. The cohomology sequence gives the isomorphisms and injections. Taking direct sums gives the algebra map and its fibre, where $i^*$ is represented by a dg functor after replacing the enhancements, or by an $A_\infty$ functor. The left $E$-action on $\RHom_X(T,QT)$ applies $Q$ to endomorphisms of $T$, the right action is precomposition, and naturality of the counit makes the identification of the fibre $E$-bilinear.
\end{proof}

\begin{corollary}\label{cor:reconstruct}
If $G=\bigoplus_aT_a$ is a classical generator of $\Db(X)$, then $W=i^*G$ is a classical generator of $\cA$, and $\cA\simeq\Perf(\RHom_\cA(W,W))$. In particular, $\cB=\Db(X)$ if and only if every comparison $N_{ab}\to C_{ab}$ is a quasi-isomorphism.
\end{corollary}

\begin{proof}
For $A\in\cA$, applying $i^*$ to $iA\in\thick(G)$ gives $A\simeq i^*iA\in\thick(W)$, and derived Morita theory \cite[Theorem~3.8(b)]{Kel06} gives the equivalence. The comparisons are all quasi-isomorphisms exactly when $\RHom_\cA(W,W)=0$, that is, when $W=0$ and hence $\cA=0$.
\end{proof}

\begin{remark}\label{rem:dg}
For a smooth proper decomposition $\cD=\langle\cA,\cB\rangle$, the perfect endofunctors $F_{ab}(E)=\RHom_\cD(T_a,E)\otimes_kS_\cD(T_b)$ give the same three complexes, by the trace computation in the proof of Theorem~\ref{thm:external}. So the statements of this subsection hold for arbitrary smooth proper dg categories.
\end{remark}

\begin{remark}\label{rem:addition}
Adding coefficients does not improve the height: by exactness $N_{\bigoplus_\nu K_\nu}=\bigoplus_\nu N_{K_\nu}$ and $h_{\bigoplus_\nu K_\nu}=\min_\nu h_{K_\nu}$. The additional information lies in the mixed entries and in the morphisms between coefficients, which is what the rest of this section exploits.
\end{remark}

\begin{corollary}[Preservation of a Yoneda product]\label{cor:yoneda}
Let $u_\nu\in\Ext^{p_\nu}_X(T_{a_{\nu-1}},T_{a_\nu})$ for $1\le\nu\le\ell$, and suppose that the composite $u_\ell\cdots u_1\in\Ext^p_X(T_{a_0},T_{a_\ell})$, $p=\sum_\nu p_\nu$, is nonzero. If $h_{a_0a_\ell}\ge p+1$, then the composite of the images in $\cA$ is nonzero; in particular none of the images $i^*u_\nu$ is zero.
\end{corollary}

\begin{proof}
The functor $i^*$ preserves the product, and the comparison is injective in degree $p$ under the stated bound.
\end{proof}

\subsection{The dg category of coefficient kernels}

For the fixed projection triangle, the assignments
\[
  K\mapsto N_K=\RHom(Q^R,K),\qquad K\mapsto C_K=\RHom(\cO_\Delta,K),\qquad K\mapsto R_K=\RHom(P^R,K)
\]
are exact enhanced functors on $\Perf(X\times X)$ forming a natural triangle $N\to C\to R$, and morphisms of kernels act on this triangle.

\begin{proposition}[A coefficient dg algebra]\label{prop:coefficient-algebra}
Let $K_1,\dots,K_m\in\Perf(X\times X)$, put $K_\Sigma=\bigoplus_aK_a$, and set
\[
  \cH=\RHom(K_\Sigma,K_\Sigma),\qquad M_N=N_{K_\Sigma},\qquad M_C=C_{K_\Sigma},\qquad M_R=R_{K_\Sigma}.
\]
Postcomposition makes these three complexes left dg $\cH$-modules, and their comparison is a triangle of $\cH$-modules. If $V$ is a perfect right $\cH$-module and $K_V=V\otimes^L_\cH K_\Sigma$, then the coefficient triangle of $K_V$ is naturally
\begin{equation}\label{eq:module-triangle}
  V\otimes^L_\cH M_N\longrightarrow V\otimes^L_\cH M_C\longrightarrow V\otimes^L_\cH M_R .
\end{equation}
Every kernel in the thick subcategory generated by the $K_a$ is of this form. In particular
\[
  h_{K_V}=\inf\bigl(V\otimes^L_\cH M_N\bigr),\qquad\alpha_{K_V}=\inf\bigl(V\otimes^L_\cH M_C\bigr),
\]
so that on this subcategory the computation of heights is a computation with dg modules over $\cH$.
\end{proposition}

\begin{proof}
Choose enhanced representatives. Postcomposition gives the left $\cH$-actions, and naturality makes the comparison $\cH$-linear. Formula \eqref{eq:module-triangle} holds for $V=\cH$ and is preserved by shifts, finite sums, cones and retracts, hence holds for every perfect right module. These operations also build $\thick(K_\Sigma)$, and derived Morita theory \cite[Theorem~3.8(b)]{Kel06} shows that every object of $\thick(K_\Sigma)$ is of the form $K_V$.
\end{proof}

Morphisms of coefficients can cancel low normal cohomology without cancelling the corresponding ambient cohomology.

\begin{proposition}[Cancellation by a morphism of kernels]\label{prop:cancellation}
Let $f\colon K_1\to K_2$ be a morphism of kernels and put $K_f=\Cone(f)$. For $M=N$ or $M=C$, write $f_M^q\colon H^q(M_{K_1})\to H^q(M_{K_2})$. For every $q$ there is a natural short exact sequence
\begin{equation}\label{eq:cone-sequence}
  0\longrightarrow\coker f^q_M\longrightarrow H^q(M_{K_f})\longrightarrow\ker f^{q+1}_M\longrightarrow0 .
\end{equation}
In particular, if $f^q_N$ is an isomorphism for every $q<d$ and injective for $q=d$, then $h_{K_f}\ge d$. If in addition either
\begin{enumerate}[(i)]
  \item $f^r_C$ is not surjective for some $r<d$, or
  \item $f^r_C$ is not injective for some $r\le d$,
\end{enumerate}
then $C_{K_f}\ne0$, $\alpha_{K_f}<d\le h_{K_f}$, and $\cB$ is not full.
\end{proposition}

\begin{proof}
Apply the long exact cohomology sequence to $M_{K_1}\to M_{K_2}\to M_{K_f}$. The hypotheses on $f_N$ make both outer terms of \eqref{eq:cone-sequence} vanish for $q<d$. A cokernel of $f_C$ in degree $r$ contributes to $H^r(C_{K_f})$, and a kernel in degree $r$ contributes to $H^{r-1}(C_{K_f})$; the two alternatives place this cohomology below $d$.
\end{proof}

\begin{remark}\label{rem:cone-shift}
For $N_{K_1},N_{K_2}\in D^{\ge a}$, the cone lies a priori only in $D^{\ge a-1}$. To improve this to $D^{\ge a+1}$ one needs an isomorphism on $H^a$ of the normal complexes and injectivity on $H^{a+1}$; separate height bounds for $K_1$ and $K_2$ are not enough.
\end{remark}

\begin{proposition}[A finite complex of coefficients]\label{prop:finite-complex}
Let $K$ be the totalization of a finite complex, or a finite twisted complex, of kernels $K^p$, with the finite filtration satisfying $\gr^pK\simeq K^p[-p]$. There are compatible finite spectral sequences
\begin{equation}\label{eq:complex-ss}
  E_1^{p,q}(N)=H^q(N_{K^p})\Longrightarrow H^{p+q}(N_K),\qquad E_1^{p,q}(C)=H^q(C_{K^p})\Longrightarrow H^{p+q}(C_K),
\end{equation}
whose first differential is induced by the part of the differential of $K$ that raises the filtration degree by one. In particular,
\begin{equation}\label{eq:complex-bound}
  h_K\ge\min_p\{p+h_{K^p}\}.
\end{equation}
More strongly, if some page $E_r(N)$ vanishes in all total degrees below $d$, and the ambient abutment is nonzero in one of those degrees, then $\alpha_K<d\le h_K$ and $\cB$ is not full.
\end{proposition}

\begin{proof}
The exact functors $N$ and $C$ carry the graded pieces to $N_{K^p}[-p]$ and $C_{K^p}[-p]$, which gives the spectral sequences of the filtered complexes \eqref{eq:complex-ss}; finiteness ensures convergence \cite[Section~5.5]{Wei94}. The first page gives \eqref{eq:complex-bound}, and vanishing on a later page persists to the abutment.
\end{proof}

\begin{example}[Cancellation on projective three-space]\label{ex:P3}
Let $X=\PP^3$, $\cB=\langle\cO_X\rangle$ and $T=\Omega^1_X(1)$. The Euler sequence is
\begin{equation}\label{eq:euler}
  0\longrightarrow T\longrightarrow\cO_X^{\oplus4}\xrightarrow{\ \varepsilon\ }\cO_X(1)\longrightarrow0 .
\end{equation}
Put $Y_1=\cO_X^{\oplus4}$, $Y_2=\cO_X(1)$, and choose the external coefficients $K_j=T^\vee\boxtimes S_X(Y_j)$ with the morphism $f=1\boxtimes S_X(\varepsilon)\colon K_1\to K_2$. Since $Q_\cB Y=R\Gamma(Y)\otimes\cO_X$, one finds
\begin{equation}\label{eq:P3-dims}
  N_{K_1}\simeq\C^{16},\qquad N_{K_2}\simeq\C^{16},\qquad C_{K_1}\simeq\C^{16},\qquad C_{K_2}\simeq\C^{15},
\end{equation}
all in degree zero: $T^\vee=T_X(-1)$ has four sections, $T^\vee(1)=T_X$ has fifteen, and neither has higher cohomology. Thus $h_{K_j}=\alpha_{K_j}=0$ for $j=1,2$, and neither coefficient detects anything by itself. But $R\Gamma(\varepsilon)$ is an isomorphism, so the normal map $f_N$ is an isomorphism, while the ambient map $f_C$ is $\mathfrak{gl}_4\to\mathfrak{pgl}_4$, with one-dimensional kernel. Since $\Cone(\varepsilon)=T[1]$,
\[
  N_{K_f}=0,\qquad C_{K_f}\simeq\RHom_X(T,T)[1]\simeq\C[1].
\]
After the shift by $[-1]$ the kernel is $T^\vee\boxtimes S_X(T)$, with
\begin{equation}\label{eq:P3-result}
  \alpha_{K_f[-1]}=0,\qquad h_{K_f[-1]}=+\infty .
\end{equation}
So the morphism of coefficients detects that $\langle\cO_X\rangle$ is not full, by producing the exceptional bundle $T$ of Beilinson's collection \cite{Bei78} in its orthogonal.
\end{example}

\begin{remark}[Convolution is not additive for height]\label{rem:convolution}
The coefficient functors are exact but not monoidal, and the inequality $h_{L\star K}\ge h_L+h_K$ fails in general. Take $X=\PP^1$ and $\cB=\Db(X)$, so that $N=C$, and put
\[
  K_0=\cO_X\boxtimes\cO_X(-3)[1],\qquad K=\cO_\Delta\oplus K_0,\qquad L=\cO_X(3)\boxtimes\cO_X(-2)[1].
\]
Then $C_{K_0}=R\Gamma(\cO_X(-1))=0$, so $h_K=0$, and $h_L=0$. The convolution $L\star K_0$ is $\cO_X\boxtimes\cO_X(-2)[2]$, whose ambient complex is $\C[1]$. Hence $h_{L\star K}=-1<0=h_L+h_K$.
\end{remark}

\subsection{Hereditary orthogonal components}\label{subsec:hereditary}

Mixed coefficients also constrain the orthogonal itself. We record one such constraint, for orthogonals which are derived categories of hereditary abelian categories, for instance of representations of a quiver; it is applied to point coefficients on surfaces in Corollary~\ref{cor:hereditary-surface}. The proof uses the Euler form.

\begin{theorem}[Two numerically equivalent probes]\label{thm:two-probes}
Let $\Db(X)=\langle\cA,\cB\rangle$ and let $T_1,T_2\in\Perf(X)$ satisfy
\[
  \RHom_X(T_1,T_2)=\RHom_X(T_2,T_1)=0,\qquad\End_X(T_a)=\C,\qquad\dim\Ext^1_X(T_a,T_a)\ge2,
\]
and have the same left Euler functional:
\begin{equation}\label{eq:same-euler}
  \chi_X(T_1,E)=\chi_X(T_2,E)\qquad\text{for every }E\in\Perf(X).
\end{equation}
Suppose that the four coefficients $K_{ab}=T_a^\vee\boxtimes S_X(T_b)$ have heights $h_{11},h_{22},h_{12},h_{21}\ge2$, and put $W_a=i^*T_a$. If $\cA$ is equivalent to the bounded derived category of an Ext-finite hereditary abelian category, then
\begin{equation}\label{eq:quantitative}
\begin{aligned}
  \dim H^2N_{K_{12}}&=\dim\Ext^1_\cA(W_2,W_2)-1\ge\dim\Ext^1_X(T_2,T_2)-1,\\
  \dim H^2N_{K_{21}}&=\dim\Ext^1_\cA(W_1,W_1)-1\ge\dim\Ext^1_X(T_1,T_1)-1 .
\end{aligned}
\end{equation}
In particular, if $h_{12}\ge3$ or $h_{21}\ge3$, then $\cA$ is not of this form.
\end{theorem}

\begin{proof}
Suppose $\cA\simeq\Db(\cH')$ with $\cH'$ hereditary. The diagonal comparison of Theorem~\ref{thm:simultaneous} gives $\End_\cA(W_a)=\C$ and $\Ext^1_X(T_a,T_a)\hookrightarrow\Ext^1_\cA(W_a,W_a)$. In the derived category of a hereditary category every bounded complex is the direct sum of its shifted cohomology objects \cite[Chapter~I]{Hap88}, so the indecomposable object $W_a$ is of the form $N_a[s_a]$ with $N_a\in\cH'$, and
\[
  \chi_\cA(W_a,W_a)=1-\dim\Ext^1_\cA(W_a,W_a)<0 .
\]
Adjunction and \eqref{eq:same-euler} give $\chi_\cA(W_1,Z)=\chi_X(T_1,iZ)=\chi_X(T_2,iZ)=\chi_\cA(W_2,Z)$ for $Z\in\cA$, hence
\begin{equation}\label{eq:mixed-negative}
  \chi_\cA(W_1,W_2)=\chi_\cA(W_2,W_2)<0,\qquad\chi_\cA(W_2,W_1)=\chi_\cA(W_1,W_1)<0 .
\end{equation}
For $a\ne b$ the ambient complex vanishes, so $\RHom_\cA(W_a,W_b)\simeq N_{K_{ab}}[1]\in D^{\ge1}(\C)$, while $\Ext^q_\cA(W_a,W_b)=\Ext^{q+s_b-s_a}_{\cH'}(N_a,N_b)$ can be nonzero only for $q=s_a-s_b$ and $q=s_a-s_b+1$. If $s_1<s_2$, these degrees are at most zero for $(a,b)=(1,2)$, so $\RHom_\cA(W_1,W_2)=0$, contradicting \eqref{eq:mixed-negative}; by symmetry $s_2<s_1$ is excluded as well. So $s_1=s_2$, and both mixed complexes are concentrated in degree one. Therefore
\[
  -\dim H^2N_{K_{12}}=\chi_\cA(W_1,W_2)=\chi_\cA(W_2,W_2)=1-\dim\Ext^1_\cA(W_2,W_2),
\]
and similarly with the indices exchanged; the inequalities follow from the diagonal injections. Since $\dim\Ext^1_X(T_a,T_a)\ge2$, both spaces $H^2N_{K_{12}}$ and $H^2N_{K_{21}}$ are nonzero, which gives the last assertion.
\end{proof}

\begin{remark}[Point probes]\label{rem:point-probes}
Skyscraper sheaves of distinct points have disjoint supports; the Koszul resolution gives scalar endomorphisms and $\Ext^1_X(\cO_x,\cO_x)=T_xX$, and local duality gives $\chi_X(\cO_x,E)=(-1)^n\rk(E)$, independently of $x$. So point probes satisfy the hypotheses of Theorem~\ref{thm:two-probes} on every smooth projective variety of dimension $n\ge2$.
\end{remark}

\begin{corollary}[Two point coefficients]\label{cor:two-points}
Let $X$ be a connected smooth projective $n$-fold over $\C$ with $n\ge2$, and let $x\ne y$ be closed points whose four kernels $K_{uv}=\cO_u^\vee\boxtimes S_X(\cO_v)$, $u,v\in\{x,y\}$, have height at least two. If $\cA$ is the bounded derived category of an Ext-finite hereditary abelian category, then
\[
  \dim H^2N_{K_{xy}}=\dim\Ext^1_\cA(W_y,W_y)-1\ge n-1,
\]
and similarly with $x$ and $y$ exchanged.
\end{corollary}

\begin{proof}
Apply Theorem~\ref{thm:two-probes} to the point probes of Remark~\ref{rem:point-probes}.
\end{proof}

\section{Special coefficients and height calculations}\label{sec:special}

Let $X$ be a smooth connected projective variety of dimension $d$ over $k$. We collect the coefficient calculations used later. The ambient complexes come from standard duality and Hochschild--Kostant--Rosenberg formulas, and their terminal factors give explicit instances of the pseudoheight bound.

\subsection{Tensor coefficients and Serre powers}

Assume in this subsection that $\operatorname{char}k=0$. For $M\in\Pic(X)$, write $T_M=(-)\otimes M$, with kernel $\Delta_*M$. Adjunction and the HKR isomorphism \cite[Theorem~4.1 and Corollary~4.2]{Cal05} give $C_{T_M}(X)\simeq\RHom_X(\bigoplus_p\Omega^p_X[p],M)$, hence
\begin{equation}\label{eq:tensor-ambient}
  H^tC_{T_M}(X)\simeq\bigoplus_{p+q=t}H^q\bigl(X,\textstyle\bigwedge^pT_X\otimes M\bigr).
\end{equation}
For line bundles $E_i=\cO_X(D_i)$, the terminal groups are
\begin{equation}\label{eq:tensor-terminal}
  \Ext^k\bigl(E_j,S_X^{-1}T_M(E_i)\bigr)=H^{k-d}\bigl(X,\cO_X(D_i-D_j)\otimes\omega_X^{-1}\otimes M\bigr),
\end{equation}
since $S_X^{-1}=(-)\otimes\omega_X^{-1}[-d]$.

On a surface, put $F_{2K}=(-)\otimes\omega_X^2$. For coherent sheaves $E_i,E_j$, let $u(E_i,E_j)=\sup\{m:\Ext^m(E_i,E_j)\ne0\}$, with $u=-\infty$ when the morphism complex is acyclic. Serre duality yields
\begin{equation}\label{eq:2K-terminal}
  \Ext^k\bigl(E_j,S_X^{-1}F_{2K}(E_i)\bigr)=\Ext^{k-2}(E_j,E_i\otimes\omega_X)\simeq\Ext^{4-k}(E_i,E_j)^\vee .
\end{equation}
For an internally active path $P=(i_0<\dots<i_p)$ with $\RHom(E_{i_0},E_{i_p})\ne0$, define $\Delta_2(P)=u(E_{i_0},E_{i_p})+\kappa(P)$ and $\Delta_2(\bE)=\max_P\Delta_2(P)$.

\begin{proposition}\label{prop:2K}
For a nonempty exceptional collection of coherent sheaves on a surface, $\ph_{F_{2K}}(\bE)=4-\Delta_2(\bE)$. If $X$ is geometrically connected, minimal of general type, with $p_g=q=0$, then
\[
  \alpha_{F_{2K}}=0,\qquad H^0C_{F_{2K}}=H^0(X,2K_X),\qquad h^0(X,2K_X)=1+K_X^2 .
\]
Restriction on $H^0(X,2K_X)$ is injective if $\Delta_2(\bE)\le3$ and an isomorphism if $\Delta_2(\bE)\le2$; in either case the orthogonal is nonzero.
\end{proposition}

\begin{proof}
Equation \eqref{eq:2K-terminal} gives terminal degree $4-u(E_{i_0},E_{i_p})$, which proves the formula for the pseudoheight. By \eqref{eq:tensor-ambient} the ambient complex is concentrated in nonnegative degrees, with $H^0(X,2K_X)$ in degree zero. Since $K_X$ is nef and big, Kawamata--Viehweg vanishing \cite{Kaw82} and Riemann--Roch give
\[
  h^0(X,2K_X)=\chi(\cO_X)+\tfrac12(2K_X)(2K_X-K_X)=1+K_X^2>0 .
\]
The two bounds give height at least one and two, respectively; apply Proposition~\ref{prop:range}.
\end{proof}

For coherent sheaves on a surface $u(E_i,E_j)\le2$, so vanishing of all forward Hom groups gives $\Delta_2\le2$; more generally $\Delta_2\le3$ if each internally active path has at most one edge with nonzero Hom in degree zero.

In arbitrary dimension and for $a\in\Z$,
\[
  C_{S_X^a}=\RHom(\id,S_X^a),\qquad N_{S_X^a}=\Tr_X(Q_\cB S_X^{a-1}),
\]
and since $S_X^a=(-)\otimes\omega_X^a[ad]$, equation \eqref{eq:tensor-ambient} gives
\begin{equation}\label{eq:serre-ambient}
  H^tC_{S_X^a}\simeq\bigoplus_{p+q=t+ad}H^q\bigl(X,\textstyle\bigwedge^pT_X\otimes\omega_X^a\bigr).
\end{equation}
On a surface $F_{2K}=S_X^2[-4]$, so these two coefficients have the same excess height. Although $i^!S_Xi=S_\cA$, the functor $i^!S_X^ai$ is in general different from $S_\cA^a$; Proposition~\ref{prop:compression} describes the difference for positive powers.

\subsection{The point coefficient}

For a closed $k$-rational point $x\in X$, put
\[
  P_x(E)=E\otimes^L\cO_x,\qquad K_x=\Delta_*\cO_x,\qquad\tau_x=\det T_xX=\omega_X^{-1}|_x .
\]

\begin{proposition}\label{prop:point}
There are natural identifications
\begin{equation}\label{eq:point-ambient}
  K_x\simeq K_{\cO_x,\cO_x},\qquad C_{P_x}(X)\simeq\RHom_X(\cO_x,\cO_x),\qquad H^tC_{P_x}(X)\simeq\textstyle\bigwedge^tT_xX .
\end{equation}
Thus $\alpha_{P_x}=0$. For vector bundles $E_i,E_j$,
\begin{equation}\label{eq:point-terminal}
  \RHom_X\bigl(E_j,S_X^{-1}P_x(E_i)\bigr)\simeq(E_j|_x)^\vee\otimes_kE_i|_x\otimes_k\tau_x[-d].
\end{equation}
\end{proposition}

\begin{proof}
The Koszul resolution at $x$ gives $\cO_x^\vee\simeq\cO_x\otimes_k\tau_x[-d]$ and $S_X(\cO_x)\simeq\cO_x\otimes_k\tau_x^\vee[d]$, whose external product is $K_x$. Theorem~\ref{thm:external} identifies its ambient complex with the self-Ext complex of $\cO_x$; applying $\Hom(-,\cO_x)$ to the Koszul resolution gives zero differential and cohomology $\bigwedge^tT_xX$. Finally $S_X^{-1}P_x(E_i)=E_i|_x\otimes_k\tau_x\otimes_k\cO_x[-d]$ and $\RHom_X(E_j,\cO_x)=(E_j|_x)^\vee$.
\end{proof}

\begin{proposition}\label{prop:point-pseudoheight}
For a nonempty exceptional collection of vector bundles,
\begin{equation}\label{eq:point-pseudoheight}
  \ph_{P_x}(\bE)=d-\max_{P\text{ internally active}}\kappa(P).
\end{equation}
In particular, the collection is not full if every internally active path satisfies $\kappa(P)\le d-1$.
\end{proposition}

\begin{proof}
Each terminal factor in \eqref{eq:point-terminal} is nonzero exactly in degree $d$. Apply \eqref{eq:pseudoheight} and $\alpha_{P_x}=0$.
\end{proof}

On a surface, vanishing of all forward Hom spaces gives $\kappa(P)\le0$, with equality on paths of length zero, hence pseudoheight two; allowing at most one edge with nonzero Hom on each path gives pseudoheight at least one.

\begin{lemma}\label{lem:no-hom}
Let $\bE=(L_1,\dots,L_N)$ be a nonempty exceptional collection of line bundles on a smooth projective surface with $\Hom(L_i,L_j)=0$ for $i<j$. For every closed $k$-rational point $x$,
\[
  \alpha_{P_x}=0,\qquad\ph_{P_x}(\bE)=h_{P_x}(\bE)=2 .
\]
\end{lemma}

\begin{proof}
Proposition~\ref{prop:point-pseudoheight} gives pseudoheight two. Moreover $E_1^{0,2}(P_x)=\bigoplus_{i=1}^N\tau_x\ne0$. All differentials into this term have sources of total degree one, which vanish on $E_1$, and all differentials out of it have targets with positive first index, which are zero. So this term survives to $E_\infty$.
\end{proof}

\subsection{Curve coefficients}

For $A,B\in\Perf(X)$, projection and duality along the diagonal give
\begin{equation}\label{eq:box-ambient}
  F_{A\boxtimes B}(E)=R\Gamma(X,E\otimes^LA)\otimes_kB,\qquad C_{A\boxtimes B}(X)\simeq R\Gamma\bigl(X,A\otimes^LB\otimes\omega_X^{-1}\bigr)[-d],
\end{equation}
since the diagonal is a regular immersion with normal bundle $T_X$, so that $\Delta^!K=L\Delta^*K\otimes\omega_X^{-1}[-d]$. For vector bundles $E_i,E_j$ the terminal factor is
\begin{equation}\label{eq:box-terminal}
  R\Gamma(X,E_i\otimes^LA)\otimes_kR\Gamma\bigl(X,E_j^\vee\otimes^LB\otimes\omega_X^{-1}\bigr)[-d].
\end{equation}
On a surface, put $a_i=\inf R\Gamma(X,E_i\otimes^LA)$ and $b_j=\inf R\Gamma(X,E_j^\vee\otimes^LB\otimes\omega_X^{-1})$. Then
\begin{equation}\label{eq:box-pseudoheight}
  \ph_{A\boxtimes B}(\bE)=\min_{P=(i_0<\dots<i_p)}\{2+a_{i_0}+b_{i_p}-\kappa(P)\},
\end{equation}
the minimum running over internally active paths with nonacyclic endpoints.

Let $\iota_\nu\colon C_\nu\hookrightarrow X$ be effective Cartier divisors on a surface, $M_\nu\in\Pic(C_\nu)$, and $K=\iota_{1*}M_1\boxtimes\iota_{2*}M_2$. Its endpoint complexes are
\[
  R\Gamma(C_1,E_i|_{C_1}\otimes M_1),\qquad R\Gamma\bigl(C_2,E_j^\vee|_{C_2}\otimes M_2\otimes\omega_X^{-1}|_{C_2}\bigr),
\]
whose initial degrees, when finite, are zero or one. If the curves have no common component and $Z=C_1\cap C_2\ne\emptyset$, their local equations form a regular sequence, and
\begin{equation}\label{eq:two-curves}
  C_K\simeq H^0\bigl(Z,M_1|_Z\otimes M_2|_Z\otimes\omega_X^{-1}|_Z\bigr)[-2].
\end{equation}
Thus $\alpha_K=2$, and $\dim_kH^2C_K$ is the degree of the intersection scheme. For a closed $k$-rational point $x\in C$ and $K=\cO_x\boxtimes\iota_*M$, the Cartier resolution has, after tensoring with $\cO_x$, zero differential, and
\[
  H^1C_K=M_x\otimes\omega_C^{-1}|_x,\qquad H^2C_K=M_x\otimes\omega_X^{-1}|_x,
\]
with all other groups zero; if $x\notin C$ the ambient complex vanishes.

\begin{proposition}\label{prop:curve}
Let $\iota\colon C\hookrightarrow X$ be a nonempty effective Cartier divisor on a smooth projective surface. Put $\nu_C=\cO_C(C)$, choose $M\in\Pic(C)$, and set
\[
  K_{C,M}=\iota_*M\boxtimes\iota_*(\omega_C\otimes M^{-1}),\qquad V=\iota_*(\nu_C\otimes M^{-1}).
\]
Then
\[
  K_{C,M}\simeq K_{V,V}[-1],\qquad C_{K_{C,M}}\simeq\RHom_X(V,V)[-1],\qquad N_{K_{C,M}}\simeq\RHom_X(V,Q_\cB V)[-1].
\]
In particular $\alpha_{K_{C,M}}=1$, and the identity of $V$ defines a nonzero ambient class of degree one, whose restriction is $\id_{i^*V}$.
\end{proposition}

\begin{proof}
Duality for the Cartier divisor and adjunction give $V^\vee\simeq\iota_*M[-1]$ and $S_X(V)\simeq\iota_*(\omega_C\otimes M^{-1})[2]$. Their external product is $K_{C,M}[1]$, so Theorem~\ref{thm:external} applies. Since $V$ is a nonzero sheaf, its self-Ext complex begins in degree zero with $\Hom_X(V,V)=H^0(C,\cO_C)$; the shift places its identity in degree one.
\end{proof}

The local self-Ext sheaves of $V$ are $\cO_C$ and $\nu_C$, in degrees zero and one, so the local-to-global spectral sequence gives
\begin{gather*}
  H^1C_{K_{C,M}}=H^0(C,\cO_C),\qquad H^3C_{K_{C,M}}=H^1(C,\nu_C),\\
  0\to H^1(C,\cO_C)\to H^2C_{K_{C,M}}\to H^0(C,\nu_C)\to0 .
\end{gather*}
Write the endpoint complexes as $V_i(M)=R\Gamma(C,E_i|_C\otimes M)$ and $W_j(M)=R\Gamma(C,E_j^\vee|_C\otimes\nu_C\otimes M^{-1})$. If every internally active path has an acyclic endpoint or satisfies
\begin{equation}\label{eq:curve-criterion}
  \inf V_{i_0}(M)+\inf W_{i_p}(M)\ge\kappa(P),
\end{equation}
then \eqref{eq:box-pseudoheight} gives $h_{K_{C,M}}\ge2>1=\alpha_{K_{C,M}}$, and the collection is not full. For $\kappa(P)=1$ this requires one endpoint to have vanishing $H^0$; for $\kappa(P)=2$ both; for larger values an endpoint must be acyclic. For two properly intersecting curves the corresponding condition is $a_{i_0}+b_{i_p}\ge\kappa(P)+1$, because the ambient degree is two.

\subsection{Graph coefficients and symmetries}

Return to arbitrary dimension. For an automorphism $g$ of $X$ and $M\in\Pic(X)$, put
\[
  \gamma_g(x)=(x,g(x)),\qquad K_{g,M}=\gamma_{g*}(g^*M),
\]
with transform $F_{g,M}(E)=g_*E\otimes M$. The next calculation is the form with line-bundle coefficients of the fixed-locus HKR formula of Arinkin, C\u{a}ld\u{a}raru and Hablicsek \cite{ACH19}.

\begin{proposition}\label{prop:graph}
Suppose $g$ has finite order and $\operatorname{char}k=0$. For a connected component $Z$ of the fixed locus $X^g$, put $c_Z=\codim_XZ$. Then
\begin{equation}\label{eq:graph-ambient}
  H^tC_{K_{g,M}}\simeq\bigoplus_{Z\subset X^g}\ \bigoplus_{p+q+c_Z=t}H^q\bigl(Z,\textstyle\bigwedge^pT_Z\otimes\det N_{Z/X}\otimes M|_Z\bigr).
\end{equation}
\end{proposition}

\begin{proof}
The fixed locus is smooth, and its derived formality \cite[Corollary~1.12 and Theorem~4.4]{ACH19} decomposes the derived restriction to the diagonal into the summands $\Omega^a_Z[a]\otimes M|_Z$. To compute $\Delta^!K_{g,M}$, tensor with $\omega_X^{-1}|_Z[-d]$. With $z=\dim Z=d-c_Z$ and $p=z-a$, adjunction gives $\Omega^a_Z\otimes\omega_X^{-1}|_Z\simeq\bigwedge^pT_Z\otimes\det N_{Z/X}$, and the shift $[a-d]=[-p-c_Z]$ gives the stated total degree.
\end{proof}

For vector bundles, let $\lambda_{g,M}(i,j)=\inf\{q:\Ext^q(E_j,g_*E_i\otimes M\otimes\omega_X^{-1})\ne0\}$. The terminal calculation gives
\[
  \ph_{g,M}(\bE)=\min_P\{d+\lambda_{g,M}(i_0,i_p)-\kappa(P)\},
\]
and a nonzero fixed-locus class of degree $t$ survives restriction if $h_{g,M}>t$. On a surface, an isolated fixed point contributes in degree two, and a fixed curve contributes in degree one if $H^0(Z,N_{Z/X}\otimes M|_Z)\ne0$.

The following invariance is a consequence of the rigidity of admissible subcategories under connected groups of autoequivalences \cite[Corollary~3.15]{KO18}.

\begin{proposition}\label{prop:conjugation}
Let $X$ be smooth projective over an algebraically closed field of characteristic zero. If $\Phi$ lies in the identity component of the group of enhanced derived autoequivalences and $\bE$ is an exceptional collection, then
\[
  h_{\Phi F\Phi^{-1}}(\bE)=h_F(\bE),\qquad\ph_{\Phi F\Phi^{-1}}(\bE)=\ph_F(\bE),
\]
and the initial ambient degrees and excess heights agree as well. In particular $h_{P_{g(x)}}(\bE)=h_{P_x}(\bE)$ for $g\in\Aut^0(X)$.
\end{proposition}

\begin{proof}
The identity component consists of automorphisms in $\Aut^0(X)$ composed with tensor products by elements of $\Pic^0(X)$ \cite[Theorem~2.12]{Ros09}. Each admissible subcategory $\langle E_i\rangle$ is preserved \cite[Corollary~3.15]{KO18}, so $\Phi(E_i)\simeq E_i[n_i]$ for integers $n_i$. Conjugation transports the localization triangle and preserves the generated subcategory, so Proposition~\ref{prop:formal} applies; the shifts cancel along every bar path, as after Proposition~\ref{prop:pseudoheight}. Finally $g_*P_x(g^{-1})_*=P_{g(x)}$.
\end{proof}

\subsection{Compression and powers}

Let $i\colon\cA\hookrightarrow\cD$ be admissible, and put
\[
  \Pi=ii^!,\qquad C=\Cone(\Pi\to\id_\cD),\qquad F_\cA=i^!Fi .
\]
Then $(F_\cA)^m=i^!F(\Pi F)^{m-1}i$, and applying the counit at the intermediate positions defines
\[
  \mu_m\colon(F_\cA)^m\longrightarrow i^!F^mi,\qquad D_m(F)=\Cone(\mu_m).
\]

\begin{proposition}\label{prop:compression}
For $m\ge2$, the cofibre $D_m(F)$ has a finite filtration whose graded terms, for $1\le s\le m-1$, are
\begin{equation}\label{eq:compression}
  \bigoplus_{\substack{J\subset\{1,\dots,m-1\}\\|J|=s}}i^!FX_1F\cdots X_{m-1}Fi,\qquad X_r=\begin{cases}C,&r\in J,\\ \Pi,&r\notin J.\end{cases}
\end{equation}
In particular there is a distinguished triangle $(F_\cA)^2\to i^!F^2i\to i^!FCFi$. If $F(\cA)\subset\cA$, every $\mu_m$ is an equivalence.
\end{proposition}

\begin{proof}
Filter each intermediate identity in $i^!F\,\id_\cD\,F\cdots\id_\cD\,Fi$ by $\Pi\to\id_\cD$, with graded terms $\Pi$ and $C$. Exactness of composition gives the total filtration; its initial term, with all $X_r=\Pi$, is $(F_\cA)^m$, and taking the cofibre gives \eqref{eq:compression}. If $F$ preserves $\cA$, read each graded summand from the right: its values lie in $\cA$ until the first occurrence of $C$, which annihilates them.
\end{proof}

Put $\delta_m(F)=\inf\RHom(\id_\cA,D_m(F))$. Normal localization and the triangle of $\mu_m$ give
\begin{equation}\label{eq:compression-range}
  H^tC_{F^m}(\cD)\simeq H^t\RHom\bigl(\id_\cA,(F_\cA)^m\bigr)\qquad\text{for }t\le\min\{h_{F^m}-2,\ \delta_m(F)-1\}.
\end{equation}
For $F=S_\cD$ this compares the ambient groups $\Hom(\id,S_\cD^m[t])$ with the intrinsic groups $\Hom(\id_\cA,S_\cA^m[t])$ of the orthogonal.

\section{The two-sided bar calculation and point height}\label{sec:point-height}

In this section we compute the height of the point coefficient for exceptional collections of line bundles on surfaces. The geometric input, on left-orthogonal divisors, is due to Elagin and to \cite{LZ26}; the new part is a two-sided bar calculation which identifies the normal complex of the point coefficient and bounds it.

\subsection{Geometric input}

Let $X$ be a smooth connected projective complex surface with $p_g(X)=q(X)=0$. A divisor class $D$ is \emph{left-orthogonal} if $H^q(X,\cO_X(-D))=0$ for every $q$. Riemann--Roch gives
\begin{equation}\label{eq:lo}
  D^2+K_XD=-2,\qquad\chi(\cO_X(D))=-K_XD=D^2+2 .
\end{equation}

\begin{lemma}[Elagin, {\cite[Theorem~1.1]{Ela16}}]\label{lem:elagin}
The support of a nonzero effective left-orthogonal divisor is connected, and every prime component is a smooth rational curve.
\end{lemma}

Let $\bE=(L_1,\dots,L_N)$ be a nonempty exceptional collection of line bundles. Exceptionality implies $H^{>0}(X,\cO_X)=0$, so the preceding hypotheses hold. Each difference $L_j\otimes L_i^{-1}$ with $i<j$ is left-orthogonal, and a nonzero map $L_i\to L_j$, $i<j$, has a nonzero effective left-orthogonal zero divisor. We assume
\begin{equation}\label{eq:nonneg}
\begin{gathered}
  K_XC\ge0\ \text{ for every prime component $C$ of the zero divisor}\\
  \text{of every nonzero map }L_i\to L_j,\ i<j .
\end{gathered}
\end{equation}
This holds, for instance, whenever $K_X$ is nef. Put
\[
  A=\Ext^*_X\Bigl(\bigoplus_iL_i,\bigoplus_iL_i\Bigr),\qquad e_jA^qe_i=\Ext^q(L_i,L_j),\qquad R=A^0 .
\]
Thus $A$ is concentrated in degrees $0,1,2$; it is the cohomology algebra of the directed dg algebra $\cE$ of Section~\ref{sec:bar}. The \emph{Hom graph} has vertices $1,\dots,N$ and an edge $\{i,j\}$ whenever $i<j$ and $\Hom(L_i,L_j)\ne0$.

\begin{proposition}[{\cite[Lemma~3.2 and Proposition~3.8]{LZ26}}]\label{prop:LZ}
Under \eqref{eq:nonneg}, every nonzero forward Hom space is one-dimensional. For each connected component $\Lambda=\{a_1<\dots<a_m\}$ of the Hom graph, put $R_\Lambda=(\sum_{i\in\Lambda}e_i)R$. Then
\[
  R=\prod_\Lambda R_\Lambda,\qquad R_\Lambda\simeq\C Q_\Lambda,\qquad Q_\Lambda=(a_1\to\dots\to a_m),
\]
with the indicated vertex idempotents. Moreover $A^1$ is projective as a left $R$-module.
\end{proposition}

For each nonzero $\Hom(L_i,L_j)$, let $C_{ij}$ be the zero divisor of a generator. Our convention identifies $e_{a_s}R_\Lambda e_{a_r}$ with the paths from $a_r$ to $a_s$.

\subsection{Evaluation modules}

The complement of the Hom divisors is the nonempty open subset
\begin{equation}\label{eq:open}
  X^\circ_\bE=X\setminus\bigcup_{\substack{i<j\\ \Hom(L_i,L_j)\ne0}}\Supp C_{ij}.
\end{equation}
For $x\in X^\circ_\bE$, define
\[
  V_x=\bigoplus_iL_i|_x,\qquad U_x=\bigoplus_i(L_i|_x)^\vee .
\]
Evaluation at $x$ makes $V_x$ a left $R$-module and $U_x$ a right $R$-module: for $a\in e_jRe_i$ the actions are $v\mapsto a_xv$ and $u\mapsto u\circ a_x$.

\begin{lemma}\label{lem:evaluation}
Both evaluation modules are projective. If $c$ is the number of connected components of the Hom graph, then $\dim(U_x\otimes_RV_x)=c$.
\end{lemma}

\begin{proof}
Since $x$ lies on no Hom divisor, every nonzero Hom arrow evaluates to an isomorphism at $x$. On $\Lambda=\{a_1<\dots<a_m\}$ choose bases so that consecutive arrows are identities. Then
\[
  V_x|_\Lambda\simeq R_\Lambda e_{a_1},\qquad U_x|_\Lambda\simeq e_{a_m}R_\Lambda,
\]
which are projective on the indicated sides, and
\[
  e_{a_m}R_\Lambda\otimes_{R_\Lambda}R_\Lambda e_{a_1}\simeq e_{a_m}R_\Lambda e_{a_1}\simeq\C,
\]
spanned by the full path from $a_1$ to $a_m$. Different components have zero mixed tensor products, and summing gives dimension $c$.
\end{proof}

\subsection{Relative resolutions}

Write $R=S_0\oplus J_0$, where $S_0=\bigoplus_i\C e_i$ and $J_0$ is the ideal of nonstationary paths.

\begin{lemma}[Normalized relative bar resolution]\label{lem:bar-resolution}
For a left $R$-module $V$, the augmented normalized complex
\begin{equation}\label{eq:bar-resolution}
  \cdots\longrightarrow R\otimes_{S_0}J_0^{\otimes2}\otimes_{S_0}V\longrightarrow R\otimes_{S_0}J_0\otimes_{S_0}V\longrightarrow R\otimes_{S_0}V\longrightarrow V
\end{equation}
is a projective resolution, with the term of length $p$ in cohomological degree $-p$ and differential the alternating sum of the multiplications of adjacent factors, including the final action on $V$.
\end{lemma}

\begin{proof}
This is the normalized bar resolution for $S_0\hookrightarrow R$ \cite[8.6.12 and Exercise~8.6.4]{Wei94}. The augmentation is $r\otimes v\mapsto rv$, and
\begin{multline*}
  d(r\otimes a_1\otimes\dots\otimes a_p\otimes v)=ra_1\otimes a_2\otimes\dots\otimes a_p\otimes v+\sum_{s=1}^{p-1}(-1)^sr\otimes\dots\otimes a_sa_{s+1}\otimes\dots\otimes v\\
  +(-1)^pr\otimes a_1\otimes\dots\otimes a_{p-1}\otimes a_pv .
\end{multline*}
After restriction to $S_0$, a contracting homotopy is $s_{-1}(v)=1\otimes v$ and $s_p(r\otimes a_1\otimes\dots\otimes a_p\otimes v)=1\otimes\bar r\otimes a_1\otimes\dots\otimes a_p\otimes v$, where $\bar r$ is the $J_0$-component of $r$; the identity $ds+sd=\id$ follows by cancellation of adjacent faces. Each term is induced from a module over the semisimple algebra $S_0$, hence is projective, and directedness bounds the length.
\end{proof}

Tensoring \eqref{eq:bar-resolution} on the left with a right module $U$ computes $U\otimes^L_RV$. The same construction works in the presence of auxiliary module actions, and can be inserted successively between bimodules.

\begin{lemma}\label{lem:hereditary}
Let $R$ be a finite product of path algebras of finite acyclic quivers. Every left $R$-module has a projective resolution of length at most one. Consequently, if $M\in D^{\ge d}(R^{\op})$ and $N$ is a left $R$-module in degree zero, then
\begin{equation}\label{eq:tensor-bound}
  M\otimes^L_RN\in D^{\ge d-1}(\C),
\end{equation}
and the bound is $d$ if $N$ is projective. These statements remain valid in the presence of auxiliary bimodule actions.
\end{lemma}

\begin{proof}
For a representation $V=(V_i,V_a)$, the standard resolution is
\begin{equation}\label{eq:quiver-resolution}
  0\to\bigoplus_{a\colon i\to j}Re_j\otimes_\C V_i\longrightarrow\bigoplus_iRe_i\otimes_\C V_i\longrightarrow V\to0,
\end{equation}
whose first map sends $b\otimes v$ to $ba\otimes v-b\otimes V_a(v)$ and whose last map is evaluation \cite[Proposition~1.4.1]{Bri12}. Represent $M$ by its canonical truncation in degrees at least $d$ and resolve $N$ by projectives in degrees $-1,0$; the totalization of their tensor product vanishes below $d-1$ and computes $M\otimes^L_RN$. If $N$ is projective, only degree zero is needed. The resolution \eqref{eq:quiver-resolution} is functorial, so it preserves auxiliary actions.
\end{proof}

\subsection{Positive-degree arrows and the gap filtration}

The following argument adapts the positive-degree filtration of \cite[Lemmas~3.9--3.15]{LZ26} to a two-sided bar complex. The projectivity of the final module replaces the derived trace in that calculation and improves the lower bound by one.

Put $A=R\oplus A^1\oplus A^2$ and $J=J_0\oplus A^1\oplus A^2$, and regard $A$ as a graded algebra with zero differential. If a right $R$-module $U$ and a left $R$-module $V$ are placed in degree zero, extend them to $A$-modules by letting $A^{>0}$ act by zero. Define
\begin{equation}\label{eq:two-sided}
  B_A(U,V)=\bigoplus_{p\ge0}\bigl(U\otimes_{S_0}J^{\otimes_{S_0}p}\otimes_{S_0}V\bigr)[p],
\end{equation}
with the whole-summand shift convention of Section~\ref{sec:bar}. For $z=u\otimes a_1\otimes\dots\otimes a_p\otimes v$ the differential is $bz=\sum_{\nu=0}^p(-1)^\nu d_\nu z$ for $p>0$ and $b=0$ for $p=0$, where
\begin{equation}\label{eq:two-sided-faces}
\begin{aligned}
  d_0z&=ua_1\otimes a_2\otimes\dots\otimes a_p\otimes v,\\
  d_\nu z&=u\otimes\dots\otimes a_\nu a_{\nu+1}\otimes\dots\otimes v\qquad(0<\nu<p),\\
  d_pz&=u\otimes a_1\otimes\dots\otimes a_{p-1}\otimes a_pv .
\end{aligned}
\end{equation}
In contrast with \eqref{eq:faces} there is no cyclic face: the last arrow acts on $V$.

\begin{lemma}[Connectivity of the two-sided bar complex]\label{lem:connectivity}
Suppose $R$ is a product of path algebras of linear chains, $A^1$ is projective as a left $R$-module, and $V$ is projective as a left $R$-module. Then $B_A(U,V)\in D^{\ge0}(\C)$.
\end{lemma}

\begin{proof}
Let $\rho(z)$ be the number of factors of positive degree among $a_1,\dots,a_p$, and filter by $F_\ell=\langle z:\rho(z)\le\ell\rangle$. Products never increase $\rho$; a product of two positive factors decreases it, and the actions of positive factors on the endpoints vanish. Hence the associated graded complex keeps exactly the faces involving a factor of degree zero.

Fix the degrees $\mathbf q=(q_1,\dots,q_\ell)\in\{1,2\}^\ell$ of the positive factors, in tensor order, and put $Q=\sum_jq_j$ and $\sigma=\ell-Q$. There are $\ell+1$ gaps of degree-zero factors. For $\mathbf h=(h_0,\dots,h_\ell)$ with all $h_k\ge0$, put $|\mathbf h|=\sum_kh_k$, $H_k(\mathbf h)=\sum_{j<k}h_j$, and
\[
  U^{\mathbf q}_{\mathbf h}=U\otimes_{S_0}J_0^{\otimes h_0}\otimes_{S_0}A^{q_1}\otimes_{S_0}J_0^{\otimes h_1}\otimes_{S_0}\dots\otimes_{S_0}A^{q_\ell}\otimes_{S_0}J_0^{\otimes h_\ell}\otimes_{S_0}V .
\]
Temporarily place all these factors in degree zero. The summand of the associated graded complex corresponding to $\mathbf q$ is
\begin{equation}\label{eq:graded-piece}
  G^{\mathbf q}=\bigoplus_{\mathbf h}U^{\mathbf q}_{\mathbf h}\bigl[|\mathbf h|+\sigma\bigr],
\end{equation}
since the original internal degree is $Q$ and the length is $\ell+|\mathbf h|$, giving total degree $Q-\ell-|\mathbf h|$.

For the $k$-th gap let its neighbours be $M_k,M_{k+1}$, where $M_0=U$, $M_j=A^{q_j}$ for $1\le j\le\ell$, and $M_{\ell+1}=V$, and write its degree-zero factors as $c_{k,1},\dots,c_{k,h_k}$. For $h_k>0$ and $0\le s\le h_k$, let $\partial^{\mathbf h}_{k,s}$ be the face which multiplies $m_k\otimes c_{k,1}$ (for $s=0$), $c_{k,s}\otimes c_{k,s+1}$ (for $0<s<h_k$), or $c_{k,h_k}\otimes m_{k+1}$ (for $s=h_k$), leaving every other factor fixed; its target is $U^{\mathbf q}_{\mathbf h-e_k}$. Put $\partial^{\mathbf h}_k=\sum_{s=0}^{h_k}(-1)^s\partial^{\mathbf h}_{k,s}$, and $\partial^{\mathbf h}_k=0$ when $h_k=0$. There are $k+H_k(\mathbf h)$ factors before the $k$-th gap, so the differential of the associated graded complex is
\begin{equation}\label{eq:graded-differential}
  b|_{U^{\mathbf q}_{\mathbf h}}=\sum_{k=0}^\ell(-1)^{k+H_k(\mathbf h)}\partial^{\mathbf h}_k .
\end{equation}

As in \cite[Lemma~3.10]{LZ26}, inserting the normalized resolution of Lemma~\ref{lem:bar-resolution} in every gap computes
\begin{equation}\label{eq:iterated-tensor}
  D^{\mathbf q}=U\otimes^L_RA^{q_1}\otimes^L_R\dots\otimes^L_RA^{q_\ell}\otimes^L_RV,
\end{equation}
with the $A^{q_j}$ still placed in degree zero: the resulting terms are $U^{\mathbf q}_{\mathbf h}[|\mathbf h|]$, with differential $\delta=\sum_k(-1)^{H_k(\mathbf h)}\partial^{\mathbf h}_k$. The sign change $\Phi|_{U^{\mathbf q}_{\mathbf h}}=(-1)^{\sum_k(k+\sigma)h_k}\id$ satisfies $\Phi b=(-1)^\sigma\delta\Phi$, and therefore
\begin{equation}\label{eq:graded-identification}
  G^{\mathbf q}\simeq D^{\mathbf q}[\sigma].
\end{equation}

Let $s$ be the number of indices with $q_j=2$. Each derived tensor product with $A^1$ preserves a lower cohomological bound, because $A^1$ is left projective; each tensor product with $A^2$ lowers it by at most one, by Lemma~\ref{lem:hereditary}; and the final tensor product with $V$ preserves it, because $V$ is projective. Starting from $U$ in degree zero we get $D^{\mathbf q}\in D^{\ge-s}$. But $Q=\ell+s$, so $\sigma=-s$, and \eqref{eq:graded-identification} gives $G^{\mathbf q}\in D^{\ge0}$. For $\ell=0$ the same statement is $U\otimes^L_RV=U\otimes_RV$ in degree zero. All bar lengths are bounded by $N-1$, so the filtration is finite, and the long exact sequences of its filtration triangles prove the assertion for $B_A(U,V)$.
\end{proof}

\subsection{Identification with the point normal complex}

Let $\cE$ be the directed dg algebra of the $L_i$, and define its right and left dg modules
\[
  \mathcal U_x=\bigoplus_i\RHom_X(L_i,\cO_x),\qquad\mathcal W_x=\bigoplus_i\RHom_X(\cO_x,L_i).
\]

\begin{proposition}\label{prop:point-normal}
Let $\cB=\langle L_1,\dots,L_N\rangle$ and let $x,y\in X$ be closed points. Then
\begin{equation}\label{eq:point-normal}
  N_{K_{\cO_x,\cO_y}}(\cB,X)\simeq\RHom_X(\cO_x,Q_\cB\cO_y)\simeq\mathcal U_y\otimes^L_\cE\mathcal W_x ,
\end{equation}
and the total complex of the first page of its bar spectral sequence is
\begin{equation}\label{eq:first-page}
  \bigl(E_1,d_1\bigr)\simeq B_A(U_y,V_x)\otimes_\C\tau_x[-2].
\end{equation}
For $x=y$ this is the normal complex $N_{P_x}(\cB,X)$ of the point coefficient.
\end{proposition}

\begin{proof}
The first identification follows from Theorem~\ref{thm:external}, and for $x=y$ from Proposition~\ref{prop:point}. For the second, let $L=\bigoplus_iL_i$ with its left $\cE$-action. The evaluation morphism
\[
  \RHom_X(L,F)\otimes^L_\cE L\longrightarrow F
\]
becomes an equivalence after applying each $\RHom_X(L_i,-)$, by the dg Yoneda lemma. Its source belongs to $\cB$, since the finite directed bar resolution expresses it by cones of finite sums of the $L_i$, and its cone belongs to $\cB^\perp$; it is therefore the projection $Q_\cB F\to F$ \cite{Kel06}. For $F=\cO_y$, applying $\RHom_X(\cO_x,-)$ to that finite resolution replaces $L$ by $\mathcal W_x$, proving \eqref{eq:point-normal}.

Local freeness and local duality give $H^*(\mathcal U_y)=U_y$ and $H^*(\mathcal W_x)=V_x\otimes\tau_x[-2]$. Each endpoint has cohomology in a single degree, so a positive-degree element of $A$ acts by zero on endpoint cohomology, and degree-zero elements act by evaluation and dual evaluation, as in Lemma~\ref{lem:evaluation}. Taking internal cohomology of the dg bar model proves \eqref{eq:first-page}; moving the line $\tau_x$ of even degree past other factors produces no sign.
\end{proof}

\begin{theorem}[Point height]\label{thm:point-height}
Let $\bE$ be a nonempty exceptional collection of line bundles on a smooth projective complex surface satisfying \eqref{eq:nonneg}, for example on a surface with $K_X$ nef. For every $x\in X^\circ_\bE$,
\[
  \alpha_{P_x}(X)=0,\qquad h_{P_x}(\bE)=2 .
\]
If $\Db(X)=\langle\cA,\bE\rangle$, restriction gives
\begin{equation}\label{eq:point-restriction}
  \C\xrightarrow{\ \sim\ }H^0R_{P_x}(\cA\subset X),\qquad T_xX\hookrightarrow H^1R_{P_x}(\cA\subset X).
\end{equation}
\end{theorem}

\begin{proof}
Proposition~\ref{prop:LZ} supplies the description of $R$ and the left projectivity of $A^1$, and Lemma~\ref{lem:evaluation} the left projectivity of $V_x$. Lemma~\ref{lem:connectivity} and \eqref{eq:first-page} with $y=x$ give
\[
  E_2^{-p,q}(P_x)=0\qquad(q-p<2).
\]
All later pages vanish in these bidegrees, so $h_{P_x}\ge2$. Because the vanishing is established on the second page, it holds whatever the later differentials are; in particular the argument is independent of any formality property of $\cE$.

To prove equality, consider the bidegree $(0,2)$. The first differential into this term imposes exactly the balancing relations $(ua)\otimes v-u\otimes(av)$ for $a\in R$; only arrows of degree zero occur, since the terminal factor already has degree two. Hence
\[
  E_2^{0,2}(P_x)=(U_x\otimes_RV_x)\otimes\tau_x\ne0,
\]
of dimension the number of components of the Hom graph. For $r\ge2$, a differential into this term comes from $E_r^{-r,r+1}$, of total degree one, which is zero; a differential out of it has target with positive first index, which is also zero. Thus $E_\infty^{0,2}=E_2^{0,2}\ne0$, and $h_{P_x}=2$. Proposition~\ref{prop:point} identifies the ambient groups in degrees zero and one, and Proposition~\ref{prop:range} gives \eqref{eq:point-restriction}.
\end{proof}

\subsection{Mixed point heights and the first normal group}
For points $x,y$ put $N_{xy}=\RHom_X(\cO_x,Q_\cB\cO_y)$ and
$m_{xy}=\dim H^2N_{xy}$. With our external coefficient convention this
is the normal complex of $K_{\cO_x,\cO_y}$.

\begin{proposition}[Mixed point heights]\label{prop:mixed}
Let $\bE$ satisfy the hypotheses of Theorem~\ref{thm:point-height}, and let
$c$ be the number of components of its Hom graph. For all $x,y\in
X^\circ_{\bE}$ we have $\inf N_{xy}=2$ and $m_{xy}\ge c$.
\end{proposition}
\begin{proof}
The cohomological bar page in Proposition~\ref{prop:point-normal}
has endpoints $U_y$ and $V_x$ with the common factor $\tau_x[-2]$.
Lemma~\ref{lem:connectivity} gives vanishing below total degree two
on the second page. The term
\[
 E_2^{0,2}=(U_y\otimes_R V_x)\otimes\tau_x
\]
has one dimension per chain by Lemma~\ref{lem:evaluation}. No higher
differential can enter it from total degree one or leave it for a
positive first index. This proves both assertions. Contributions
$E_\infty^{-p,p+2}$ for $p>0$ enter the same total degree.
\end{proof}
\begin{theorem}[Point-like objects in the orthogonal]\label{thm:point-like}
Let $\Db(X)=\langle\cA,\bE\rangle$ as in Proposition~\ref{prop:mixed},
and put $W_x=i^*\cO_x$. On $X^\circ_{\bE}$ one has
\[
 \End(W_x)=\C,\qquad T_xX\hookrightarrow\Ext^1(W_x,W_x),
 \qquad \dim\Ext^1(W_x,W_x)=1+m_{xx}.
\]
For $x\ne y$, one has $\Ext^{\le0}(W_x,W_y)=0$ and
$\dim\Ext^1(W_x,W_y)=m_{xy}$. The product of the images of
any two tangent vectors in $\Ext^2(W_x,W_x)$ vanishes.
\end{theorem}
\begin{proof}
Apply the triangle of Theorem~\ref{thm:external} to $\cO_x,\cO_y$.
For distinct points its ambient term vanishes. On the diagonal, height
two gives the assertions in degrees zero and one apart from the dimension.
Choose a nonzero evaluation $t:L_i\to\cO_x$. It factors through
$Q_\cB\cO_x$ by adjunction. By Serre duality, composition with $t[2]$
maps $\Ext^2(\cO_x,L_i)$ onto $\Ext^2(\cO_x,\cO_x)$. Its dual
sends the identity, after trivializing $\omega_X|_x$, to the evaluation
$t$ up to a nonzero scalar. Hence the map $H^2N_{xx}\to\bigwedge^2T_xX$
is surjective. The exact sequence of the triangle gives the dimension,
and the map from $\Ext^2(\cO_x,\cO_x)$ to $\Ext^2(W_x,W_x)$ is zero.
Apply compatibility with products from Theorem~\ref{thm:simultaneous}.
\end{proof}

\begin{example}\label{ex:point-multiplicity}
On $X=\PP^1\times\PP^1$ take $L_1=\cO$ and $L_2=\cO(1,-n)$
with $n\ge2$. The reverse line bundle $\cO(-1,n)$ is acyclic.
The forward cohomology is confined to degree one and has dimension
$2(n-1)$ by K\"unneth. Thus the nonnegativity condition is vacuous
and $c=2$. In the normal bar complex the length-zero part contributes
two dimensions in degree two, and the single-arrow part contributes
$2(n-1)$ more in the same degree. All other terms vanish.
Consequently $m_{xy}=2n$ for all $x,y$. In particular, for $n=2$
the mixed dimension is four and $\dim\Ext^1(W_x,W_x)=5$.
\end{example}

\begin{corollary}\label{cor:hereditary-surface}
In Theorem~\ref{thm:point-like} put $G_{ij}=\chi(L_i,L_j)$ and let
$\mathbf1$ be the column vector of ones. If $\cA$ is the bounded
derived category of an Ext-finite hereditary abelian category, then
$N_{xy}$ is concentrated in degree two and
\[
  m_{xy}=\mathbf1^{t}G^{-1}\mathbf1,\qquad
  \dim\Ext^1(W_x,W_x)=1+\mathbf1^{t}G^{-1}\mathbf1.
\]
\end{corollary}
\begin{proof}
The column of coefficients of $[Q_\cB\cO_y]$ in the line-bundle basis
is $G^{-1}\mathbf1$, by adjunction and $\chi(L_i,\cO_y)=1$.
Since $\chi(\cO_x,L_i)=1$, pairing on the left gives
$\chi(N_{xy})=\mathbf1^tG^{-1}\mathbf1$. For $x\ne y$, the proof
of Theorem~\ref{thm:two-probes} gives concentration of
$\RHom(W_x,W_y)$ in degree one. For $x=y$, heredity and
indecomposability give vanishing of $\Ext^{>1}(W_x,W_x)$. The
projection triangle and Proposition~\ref{prop:mixed} now give the
stated concentration on and off the diagonal.
\end{proof}

\subsection{Concentration of the normal point complex}

\begin{proposition}\label{prop:normal-pure-two}
In the setting of Theorem~\ref{thm:point-height}, assume in addition
that $A^2=0$. This holds, in particular, for a connected Hom graph.
For $x,y\in X^\circ_{\bE}$ the normal complex is concentrated in degree
two, of dimension
\[
 M=\sum_{j=0}^{n-1}\dim\bigl(U_y\otimes_R(A^1)^{\otimes_R j}
                                  \otimes_R V_x\bigr)\ge1,
\]
where $n$ is the number of lines and the zeroth tensor power is $R$.
This number is independent of $x,y$ on this open set.
\end{proposition}
\begin{proof}
If the Hom graph is connected, every forward pair in its chain has
nonzero Hom by Proposition~\ref{prop:LZ}. Its difference is an effective
divisor $D$, and $H^0(K_X-D)=0$ since $p_g(X)=0$ and multiplication by
a defining section of $D$ is injective. Serre duality gives $A^2=0$.
Under this vanishing all positive-degree arrows have degree one. In
the positive-input filtration of Lemma~\ref{lem:connectivity}, its
$j$th piece is $U_y\otimes_R^L(A^1)^{\otimes_R^L j}\otimes_R^L V_x$,
with no net shift. Left projectivity of $A^1$ and $V_x$ makes every
tensor underived, so all pieces are in degree zero. Both the resulting
spectral sequence and the remaining bar spectral sequence collapse
by total degree. The external shift is $[-2]$, proving the formula.
The vertex descriptions in Lemma~\ref{lem:evaluation} make the
dimensions independent of $x,y$, and the zeroth term is the positive
number of Hom components.
\end{proof}

\section{Rationality with one or two arbitrary exceptional objects}\label{sec:rationality}

Orlov's rationality conjecture asks whether a smooth projective complex surface admitting a full exceptional collection is rational \cite[Introduction]{Via17}. It is known for collections of line bundles \cite[Theorem~1.2]{LZ26}. This section proves part~(ii) of Theorem~\ref{thm:intro-rational}, and then treats a pair of arbitrary objects next to a collection of line bundles.

\begin{theorem}[Rationality of surfaces]\label{thm:rational}
Let $X$ be a smooth connected projective surface over $\C$. Suppose that $\Db(X)$ admits a full exceptional collection all but at most one of whose members are line bundles. Then $X$ is rational.
\end{theorem}

The proof has two parts. First, the point-height estimate excludes an orthogonal generated by a single exceptional object (Section~\ref{subsec:one-object}). Second, on a nonminimal surface we deform the centres of the blowups and contract exceptional curves, preserving the form of the collection; this geometric reduction follows \cite[Sections~3.4--3.7]{LZ26}, with additional care for the arbitrary member.

\subsection{An orthogonal generated by one object}\label{subsec:one-object}

\begin{lemma}\label{lem:one-object-complex}
Suppose $\Db(X)=\langle\cA,\cB\rangle$ with $\cA=\langle G\rangle$ generated by an exceptional perfect object. For $x\in X$, put $V_x=L\iota_x^*G$ and $b_a=\dim_\C H^a(V_x)$. Then
\begin{equation}\label{eq:one-object}
  R_{P_x}(\cA\subset X)\simeq\RHom_\C(V_x,V_x),\qquad\dim H^tR_{P_x}=\sum_ab_ab_{a+t}.
\end{equation}
\end{lemma}

\begin{proof}
The equivalence $\Perf(\C)\to\cA$ sends $W$ to $W\otimes_\C G$. Under this equivalence, $i^!P_xi$ is the tensor product with
\[
  \RHom_X(G,G\otimes^L\cO_x)\simeq\RHom_X(G,\iota_{x*}L\iota_x^*G)\simeq\RHom_\C(V_x,V_x).
\]
The complex of transformations from the identity of $\Perf(\C)$ to the tensor product with a complex $C$ is $C$ itself, by evaluation on the generator $\C$; this gives the first identity. Every complex over a field is the direct sum of its shifted cohomology, so $H^t\RHom(V_x,V_x)=\bigoplus_a\Hom(H^a(V_x),H^{a+t}(V_x))$, which gives the dimension formula; the sum is finite because $G$ is perfect.
\end{proof}

\begin{proposition}\label{prop:one-object}
Let $X$ be a positive-dimensional smooth projective complex variety and let $\cB\subset\Db(X)$ be admissible. If $h_{P_x}(\cB\subset X)\ge2$ for some point $x$, then $\cB^\perp$ is not generated by one exceptional object. This includes the case $\cB=0$, in which the normal complex is zero and the height is $+\infty$.
\end{proposition}

\begin{proof}
Assume $\cB^\perp=\langle G\rangle$. Propositions~\ref{prop:range} and~\ref{prop:point} give
\[
  H^0R_{P_x}\simeq\C,\qquad T_xX\hookrightarrow H^1R_{P_x}.
\]
By Lemma~\ref{lem:one-object-complex}, $\sum_ab_a^2=1$. The $b_a$ are nonnegative integers, so exactly one of them equals one and all the others vanish. Therefore $\sum_ab_ab_{a+1}=0$, that is, $H^1R_{P_x}=0$, contradicting $\dim T_xX>0$.
\end{proof}

\begin{remark}\label{rem:two-degrees}
The argument uses two adjacent degrees. Positivity of a height only shows that the orthogonal of the line bundles is nonzero, which is compatible with an orthogonal generated by one exceptional object. The contradiction comes from using $H^0R_{P_x}=\C$ and $T_xX\hookrightarrow H^1R_{P_x}$ together, for the derived fibre of the generator. For an orthogonal generated by several exceptional objects the restricted complex is no longer a matrix algebra of a single complex, and the sum-of-squares argument does not apply; this is why the case of several arbitrary members of Orlov's conjecture remains open.
\end{remark}

\subsection{Mutations and the order of the collection}

For an exceptional object $E$, left and right mutations are defined by the triangles
\begin{equation}\label{eq:mutations}
  \RHom(E,F)\otimes E\longrightarrow F\longrightarrow L_EF,\qquad R_EF\longrightarrow F\longrightarrow\RHom(F,E)^\vee\otimes E .
\end{equation}
Thus an exceptional pair $(E,F)$ can be replaced by $(L_EF,E)$ or by $(F,R_FE)$; both replacements preserve exceptionality and the generated subcategory. These are the standard mutations, with shifts fixed by \eqref{eq:mutations}; compare \cite[Definition~2.8 and Lemma~2.11]{Hu19}.

\begin{lemma}\label{lem:reorder}
A full exceptional collection with one designated arbitrary member and $n$ line-bundle members can be mutated to $(G,L_1,\dots,L_n)$, where the line bundles are the original ones, in their original relative order.
\end{lemma}

\begin{proof}
If a line bundle $L$ immediately precedes the designated member $G$, replace $(L,G)$ by $(L_LG,L)$. No line bundle changes, and each replacement decreases by one the number of line bundles preceding the designated member.
\end{proof}

\subsection{Deformation to distinct centres}

We use the family of ordered blowups of \cite[Propositions~1.4--1.5]{Roe04}, in the form of \cite[Lemma~3.20]{LZ26}, together with Hu's deformation theorem for exceptional objects.

\begin{lemma}[{\cite[Lemma~3.20]{LZ26}}]\label{lem:family}
For a smooth projective connected surface $S$ and $r\ge0$, there are a smooth irreducible projective variety $B_r$, a smooth projective family $f_r\colon\cX_r\to B_r$, and a morphism $g_r\colon\cX_r\to S$ with the following properties:
\begin{enumerate}[(i)]
  \item every ordered tower of $r$ point blowups over $S$ is a fibre;
  \item the locus of distinct centres is the nonempty open set $\Conf_r(S)=S^r\setminus\bigcup_{a<b}\{p_a=p_b\}$;
  \item there are line bundles on $\cX_r$ restricting on each fibre to the total transforms of the successive exceptional divisors;
  \item every line bundle on a fibre extends to $\cX_r$.
\end{enumerate}
\end{lemma}

\begin{proof}
Start with $B_0=\operatorname{Spec}\C$ and $\cX_0=S$. Inductively let $B_r=\cX_{r-1}$ and $\widetilde{\cX}_{r-1}=\cX_{r-1}\times_{B_{r-1}}B_r$. The diagonal gives a section $\sigma_r\colon B_r\to\widetilde{\cX}_{r-1}$; let
\[
  \beta_r\colon\cX_r=\Bl_{\sigma_r(B_r)}\widetilde{\cX}_{r-1}\longrightarrow\widetilde{\cX}_{r-1},\qquad f_r=\mathrm{pr}_2\circ\beta_r .
\]
The section has relative codimension two, so $\cX_r\to B_r$ is a smooth family of surfaces, whose fibre over $b\in\cX_{r-1}$ is the blowup of the fibre of $f_{r-1}$ through $b$ at the point $b$. Induction gives projectivity, smoothness and irreducibility, and realizes every ordered tower, including infinitely near centres \cite[Propositions~1.4--1.5]{Roe04}.

The images of the centres define $B_r\to S^r$. When they are distinct, every successive centre has a unique lift away from the earlier exceptional curves, so the inverse image of $\Conf_r(S)$ is isomorphic to it. Let $H_{r,r}$ be the line bundle of the last relative exceptional divisor and, for $a<r$, let $H_{a,r}=\beta_r^*\mathrm{pr}_1^*H_{a,r-1}$; on a fibre these restrict to $\cO(E_a^{\mathrm{tot}})$. For a fibre $\pi\colon X\to S$ one has
\[
  \Pic(X)=\pi^*\Pic(S)\oplus\bigoplus_{a=1}^r\Z[E_a^{\mathrm{tot}}]
\]
\cite[Section~1]{Roe04}, and the line bundle $\pi^*M\otimes\cO_X(\sum_am_aE_a^{\mathrm{tot}})$ extends to $g_r^*M\otimes\bigotimes_aH_{a,r}^{m_a}$.
\end{proof}

\begin{lemma}[{\cite[Proposition~3.9]{Hu19}}]\label{lem:open}
Let $f\colon\cX\to T$ be a smooth projective family with connected fibres over a smooth complex variety, and let $\widetilde E_0,\dots,\widetilde E_n$ be perfect complexes on $\cX$. If their restrictions form a full exceptional collection on the fibre over $t_0$, the same is true over a Zariski open neighbourhood of $t_0$.
\end{lemma}

\begin{proof}
Consider the perfect complexes on $T$
\[
  D_a=\Cone\bigl(\cO_T\to Rf_*\RHom(\widetilde E_a,\widetilde E_a)\bigr),\qquad K_{ba}=Rf_*\RHom(\widetilde E_b,\widetilde E_a)\quad(b>a).
\]
By base change their derived fibres are the cones of the units of the self-Hom complexes and the reverse morphism complexes on the fibres, which vanish at $t_0$. Supports of perfect complexes are closed, so these complexes vanish on a neighbourhood of $t_0$, over which the given objects form a relative exceptional collection. Hu's openness theorem \cite[Proposition~3.9]{Hu19} gives a smaller neighbourhood on which the collection is full.
\end{proof}

\begin{lemma}[{\cite[Lemmas~3.21--3.22]{LZ26}}]\label{lem:rational-curves}
A non-uniruled smooth projective complex surface $S$ contains at most countably many integral rational curves. Every nonempty Zariski open subset of $\Conf_r(S)$ contains a tuple none of whose coordinates lies on such a curve.
\end{lemma}

\begin{proof}
Fix a projective embedding. For each degree $d>0$, the scheme $\Hom_d(\PP^1,S)$ is of finite type. Let $T$ be one of its irreducible components. If the evaluation map $\PP^1\times T\to S$ were dominant, generic smoothness would give a point at which its differential has rank two, and choosing a curve in $T$ through the corresponding parameter in a suitable direction would produce a dominant map from $\PP^1$ times a curve to $S$, contradicting non-uniruledness. So the image of the evaluation map has one-dimensional closure, and all maps in $T$ have the same integral rational image. Every rational curve is the image of its normalization, which gives countability; this is the surface case of \cite[Theorem~25 and~27(1)]{Kol10}.

Enumerate these curves as $C_m$. Inside a nonempty open subset $W\subset\Conf_r(S)$, each $W\cap\mathrm{pr}_a^{-1}(C_m)$ is a proper closed subset. Choose a nonempty analytic ball $\Omega\subset W$; these subsets are closed and nowhere dense in $\Omega$, and by the Baire category theorem their countable union does not cover $\Omega$.
\end{proof}

\begin{proposition}\label{prop:deform}
Let $S$ be non-uniruled, and let $X\to S$ be a tower of $r$ point blowups. Suppose $X$ has a full exceptional collection $(G,L_1,\dots,L_n)$ with $G$ arbitrary and every $L_i$ a line bundle. Then there exist distinct points $p_1,\dots,p_r\in S$, lying on no integral rational curve of $S$, such that $\Bl_{p_1,\dots,p_r}S$ has a full exceptional collection of the same form.
\end{proposition}

\begin{proof}
For $r=0$ there is nothing to prove. Realize $X$ as the fibre over $b_0\in B_r$ by Lemma~\ref{lem:family}, and extend each $L_i$ to a line bundle $\widetilde L_i$ on $\cX_r$. Hu's deformation theorem \cite[Theorem~3.3]{Hu19} extends the exceptional object $G$ to a perfect complex $\widetilde G$ on $\cX_r\times_{B_r}T$ for an \'etale neighbourhood $\nu\colon T\to B_r$, $t_0\mapsto b_0$. Lemma~\ref{lem:open}, applied to $\widetilde G$ and the pullbacks of the $\widetilde L_i$, allows us to shrink $T$ so that this collection is full and exceptional on every fibre.

The image $V=\nu(T)$ is a nonempty open subset of the irreducible variety $B_r$, so its intersection with the locus of distinct centres corresponds to a nonempty open subset $W\subset\Conf_r(S)$. Lemma~\ref{lem:rational-curves} gives a tuple in $W$ avoiding all rational curves, and it lifts to a point of $T$. The corresponding fibre is the blowup at the chosen distinct centres and carries the restricted full exceptional collection.
\end{proof}

\subsection{Extraction of an exceptional-curve pair}

Let $S$ have $p_g=q=0$ and nef canonical divisor, and fix $\pi\colon Y=\Bl_{p_1,\dots,p_r}S\to S$, where the centres are distinct and avoid every integral rational curve of $S$. Denote its disjoint exceptional curves by $E_a$; then
\begin{equation}\label{eq:blowup}
  E_a^2=-1,\qquad E_aE_b=0\ (a\ne b),\qquad K_Y=\pi^*K_S+\sum_aE_a .
\end{equation}
The following arguments on divisors and adjacency are from \cite{LZ26}. They apply to any exceptional collection of line bundles on $Y$, full or not.

\begin{lemma}[{\cite[Lemma~3.24 and Corollary~3.25]{LZ26}}]\label{lem:lo-divisors}
Every effective left-orthogonal divisor $D>0$ on $Y$ is either a single $E_a$ or the pullback of an effective left-orthogonal divisor $D_0>0$ on $S$. In the latter case it is disjoint from the exceptional locus, and $\chi(\cO_Y(D))=-K_SD_0\le0$. If an exceptional pair of line bundles $(L,M)$ has nonzero $\RHom(L,M)$ concentrated in degree zero, then $M\simeq L(E_a)$ for some $a$, and $\RHom(L,M)=\C$.
\end{lemma}

\begin{proof}
By Lemma~\ref{lem:elagin}, the support of $D$ is connected and all its prime components are rational. A component that is not contracted has rational image in $S$, which avoids the centres; it is therefore the full pullback of its image and is disjoint from all the $E_a$. Connectedness leaves two possibilities: $D=mE_a$, or $D=\pi^*D_0$ disjoint from every $E_a$. In the first case \eqref{eq:lo} gives $-m^2-m=-2$, so $m=1$. In the second case $R\pi_*\cO_Y=\cO_S$ and the projection formula give $R\Gamma(S,\cO_S(-D_0))\simeq R\Gamma(Y,\cO_Y(-D))=0$, so $D_0$ is left-orthogonal, and \eqref{eq:lo} gives the formula for the Euler characteristic.

For the exceptional pair, a nonzero morphism has such a zero divisor $D$. The pullback alternative would give $0<h^0(\cO_Y(D))=\chi(\cO_Y(D))\le0$, so $D=E_a$. Finally $0\to\cO_Y\to\cO_Y(E_a)\to\cO_{E_a}(-1)\to0$ gives $R\Gamma(\cO_Y(E_a))=\C$.
\end{proof}

\begin{lemma}[{\cite[Lemma~3.26]{LZ26}}]\label{lem:curve-homs}
If $E\subset Y$ is a $(-1)$-curve, $T=\cO_E(d)$, and $M$ is a line bundle, then
\begin{align}
  \RHom_Y(M,T)&=R\Gamma\bigl(\PP^1,\cO(d-ME)\bigr),\label{eq:hom-to-curve}\\
  \RHom_Y(T,M)&=R\Gamma\bigl(\PP^1,\cO(ME-d-1)\bigr)[-1].\label{eq:hom-from-curve}
\end{align}
Consequently the first complex vanishes exactly when $ME=d+1$, and the second exactly when $ME=d$.
\end{lemma}

\begin{proof}
For the inclusion $\iota\colon E\hookrightarrow Y$ one has $L\iota^*M=\cO_E(ME)$ and $\iota^!M=M|_E\otimes\cO_E(E)[-1]=\cO_E(ME-1)[-1]$. Adjunction gives the two formulas, and $R\Gamma(\PP^1,\cO(t))=0$ exactly for $t=-1$.
\end{proof}

\begin{lemma}[cf.\ the proof of {\cite[Proposition~3.28]{LZ26}}]\label{lem:adjacent}
If an exceptional collection of line bundles $(L_1,\dots,L_n)$ on $Y$ contains a pair with $L_j\simeq L_i(E_a)$, then a pair minimizing $j-i$ is adjacent in the collection.
\end{lemma}

\begin{proof}
Take a minimizing pair, and put $E=E_a$, $\ell=L_iE$ and $T=\cO_E(\ell-1)$. The sequence $0\to L_i\to L_j\to T\to0$ and semiorthogonality imply, for $i<k<j$ and $t=L_kE-\ell$,
\begin{align*}
  \RHom(L_i,L_k)&\simeq\RHom(T,L_k)[1]\simeq R\Gamma(\PP^1,\cO(t)),\\
  \RHom(L_k,L_j)&\simeq\RHom(L_k,T)\simeq R\Gamma(\PP^1,\cO(-t-1)).
\end{align*}
For $t\ge0$ the first is nonzero and concentrated in degree zero; for $t\le-1$ the second is. Lemma~\ref{lem:lo-divisors} then gives a pair differing by an exceptional curve with smaller gap $k-i$ or $j-k$, contradicting minimality. So there is no intermediate $k$.
\end{proof}

\subsection{Contraction and the arbitrary member}

The next proposition extends \cite[Proposition~3.31]{LZ26}. The mutations of line bundles are those of \cite[Lemmas~3.29--3.30]{LZ26}; the new step places the arbitrary member in the pullback component of Orlov's blowup decomposition.

\begin{proposition}\label{prop:contract}
Let $f\colon Y\to Z$ be the contraction of a $(-1)$-curve $E$ on a smooth projective complex surface. Suppose $(G,L_1,\dots,L_n)$ is a full exceptional collection whose line bundles contain an adjacent pair $(L,L(E))$. Then $Z$ has a full exceptional collection with one arbitrary member and $n-1$ line bundles.
\end{proposition}

\begin{proof}
Exceptionality of a line bundle gives $p_g(Y)=q(Y)=0$. Put $d=LE-1$ and $T=\cO_E(d)$. The local self-Ext sheaves of $T$ are $\cO_E$ in degree zero and $\cO_E(-1)$ in degree one, so $\RHom_Y(T,T)=\C$. The exact sequence
\begin{equation}\label{eq:pair-sequence}
  0\longrightarrow L\longrightarrow L(E)\longrightarrow T\longrightarrow0
\end{equation}
and $\RHom(L,L(E))=\C$ identify $L_L(L(E))$ with $T$. Replace $(L,L(E))$ by $(T,L)$.

For a line bundle $M$ of the collection preceding $T$, semiorthogonality gives $\RHom(T,M)=0$, hence $ME=d$ by Lemma~\ref{lem:curve-homs}. The sequence
\begin{equation}\label{eq:twist-sequence}
  0\longrightarrow M(-E)\longrightarrow M\longrightarrow T\longrightarrow0
\end{equation}
identifies $R_TM$ with $M(-E)$. Thus $(M,T)$ becomes $(T,M(-E))$, and $M(-E)E=d+1$. A line bundle following $T$ already satisfies $ME=d+1$, again by semiorthogonality. Moving $T$ to the left across all preceding line bundles therefore gives
\[
  (G,T,M_1,\dots,M_{n-1}),\qquad M_aE=d+1 .
\]
Now mutate the arbitrary member: define
\begin{equation}\label{eq:G-prime}
  G'=R_TG=\Cone\bigl(G\to\RHom_Y(G,T)^\vee\otimes T\bigr)[-1],
\end{equation}
so that the pair $(G,T)$ becomes $(T,G')$. Applying $\RHom_Y(-,T)$ to the defining triangle, the coevaluation identifies $\RHom_Y(\RHom_Y(G,T)^\vee\otimes T,T)\to\RHom_Y(G,T)$ with the identity, so $\RHom_Y(G',T)=0$. Since $\RHom_Y(T,G)=0$, the same triangle gives
\[
  \RHom_Y(G',G')\simeq\RHom_Y(G',G)\simeq\RHom_Y(G,G)=\C .
\]
The two pairs have the same thick span, and the remaining semiorthogonality conditions follow from the triangle. So the new collection is full and exceptional.

Tensor the whole collection with $\cO_Y((d+1)E)$, and write
\[
  \epsilon=\cO_E(-1),\qquad H=G'((d+1)E),\qquad P_a=M_a((d+1)E).
\]
This gives the full exceptional collection $(\epsilon,H,P_1,\dots,P_{n-1})$, with $\RHom_Y(H,\epsilon)=0$ and $P_aE=0$. Hence $P_a=f^*N_a$ for line bundles $N_a$ on $Z$. Orlov's blowup formula \cite[Theorem~4.3]{Orl93} gives
\begin{equation}\label{eq:orlov}
  \Db(Y)=\bigl\langle\cO_E(-1),Lf^*\Db(Z)\bigr\rangle,\qquad Lf^*\Db(Z)=\{F:\RHom_Y(F,\cO_E(-1))=0\}.
\end{equation}
Therefore $H=Lf^*G_Z$ for a perfect object $G_Z$ on $Z$. Full faithfulness of $Lf^*$ \cite[Lemma~4.1]{Orl93} identifies every morphism complex of $(G_Z,N_1,\dots,N_{n-1})$ with that of the corresponding collection on $Y$, so this collection is exceptional. Finally $Rf_*\epsilon=0$ and $Rf_*Lf^*=\id$: for $F\in\Db(Z)$, the object $Lf^*F$ is generated by $(\epsilon,H,P_1,\dots)$, and applying $Rf_*$ shows that $F$ lies in $\langle G_Z,N_1,\dots,N_{n-1}\rangle$. This proves fullness.
\end{proof}

\subsection{Proof of the rationality theorem}

\begin{proof}[Proof of Theorem~\ref{thm:rational}]
If all members are line bundles, apply \cite[Theorem~1.2]{LZ26}. Otherwise Lemma~\ref{lem:reorder} gives a full collection $(G,L_1,\dots,L_n)$. If $n=0$, Proposition~\ref{prop:one-object} with $\cB=0$ gives a contradiction. Hence $n>0$, and exceptionality of $L_1$ implies $p_g(X)=q(X)=0$.

Assume $X$ is not rational, and let $X\to S$ be a tower of $r$ point blowups with $S$ minimal. Then $p_g(S)=q(S)=0$. By the classification of minimal surfaces \cite[Section~5.6]{Deb11}, $S$ is $\PP^2$, ruled, or has nef canonical divisor; a ruled surface with $q=0$ is rational, so $K_S$ is nef, and $S$ is not uniruled.

\emph{The minimal case $r=0$.} Apply Theorem~\ref{thm:point-height} to $\cB=\langle L_1,\dots,L_n\rangle$ on $S$. Its right orthogonal is $\langle G\rangle$, because the collection has the order $(G,L_1,\dots,L_n)$. The point-height theorem gives $h_{P_x}(\cB\subset S)\ge2$ for every point $x$ of the nonempty open set $S^\circ_\bE$, contradicting Proposition~\ref{prop:one-object}.

\emph{The case $r>0$.} Proposition~\ref{prop:deform} produces a surface $Y=\Bl_{p_1,\dots,p_r}S$ carrying a full collection of the same form, with distinct centres avoiding every integral rational curve of $S$. If no two of its line bundles differ by an exceptional curve, Lemma~\ref{lem:lo-divisors} shows that the zero divisor of every nonzero map $L_i\to L_j$, $i<j$, is a pullback disjoint from the exceptional locus, and each of its prime components $C$ satisfies $K_YC=K_S\pi_*C\ge0$. So condition \eqref{eq:nonneg} holds on $Y$, and Theorem~\ref{thm:point-height} and Proposition~\ref{prop:one-object} give the same contradiction as before. Hence some pair of line bundles differs by an exceptional curve; Lemma~\ref{lem:adjacent} makes such a pair adjacent, and Proposition~\ref{prop:contract} contracts the curve, producing a full collection of the same form, with one line bundle fewer, on the blowup of $S$ in the remaining $r-1$ centres. These centres are still distinct and avoid the rational curves of $S$, so the argument can be repeated: at each stage an empty set of line bundles is excluded by Proposition~\ref{prop:one-object}, and a collection without a pair differing by an exceptional curve is excluded by the point-height theorem. After $r$ contractions we reach the minimal case, which has already been excluded. Therefore $X$ is rational.
\end{proof}

\subsection{An exceptional pair and its evaluation at the generic point}
\label{subsec:pair-profile}

Let $E,F$ be an exceptional pair, let $\eta$ be the generic point of $X$, and set
\[
 a_p=\dim_{\C(X)} H^p(E_\eta),\qquad b_q=\dim_{\C(X)} H^q(F_\eta),
 \qquad m_d=\dim_{\C}\Ext^d_X(E,F).
\]
For nonnegative integer vectors $u=(u_p)$, $v=(v_q)$ of finite support define
\[
 D_m(u,v)=1+\sum_{p,q}m_{q-p}u_pv_q-\sum_pu_p^2-\sum_qv_q^2 .
\]

\begin{proposition}\label{prop:pair-profile}
Suppose $\Db(X)=\langle E,F,L_1,\ldots,L_n\rangle$ with line bundles $L_i$
satisfying \eqref{eq:nonneg}. Then $D_m(a,b)\ge2$. If $D_m(a,b)=2$, then
$X$ is rational.
\end{proposition}
\begin{proof}
Let $\cA=\langle E,F\rangle$, and let $W_x$ be the left projection of
$\cO_x$ to $\cA$. For general $x$, mixed point heights at least two imply
$\End(W_x)=\C$ and $\Hom(W_x,W_y)=0$ for distinct general $x,y$.
For an empty line block the assertions hold directly for $W_x=\cO_x$.

Consider the covariant module on the generators defined by
$G\longmapsto\RHom_{\cA}(W_x,G)$. Its isomorphism class determines $W_x$
by the dual form of dg Yoneda on generators. By adjunction and local
duality its values for $G=E,F$ are
\[
 \RHom_X(\cO_x,G)\simeq G|_x^{\mathbf L}\otimes\det(T_xX)[-2].
\]
The common one-dimensional factor and common shift can be discarded.
The strictly unital minimal model of the category on $E,F$ has a single
forward arrow space $\Ext^*(E,F)$ and no composable pairs of
nonidentity arrows. Consequently every higher multiplication vanishes.
The minimal model of each module is determined by the two graded
vertex spaces and the linear actions of the individual arrows, since
a higher module action would also require at least two composable
nonidentity arrows. The action of a class in $\Ext^d(E,F)$ is the
degree-$d$ map on the cohomology of the derived fibers. In particular
a degree-preserving pair of vertex isomorphisms intertwining all arrow
maps gives an isomorphism of the modules.

On a dense open set $U\subset X$ these vertex spaces are the fibers
of vector bundles of ranks $a_p,b_q$. One sees this by generic
freeness of the cohomology sheaves and base change for perfect
complexes. Shrink $U$ to trivialize the bundles and restrict to the
open set on which the mixed point-height assertions hold. Taking
a basis for each $\Ext^d(E,F)$ gives an algebraic map
\[
 f: U\longrightarrow V=\bigoplus_{p,q}
      \Hom(\C^{a_p},\C^{b_q})^{\oplus m_{q-p}} .
\]
Two points have disjoint orbits for
$G=\prod_p\mathrm{GL}_{a_p}\times\prod_q\mathrm{GL}_{b_q}$.
Moreover the stabilizer of each $f(x)$ is the scalar group. Indeed
a degree-preserving intertwiner gives a strict module endomorphism,
and its effect on vertex cohomology is the given matrix. Thus the
algebra of such intertwiners injects into $\End(W_x)=\C$ and contains
the scalar identity.

Let $J$ be the image in $V\times U$ of $(g,x)\mapsto(gf(x),x)$.
It is irreducible constructible of dimension $\dim G+1$, since its
fibers have dimension one. Its projection to $V$ is injective on
closed points, so $\dim V\ge\dim G+1$, proving the first assertion.
If equality holds, the projection $\overline J\to V$ is dominant
and generically one-to-one. Indeed $J$ contains a dense open of its
closure, and the complement has dimension below $\dim V$. In
characteristic zero this projection is birational. The second
projection now gives a dominant rational map from the affine space
$V$ to $X$. A smooth projective unirational complex surface is
rational.
\end{proof}

\begin{theorem}\label{thm:two-generic-bundles}
Suppose a smooth connected projective complex surface has a full exceptional
collection of shifted line bundles and one consecutive exceptional pair
$(E,F)$. Suppose first that
$E_\eta\simeq\C(X)^r[s]$ and $F_\eta\simeq\C(X)^t[u]$
with $r,t>0$. Set $m=\dim\Ext^{s-u}_X(E,F)$. If
\begin{equation}\label{eq:pair-small}
               mrt-r^2-t^2\le1 ,
\end{equation}
then $X$ is rational. In particular for coherent sheaves of positive rank,
$\dim\Hom_X(E,F)\le2$ suffices with no restriction on their higher Ext
groups. Three forward degree-zero maps also suffice for ranks $(1,1)$,
$(1,2)$ or $(2,1)$.

More generally, for arbitrary $E,F$ it suffices that $D_m(u,v)\le2$
for every pair of nonnegative profiles $u\le a$, $v\le b$ satisfying
$\sum_p(-1)^pu_p=\rk E$ and $\sum_q(-1)^qv_q=\rk F$, where in this
statement $m=(m_d)_d$ is the complete list of Ext dimensions.
\end{theorem}
\begin{proof}
Rotate the line bundles preceding the pair to the end by the inverse
Serre functor and discard their shifts. Thus the pair precedes the lines.
Suppose $X$ is nonrational. As in the proof of Theorem~\ref{thm:rational},
move the blowup centers over the minimal model by extending the chosen
exceptional collection. Every dimension $m'_d$ after deformation
is at most $m_d$ after shrinking, by semicontinuity for the relative
perfect Hom complex. On a common dense open of the original fiber the
derived fiber dimensions of $E,F$ equal $a_p,b_q$. The locus on the
total space where all these dimensions are at most $a_p,b_q$ is open
and meets the original fiber. Its image on the base is open by
smoothness. On every integral fiber over this image the two generic
profiles are therefore termwise at most $a,b$. Euler ranks are locally
constant. In particular a profile initially concentrated in one degree
and of positive total dimension is unchanged, and its degree is unchanged.

Choose a fiber with distinct centers avoiding all rational curves on
the minimal model. If the line sequence satisfies \eqref{eq:nonneg},
its actual profiles give $D_{m'}(u,v)\le D_m(u,v)\le2$.
Proposition~\ref{prop:pair-profile} forces equality and rationality.
Otherwise the extraction in Lemma~\ref{lem:adjacent} gives an adjacent
exceptional-divisor pair of lines. Perform the contraction of
Proposition~\ref{prop:contract}, moving the exceptional curve object
$T$ across both $E$ and $F$ when moving it to the front. This preserves
their entire mutual Ext space: both lie in $T^\perp$, and right mutation
through $T$ is fully faithful there, as follows by applying Hom to the
two mutation triangles. Their generic fibers are unchanged since $T$
has proper support. The common twist preserves their degrees and
dimensions, and descent along the fully faithful pullback preserves the
Ext groups and the generic profiles. Repeating terminates, since each
step contracts an exceptional curve. It yields a birational surface with
the line sequence satisfying \eqref{eq:nonneg}, where the preceding
argument gives the contradiction. Finally for one-degree profiles
$D_m=1+mrt-r^2-t^2$, as asserted.
\end{proof}

\subsection{Two arbitrary objects and a numerical Euler criterion}
\begin{theorem}\label{thm:pair-numerical-height}
Suppose $\Db(X)=\langle E,F,L_1,\ldots,L_n\rangle$ is a full
exceptional collection with $E,F$ arbitrary complexes and $L_i$ line
bundles such that $\Ext^2_X(L_i,L_j)=0$ for every $i<j$. Write
$r=\rk E$, $s=\rk F$ and $m=\chi(E,F)$. If
\[
                       r^2+s^2-mrs\ge0,
\]
then $X$ is rational. In particular $|\chi(E,F)|\le2$ suffices,
even if the pair occurs in an arbitrary consecutive position and the
remaining line bundles are taken in their original order.
\end{theorem}
\begin{proof}
There must be a line, since a smooth projective surface has numerical
Grothendieck group of rank at least three. If $X$ were nonrational,
perform the deformation and the exceptional curve contractions in the
proof of Theorem~\ref{thm:rational}, keeping both arbitrary terms. The
deformation preserves Euler forms and ranks; all assumed forward Ext
vanishings persist by upper semicontinuity. In a contraction the new
terms replacing $E,F$ are their common right mutations through the
exceptional curve object $T$, followed by a common twist and descent.
Since both belong to $T^\perp$, the mutation is fully faithful on them
and preserves their mutual Euler form. Its triangle preserves their
ranks since $T$ has rank zero. Twisting and descent preserve ranks.
Among the remaining lines only the lines $P$ preceding $T$ are
replaced with $P(-D)$, where $D$ is the exceptional curve. For a pair
of surviving lines on opposite sides of $T$, the new degree-two group
is a quotient of the former one, since
\[
 H^2(X,Q\otimes P^\vee)\longrightarrow
 H^2(X,Q\otimes P^\vee(D))
\]
is surjective by the Cartier divisor sequence. The common twist and fully
faithful descent therefore preserve the asserted vanishing.

The resulting birational surface has a line sequence satisfying
\eqref{eq:nonneg}, so Proposition~\ref{prop:normal-pure-two} applies.
For its right orthogonal $\cA=\langle E,F\rangle$ put $W=i^*\cO_x$
at a general point. The normal triangle and $\chi_X(\cO_x,\cO_x)=0$
give $\chi_{\cA}(W,W)=-M<0$. On the other hand
$\chi_{\cA}(W,E)=\chi_X(\cO_x,E)=r$ and
$\chi_{\cA}(W,F)=s$ by local Riemann--Roch or duality at a general
point. The Euler matrix in the exceptional basis $[E],[F]$ then gives
$[W]=r[E]+(s-mr)[F]$ and $\chi_{\cA}(W,W)=r^2+s^2-mrs$,
a contradiction.
For the last assertion move each preceding line past the entire pair
by left-mutating the pair through that line. The equivalence induced
on the pair preserves $m$, and all lines keep their mutual order.
The inequality holds for any resulting ranks if $|m|\le2$.
\end{proof}

\begin{corollary}\label{cor:euler-coherent-pair}
Let $X$ have a full exceptional collection of shifted line bundles and
two consecutive positive-rank coherent sheaves $E,F$. Put
$r=\rk E$, $t=\rk F$, and $e=\chi(E,F)$. Then $X$ is rational if
\[
                       ert-r^2-t^2\le1 .
\]
The same conclusion applies to a consecutive pair of perfect objects
after individual shifts, provided each is generically a vector bundle
in degree zero and $\Ext^{<0}(E,F)=0$ after those shifts.
\end{corollary}

\begin{proof}
Keep the pair fixed and cyclically rotate only the preceding lines.
Assume nonrationality and perform the deformation and contractions of
Theorem~\ref{thm:two-generic-bundles}. The generic fiber profiles stay
concentrated in degree zero. Vanishing of negative forward Ext persists
by semicontinuity and then by the fully faithful mutation and descent.
The Euler number and both ranks are unchanged. On the resulting surface
the line block satisfies \eqref{eq:nonneg}. At a general point $x$,
the covariant module representing $W_x$ in
Proposition~\ref{prop:pair-profile} has both vertex spaces in degree zero,
after removing the common line and shift. In the normalized complex for
its derived endomorphisms there are only arities zero and one. Arity
zero consists of the vertex endomorphisms in degree zero. A degree-$d$
arrow from $E$ to $F$ contributes to arity one in total degree $1-d$ a
space of dimension $rt$; the internal differential is zero in the minimal
model. The differential from the vertex maps lands only in the summands
for $d=0$, because all other arrows have zero action on these vertex
spaces. There is no outgoing differential from arity one. Consequently
for $d>1$ any nonzero $\Ext^d(E,F)$ gives a negative-degree derived
endomorphism of $W_x$. For $d=1$ it gives a degree-zero endomorphism
independent of the scalar identity, which lies in arity zero.
Both cases contradict the negative-degree vanishing and
$\End(W_x)=\C$ given by height at least two. Thus the pair has no
positive forward Ext. Its only possible forward morphisms are in
degree zero, their dimension is $e$, and
Proposition~\ref{prop:pair-profile} proves the assertion.
\end{proof}

\section{Full collections in three orthogonal blocks}
\label{sec:three-object-blocks}

A block is an exceptional collection whose distinct members are completely
orthogonal in both directions. In this section the individual members of
a block may be arbitrary objects of $\Db(X)$. A full exceptional collection
is called formal if the dg category on its members is quasi-equivalent,
with the same objects, to its graded Ext category with zero differential.

\begin{theorem}\label{thm:any-three-blocks}
Let $X$ be a smooth connected projective complex surface.
If $\Db(X)$ has a full exceptional collection in at most three
consecutive blocks, then $X$ is rational. If the collection is formal,
the same conclusion holds for at most four consecutive blocks.
\end{theorem}

\begin{proof}
Suppose that $X$ is nonrational. Additivity of Hochschild homology and
HKR give $q(X)=p_g(X)=0$. We also have $H^0(X,T_X)=0$. Indeed,
$\operatorname{Aut}^0(X)$ fixes the class of a very ample line bundle
since $\Pic^0(X)=0$. Its action on the associated projective embedding
identifies it with a connected linear algebraic group. In characteristic
zero a positive-dimensional connected linear group acting faithfully has
a copy of $\mathbb G_a$ or $\mathbb G_m$ acting nontrivially. The general
orbit of this subgroup is a rational curve. The action on a curve
transverse to a general orbit gives a dominant rational map to $X$
from a ruled surface. A smooth uniruled surface with $q=0$ over
$\mathbb C$ is rational, a contradiction. The Lie algebra of
$\operatorname{Aut}^0(X)$ is $H^0(X,T_X)$, so HKR gives
\[
                 \HH^1(X)=H^1(X,\mathcal O_X)\oplus H^0(X,T_X)=0.
\]

Write $E_1,\ldots,E_N$ for the collection. On these objects take a
strictly unital minimal $A_\infty$ model over
$R=\bigoplus_i\mathbb C e_i$. Its graded category is the Ext category
$A$ with the Yoneda product $m_2$. In a collection of at most three
blocks a composable word of nonidentity classes has length at most two:
every such class goes to a strictly later block. Consequently $m_k$
vanishes on all inputs for $k\ge3$, by composability and strict unitality.
This proves formality in the three-block case. In the four-block
case it holds by hypothesis. Fullness and derived Morita invariance now
identify $\HH^1(A)$ with $\HH^1(X)$.

On the graded algebra $A$ consider the degree-zero derivation
$D(a)=p a$ for $a\in e_j A^p e_i$. As a normalized Hochschild cochain
it has arity one and total degree one. It is a cocycle since $m_2$
adds degrees. In the relative Hochschild complex over the separable
algebra $R$, the differential for this graded algebra increases arity
by one. Since the class of $D$ vanishes, its arity-one component is the
Hochschild differential of an arity-zero, internal-degree-zero cochain
$l\in\bigoplus_i e_i A^0e_i=\bigoplus_i\mathbb C e_i$.
Thus $D=[l,-]$, after choosing the sign of $l$. Writing
$l=\sum_i l_i e_i$, we obtain
\[
       e_j A^p e_i\ne0\quad\Longrightarrow\quad p=l_j-l_i.
\]
On each connected component of the graph of nonzero forward Ext spaces,
subtract the value of $l_i$ at one vertex. The resulting values $s_i$
are integers: along each edge the difference is an integer.
Then $\Ext^q(E_i[s_i],E_j[s_j])=\Ext^{q+s_j-s_i}(E_i,E_j)$ vanishes
for $q\ne0$. We have obtained a full strong exceptional collection with
the same blocks.

Its ordinary endomorphism algebra has normalized relative Hochschild
cochains only in arities less than the number of blocks, since the
radical arrows go strictly forward between blocks. All internal degrees
are zero. Thus $\HH^4(X)=0$. On a surface HKR and Serre duality
identify this group with $H^0(X,\omega_X^{\otimes2})^\vee$. Castelnuovo's
criterion now contradicts nonrationality.
\end{proof}

\begin{remark}\label{rem:formal-any-length}
The grading argument also applies to a formal full exceptional collection
of any length. If no individual shifts put its nonzero forward Ext groups
in degree zero, the Euler derivation gives a nonzero class in
$\HH^1(X)$. Since $q(X)=0$, the action argument in the proof shows that
$X$ is rational.
\end{remark}

\begin{proposition}\label{prop:three-block-square}
Every smooth connected projective complex surface with a full exceptional
collection in at most three blocks has $K_X^2>0$. If the collection is
strong up to individual shifts, then $-K_X$ is nef and big.
The latter conclusion also holds when all its members are shifts of
coherent sheaves.
\end{proposition}

\begin{proof}
Let $N$ be the length, so $N=\rho(X)+2$ by additivity of Hochschild
homology, HKR and the Lefschetz $(1,1)$ theorem. In the real numerical
Grothendieck group the exceptional classes are a basis. For a class
$v$ write $r(v)=\rk v$ and $d(v)=K_X\cdot c_1(v)$. Riemann--Roch gives
\[
       \chi(v,w)-\chi(w,v)=r(w)d(v)-r(v)d(w).
\]
Thus the map $(r,d)$ has rank at most one on the span of each block,
since the Euler pairing is symmetric there. Let $V$ be the direct sum
of the kernels on the individual blocks. Then
$\dim V\ge N-3=\rho(X)-1$.
On $V$, the Euler pairing is positive definite: its diagonal blocks
are restrictions of identity matrices, and pairings across two blocks
vanish by semiorthogonality and the displayed identity. Since $r(v)=0$,
we have $\chi(v,v)=-c_1(v)^2$. The first Chern class thus embeds $V$
as a negative definite subspace of $K_X^\perp$ in the real
N\'eron--Severi group.
If $K_X$ were numerically trivial, the full Euler matrix would be
symmetric and upper unitriangular, hence the identity. This contradicts
the isotropic nonzero point class, which pairs to one with
$[\mathcal O_X]$. The Hodge index theorem now gives $K_X^2>0$.

Suppose next that the collection is strong, after shifts, and let $T$ be
its direct sum. The ordinary algebra $B=\End_X(T)$ has global dimension
at most two. Indeed, its normalized right-module bar resolution over
$\bigoplus_i\mathbb C e_i$ has no terms of length bigger than two.
The tilting equivalence sends $T$ to $B$ and identifies
$S_X^{-1}[2]=(-)\otimes\omega_X^{-1}$ with
$\operatorname{RHom}_B(DB,-)[2]$ for right modules, where $DB$ is the
vector-space dual of $B$. The latter functor sends complexes with
cohomology in degrees at most zero to complexes of the same kind:
a projective resolution of the right module $DB$ has length at most
two. Iterating on the regular module gives
\[
                 \Ext^j_X(T,T\otimes\omega_X^{-m})=0
                 \qquad(j>0,\ m\ge1).
\]
Choose a summand $E$ of nonzero rank, which exists since the classes span
$K_0(X)$. The unit and trace for the perfect complex $E$ split off
$\mathcal O_X$ in its derived endomorphism complex, since their composite
is $r(E)\ne0$. Taking this summand gives
$H^j(X,\omega_X^{-m})=0$ for $j>0$ and $m\ge1$. Consequently
\[
               h^0(X,-mK_X)=1+\tfrac{m(m+1)}2K_X^2 .
\]
Since $K_X^2>0$, the divisor $-K_X$ is big and its volume equals its
self-intersection. For its Zariski decomposition $-K_X=P+Z$ on the
surface, one has $\operatorname{vol}(-K_X)=P^2$,
$P\cdot Z=0$ and $Z^2<0$ unless $Z=0$. The equality of the volume
and the self-intersection forces $Z=0$, proving nefness.

Finally suppose that the members can individually be shifted to
coherent sheaves. A collection with fewer than three blocks can be
refined to three nonempty blocks since $N\ge3$. The standard bounded
$t$-structure is preserved by $S_X[1-3]=(-)\otimes\omega_X$.
Bridgeland--Stern's purity theorem for a full three-block collection
therefore makes the collection strong
\cite[Theorem~4.5 and Corollary~4.6]{BSt10}.
The argument applies to any common heart of a bounded $t$-structure
preserved by tensoring with $\omega_X$.
\end{proof}

\subsection{Properly supported objects between the blocks}
\begin{corollary}\label{cor:three-blocks-supported-strings}
Let $X$ be a smooth connected projective complex surface. Suppose it has
a full exceptional collection of the form
\[
    (\mathcal T_0,\mathcal B_1,\mathcal T_1,
       \mathcal B_2,\mathcal T_2,\mathcal B_3,\mathcal T_3).
\]
Each $\mathcal B_j$ is a consecutive completely orthogonal block of
arbitrary exceptional objects. Each $\mathcal T_j$ is a possibly empty
string of exceptional objects whose supports are proper subsets of $X$.
Then $X$ is rational.
\end{corollary}

\begin{proof}
Suppose $X$ is nonrational. Its minimal model $S$ has $p_g=q=0$ and
$K_S$ nef. In the smooth family of ordered point blowups of $S$ used
in the proof of Proposition~\ref{prop:deform}, the collection extends
termwise and remains full after an \'etale base change and shrinking.
Complete orthogonality persists by semicontinuity. Proper support of a
term persists because the support of the derived restriction to a fiber
is the fiber of the support of its relative perfect extension, and the
dimension of the fibers of this proper support morphism is upper
semicontinuous. As in Proposition~\ref{prop:deform}, choose a resulting
surface $Y$ whose centers are distinct and lie on no rational curve of $S$.

Let $T$ be a term in one of the strings on $Y$. A zero-dimensional
support is impossible since it would give $\chi(T,T)=0$.
The admissible subcategory generated by $T$ has precisely its support.
Every one-dimensional component is a rational curve and at least one
component has negative $K_Y$-degree
\cite[Theorem~1.1]{PirC}.
A rational curve on $Y$ is either a blowup divisor or the strict
transform of a rational curve on $S$ disjoint from all the centers.
The latter has nonnegative canonical degree. Each blowup divisor is
therefore disjoint from all other rational curves under consideration.
Since the support of the indecomposable complex $T$ is connected, it is
one blowup divisor $D$. The classification of admissible subcategories
supported on a $(-1)$-curve
\cite[Proposition~5.7]{PirD}
identifies $T$, up to shift, with $i_*\mathcal O_D(d)$ for some $d$,
where $i$ denotes the inclusion of $D$.

Move $T$ to the front by right mutating each of its predecessors through
it. Every preceding member lies in $T^\perp$, on which the right
mutation is an equivalence onto {}$^\perp T$. Each $\mathcal B_j$ lies
wholly on one side of $T$, so each such block remains completely
orthogonal. Right mutation through a properly supported object preserves
vanishing at the generic point, and therefore preserves proper support
of the other designated string members. A common twist by
$\mathcal O_Y((d+1)D)$ sends $T$ to $i_*\mathcal O_D(-1)$.
After contracting $D$, Orlov's blowup formula identifies the remaining
sequence with a full exceptional collection on the surface with one
fewer center. Generic fibers and orthogonality are preserved under its
fully faithful pullback. Repeat until all designated strings have
disappeared. The resulting surface is birational to $X$ and has a full
collection in at most three completely orthogonal blocks. This
contradicts Theorem~\ref{thm:any-three-blocks}.
\end{proof}

\begin{theorem}\label{thm:four-pure-adjacent}
Let $X$ be a smooth connected projective complex surface with a full
exceptional collection in four consecutive completely orthogonal blocks
$\mathcal B_1,\ldots,\mathcal B_4$. Suppose the individual members can
be shifted so that, for each $i=1,2,3$, all forward Ext spaces from
$\mathcal B_i$ to $\mathcal B_{i+1}$ are concentrated in degree zero.
Then $X$ is rational. In particular a full exceptional collection of four
arbitrary objects suffices if the forward Ext spaces of each of its three
adjacent pairs are concentrated in a single degree. The three degrees may
be different, and no hypothesis is imposed on the nonadjacent pairs.
\end{theorem}

\begin{proof}
Use the shifts in the hypothesis, and take a strictly unital directed
minimal model over the vertex idempotents. Its only possible nonunit
products are $m_2$ and $m_3$. For a homogeneous off-diagonal element
$a$ define a degree-zero linear map
\[
 D(a)=
 \begin{cases}
 (\deg a)a&\text{if }\deg a\notin\{0,-1\},\\
 0&\text{if }\deg a\in\{0,-1\}.
 \end{cases}
\]
Set $D(e_i)=0$ on all vertex units. We claim that $D$ is a total-degree-one
Hochschild cocycle, including when $m_3$ is nonzero. Any composable pair
of strict arrows either consists of two adjacent-block arrows, both of
degree zero, or consists of an adjacent-block arrow and an arrow skipping
one block. In the latter case $m_2$ has the degree of the arrow skipping a
block, so $D$ commutes with this product as a derivation. In the former case
it vanishes on the inputs and the output. The three inputs of every strict
ternary product run through all four blocks and have degree zero.
Its output has degree $-1$, so the commutator with $m_3$ is zero as well.
Products with a unit respect these identities by strict unitality.

If $X$ is nonrational, $\HH^1(X)=0$ as in the proof of
Theorem~\ref{thm:any-three-blocks}. Thus $D$ is a Hochschild coboundary.
Looking at arity one gives $D=[l,-]$ on the graded Ext spaces for a
diagonal element $l=\sum_i l_ie_i$: only the bracket of $m_2$ with the
arity-zero part of a cochain can have arity one. In particular the
differences $l_j-l_i$ are integers along every nonzero Ext space
from $E_i$ to $E_j$. Subtract a constant on each connected component
of the graph of such spaces and shift the individual objects by the
resulting integers. These shifts leave the degrees $0$ and $-1$
unchanged wherever those degrees occurred, and send every other nonzero
degree $t$ to degree $t-(l_j-l_i)=0$. The adjacent-block arrows are
still in degree zero.

In the new model normalized cochains have arity at most three.
Arity-zero cochains are scalar endomorphisms in degree zero.
In arity one their total degree is at most two, since all morphism
degrees are $0$ or $-1$. Of the two inputs in arity two, at most one
can skip a block and hence have negative degree, giving total degree
at most $2+1=3$. In arity three all inputs have degree zero and the
output has degree at most zero. Thus the degree-four cochains vanish
even with the $m_3$ differential. By fullness, $\HH^4(X)=0$, and
Castelnuovo's criterion contradicts nonrationality. For four individual
objects, the three successive differences of shifts can be specified
independently, which proves the last assertion.
\end{proof}

\subsection{Nonformal collections with a prescribed higher product}
\label{sec:higher-weights}
For a directed strictly unital minimal model, call a composable input
\emph{strict} if it goes between distinct members. A word of strict inputs
visits its blocks in increasing order. Write $A^p_{ij}=\Ext^p(E_i,E_j)$.
The following argument takes account of a higher product that need not vanish.

\begin{theorem}\label{thm:fixed-arity}
Suppose that a smooth connected complex projective surface has a full
exceptional collection in $b$ orthogonal blocks. Assume there is a
directed minimal model in which, on words of strict inputs, $m_d$ is
the only possibly nonzero multiplication, for some $d\ge2$.
If $b\le 3d-2$, then $X$ is rational. In particular a full exceptional
collection in four blocks implies rationality whenever all Yoneda
products of two strict inputs vanish, even if the ternary product is nonzero.
\end{theorem}
\begin{proof}
Suppose $X$ is nonrational, so that $\HH^1(X)=0$ by the argument of
Theorem~\ref{thm:any-three-blocks}. Put $w=(d-2)/(d-1)$. Define a
normalized degree-zero endomorphism of the underlying graded spaces
by $D(a)=(p-w)a$ for $a\in A^p_{ij}$, $i<j$, and $D(e_i)=0$.
It represents a Hochschild cochain of total degree one. On $d$ strict
inputs of degrees $p_1,\ldots,p_d$, the output of $m_d$ has degree
$\sum_i p_i+2-d$. Its $D$-weight equals the sum of the input weights
because $2-d+(d-1)w=0$. The binary actions of the identities also
preserve weights. Thus $[m,D]=0$ for the full Hochschild differential.

Since $D$ is a boundary, its arity-one part is the commutator of $m_2$
with an arity-zero degree-zero cochain $\sum_i l_i e_i$. Indeed, $m_1=0$
and a bracket with $m_k$ for $k\ge3$ has arity at least two.
Consequently $p-w=l_j-l_i$ whenever $A^p_{ij}\ne0$, with the sign
of the $l_i$ fixed as in Theorem~\ref{thm:any-three-blocks}.
For a normalized cochain component with $h\ge1$ strict inputs,
input degrees $p_1,\ldots,p_h$ and output degree $q$, telescoping gives
\[
  h+q-\sum_{i=1}^h p_i
       =h+(1-h)w=\frac{h+d-2}{d-1}.
\]
Thus a nonzero component of total degree four requires $h=3d-2$.
But $h\le b-1$, and an arity-zero cochain is diagonal of degree zero.
The degree-four cochain group vanishes. This contradicts
$\HH^4(X)\simeq H^0(X,2K_X)^\vee\ne0$. For four blocks no word has
length greater than three. If the strict binary products vanish, take
$d=3$.
\end{proof}

\begin{remark}\label{rem:four-block-mixed}
There is also a four-block case allowing both binary and ternary
multiplication. Number the four blocks in order. For each pair consisting
of a vertex of block one and a vertex of block four, let $J$ be the
graded subspace of its forward Ext space spanned by the outputs of
all binary products of strict arrows. If the span of the outputs of
$m_3$ for the same endpoints is disjoint from $J$ in some minimal
model, then $X$ is rational. To see this, take a graded complement
to $J$ containing those ternary outputs. Give it weight $-1$, and
give weight zero to all other strict arrow spaces and to $J$.
Replace $p-w$ in the proof by $p$ minus this assigned weight.
Every binary product is compatible with the new derivation. For a
ternary product the inputs have weight zero, its output has weight
$-1$ and its cohomological degree drops by one. Thus the derivation
is a cocycle. A path of length at least two cannot have any input
from block one directly to block four. Its possible cochain degrees
are its length and its length minus one. For a path of length one they
are at most two. Once again no normalized cochain can have degree four.
\end{remark}

\subsection{Canonical degrees of a strong collection}

For a perfect object $E$ on a smooth connected complex projective surface set
$v(E)=(\operatorname{rk}E,K_X\cdot c_1(E))\in\mathbb Z^2$.
Write $[v(E)]\in\mathbb P^1(\mathbb Q)$ for the direction of a nonzero vector.
For a positive-rank torsion-free sheaf put
$\mu_K(E)=(K_X\cdot c_1(E))/\operatorname{rk}E$.

\begin{proposition}\label{prop:canonical-directions}
Suppose that $(E_1,\ldots,E_n)$ is a full strong exceptional collection
of arbitrary perfect objects on $X$. If the nonzero $v(E_i)$ have at
most four distinct directions, then $X$ is rational.
If $K_X$ is nef and the $E_i$ are torsion-free sheaves of positive rank,
a nonrational $X$ would have at least four members that fail slope
semistability for $K_X$.
\end{proposition}
\begin{proof}
Write $v_i=(r_i,d_i)$. Surface Riemann--Roch and backward
semiorthogonality give, for $i<j$,
\[
 h_{ij}:=\dim\operatorname{Hom}(E_i,E_j)=\chi(E_i,E_j)
 =\chi(E_i,E_j)-\chi(E_j,E_i)=d_i r_j-r_i d_j .
\]
In particular all these signed determinants are nonnegative.
The sum of the collection is tilting. For its ordinary endomorphism
algebra, compute Hochschild cohomology by normalized bar cochains
relative to the separable algebra of vertex idempotents.
A nonzero degree-four cochain requires five vertices
$i_0<\cdots<i_4$ with $h_{i_j i_{j+1}}>0$ for each $j$ and
$h_{i_0 i_4}>0$. In the oriented basis $v_{i_0},v_{i_4}$, the
nonnegative pairings with both endpoints put all three intermediate
vectors in the pointed cone between the endpoints. Their arguments
in this cone are weakly ordered by the nonnegative pairwise signed
determinants, and every consecutive argument increases strictly.
Consequently the five vectors have five distinct projective directions.
With at most four directions the fourth cochain group thus vanishes.
Derived Morita invariance and HKR give
$0=\operatorname{HH}^4(X)=H^2(X,\bigwedge^2T_X)
=H^0(X,2K_X)^*$. Moreover $q(X)=0$ by Hochschild homology and
fullness. Castelnuovo's criterion proves rationality.

For the final assertion, $h_{ij}=r_ir_j(\mu_K(E_i)-\mu_K(E_j))$.
The slopes are therefore weakly decreasing, with a nonzero Hom space
for every strict decrease. Two slope-semistable sheaves of different
slopes cannot both occur, since there is no nonzero morphism from a
larger slope to a smaller one. All the semistable terms therefore have
one common slope if they occur. At least five directions are required
on a nonrational surface, so at least four objects are unstable.
\end{proof}

\begin{remark}\label{rem:four-directions-formal}
The four-direction criterion also holds for a formal full collection of
arbitrary perfect objects. Indeed, on a hypothetical nonrational surface
the Euler derivation in Theorem~\ref{thm:any-three-blocks} makes the
collection strong after shifts. A shift changes $v_i$ by a sign, preserving
its projective direction.
\end{remark}

\begin{remark}\label{rem:semistable-bundles}
The point coefficient also excludes some collections of higher-rank
bundles without strongness. Suppose $K_X$ is nef and each $E_i$ is a
slope-semistable bundle. If, for each $i<j$ with a nonzero forward
Hom, the group $\operatorname{Ext}^1(E_i,E_j)$ vanishes, the
collection cannot be full. Indeed, for such a pair semistability
would give $\mu_K(E_i)\le\mu_K(E_j)$, and thus
$\chi(E_i,E_j)\le0$ by the skew Euler formula above. On the other hand
the alleged vanishing implies
$\chi(E_i,E_j)=h_{ij}+\dim\operatorname{Ext}^2(E_i,E_j)>0$.
There are therefore no forward Hom spaces. Each arrow on a normalized
point-bar path then has degree at least one, and the terminal
factor has degree two. Its arity cancels at most one degree per
arrow, giving point height at least two and excluding fullness.
This applies in particular if each forward pair has nonzero Ext
in at most one degree.
\end{remark}

\section{Hereditary components and exceptional generators}\label{sec:hereditary-new}
\subsection{A hereditary component and an orthogonal exceptional block}
\begin{theorem}\label{thm:gluing-hereditary}
Let $X$ be a smooth connected projective complex surface with a
semiorthogonal decomposition $\Db(X)=\langle\cA,\cB\rangle$.
Suppose $\cA$ is the perfect category of the path algebra of a finite
acyclic quiver, and $\cB$ is generated by a completely orthogonal
exceptional collection. Then $X$ is rational.
\end{theorem}
\begin{proof}
The two components give a full exceptional collection. Additivity of
Hochschild homology and HKR imply $q(X)=p_g(X)=0$.
Assume $X$ is nonrational. The automorphism argument in Theorem~\ref{thm:any-three-blocks} gives $\HH^1(X)=0$.

Write $H=\C Q$ and $R=\C^t$ for the two Morita models.
With the convention that the mixed morphisms go from $\cA$ to $\cB$,
let $M$ be their gluing bimodule. Its components $M_{aj}$, indexed by
the connected components $Q_a$ of $Q$ and the primitive factors of
$R$, are objects of $\Db(\operatorname{rep} Q_a)$ up to the opposite
module convention. The Hochschild cochain gluing triangle gives
\[
 C^*(X)\longrightarrow C^*(H)\oplus C^*(R)
    \longrightarrow \RHom_{\mathrm{bimod}}(M,M)
    \longrightarrow C^*(X)[1].
\]
This is the gluing sequence of \cite[Theorem~7.7 and Remark~7.8]{Kuz09}.
Equivalently restrict normalized cochains of the upper triangular dg
model to the two diagonal categories. A word in the kernel crosses
from $H$ to $R$ exactly once; removing that distinguished input gives
the two-sided module bar complex for the derived endomorphisms of $M$
with its cochain degree lowered by one. The connecting map from degree
zero centers to $\End(M)$ is the difference of their two actions.

Since $\HH^1(X)=0$, this difference map is surjective. The center of
each connected acyclic path algebra is $\C$: commuting with the vertex
idempotents leaves only scalar multiples of the trivial paths at the
vertices, and commuting with each arrow equates the scalars at its
endpoints. On the $(a,j)$ summand the two central actions are therefore
both scalar. Since the bimodule category is the product of the
corresponding derived module categories, surjectivity implies
$\End(M_{aj})=\C$ whenever $M_{aj}\ne0$. Every such summand is
indecomposable. Heredity identifies it with a single shifted module.
Hence $\Ext^3_{\mathrm{bimod}}(M,M)=0$.

The Hochschild cohomology of an acyclic path algebra vanishes in degrees
at least two. Indeed, writing $S$ for the product of its vertex fields
and $V$ for its arrow bimodule, its bimodule resolution is
\[
 0\to H\otimes_S V\otimes_S H\to H\otimes_S H\to H\to0.
\]
The semisimple algebra $R$ has zero positive Hochschild cohomology.
The gluing sequence now gives $\HH^4(X)=0$. The top HKR term and Serre
duality identify its dual with $H^0(X,\omega_X^2)$. Since $q=0$,
Castelnuovo's rationality criterion gives the contradiction.
\end{proof}

In particular, Theorem~\ref{thm:gluing-hereditary} applies to a full
exceptional collection consisting of one arbitrary exceptional object and
a consecutive complement of arbitrary finite acyclic quiver type, including
wild type. Rotating either component by the Serre functor interchanges
their positions and preserves its equivalence type.

\begin{corollary}\label{cor:gluing-hereditary-heart}
Theorem~\ref{thm:gluing-hereditary} still holds if $\cA$ has a full
exceptional collection and is equivalent to the bounded derived category
of an Ext-finite hereditary abelian category $\mathcal H$.
\end{corollary}
\begin{proof}
First $\cA$ has a tilting object $T$ lying in $\mathcal H$.
Indeed every exceptional generator can be shifted into $\mathcal H$ by
heredity and indecomposability. Induct on the length of the resulting
exceptional sequence. If $T_0$ is a rigid generator for the preceding
terms, let $E$ be the next exceptional object and take the universal
extension $0\to E^{\oplus d}\to T_1\to T_0\to0$, where
$d=\dim\Ext^1(T_0,E)$. The connecting map
$\Hom(E^{\oplus d},E)\to\Ext^1(T_0,E)$ is surjective, so
$\Ext^1(T_1,E)=0$. Backward semiorthogonality and exceptionality of
$E$ give $\Ext^1(E,T_1)=0$. Applying $\Hom(-,T_0)$ also gives
$\Ext^1(T_1,T_0)=0$, and then applying $\Hom(T_1,-)$ gives
$\Ext^1(T_1,T_1)=0$. The object $T_1\oplus E$ generates the enlarged
subcategory. This proves the induction.

Put $H=\End(T)$. The tilting equivalence identifies $\cA$ with
$\Perf(H)$. Since $\cA$ is an admissible subcategory of a smooth
proper category, $H$ is homologically smooth. Its perfect bimodule
resolution of the diagonal implies that every finite-dimensional
module is perfect, by tensoring that resolution with the module.
Under the inverse tilting equivalence each simple $H$-module goes to
$M[s]$ with $M\in\mathcal H$, because it is indecomposable.
Since $\RHom(T,M[s])$ is a nonzero module in degree zero and
$\mathcal H$ is hereditary, $s\in\{0,1\}$.
Thus the Ext groups in degrees $n\ge3$ between simple $H$-modules
vanish, and $H$ has global dimension at most two.

The diagonal of $H$ also has projective dimension at most two as an
$H$-bimodule. Over $\C$ the simple bimodules are
$\Hom_\C(D_i,D_j)$ for simple modules $D_i,D_j$, and bar adjunction
identifies $\Ext^n_{H^e}(H,\Hom_\C(D_i,D_j))$ with
$\Ext^n_H(D_i,D_j)$. Filtration by simple bimodules gives the bound.
In particular $\HH^4(H)=0$.

Decompose $H$ into indecomposable algebra factors $H_a$.
The center of each is a finite-dimensional commutative local algebra,
since a nontrivial central idempotent would split the factor.
In the gluing argument the combined central action on
$\End(M_{aj})$ is the action of this local center: the other term acts
by scalars, which are already in its image. Surjectivity therefore
implies that $\End(M_{aj})$ is local for each nonzero summand.
It follows that $M_{aj}$ is indecomposable and hence is a shifted
object of $\mathcal H$. For the opposite module convention the same
assertion follows by duality from the opposite hereditary category.
Its third self-Ext vanishes, so the conclusion of
Theorem~\ref{thm:gluing-hereditary} follows from the same gluing sequence.
\end{proof}

\subsection{Two hereditary components of arbitrary representation type}
\begin{theorem}\label{thm:two-hereditary-new}
Let $X$ be a smooth connected complex projective surface. Suppose that
$\Db(X)=\langle\cA,\cB\rangle$, where each component is equivalent
to the bounded derived category of representations of a finite acyclic
quiver. Fix such an equivalence on $\cB$. Assume that the images of its
simple modules restrict on a dense open subset of $X$ to vector bundles
of positive rank in a common degree. Then $X$ is rational.
\end{theorem}
\begin{proof}
Remove the common shift, denote the images of the simple modules by
$S_j$, and put $r_j=\operatorname{rk}S_j>0$. Denote the indecomposable injectives
in their quiver heart by $J_j$, indexed so that
$\dim\Hom(S_l,J_j)=\delta_{lj}$.
The objects $J_j$ admit filtrations with factors among the $S_l$.
After restriction to the chosen open set, the successive extensions
are vector bundles in degree zero. Fix a point $x$ in this open set.

Let $Q_\cB$ be the right projection and
$W=i^*\cO_x$ the left projection to $\cA$. We first show that
\[
          Q_\cB\cO_x\simeq I:=\bigoplus_j J_j^{\oplus r_j}.
\]
Indeed, $\RHom_{\cB}(S_j,Q_\cB\cO_x)
=\RHom_X(S_j,\cO_x)$ is $\C^{r_j}$ in degree zero.
Every indecomposable summand of $Q_\cB\cO_x$ is $M[t]$ for an
indecomposable module $M$, since the quiver heart is hereditary.
Any nonzero $M$ has a simple submodule, forcing $t=0$.
The equality $\Ext^1(S_j,M)=0$ for every $j$ then gives injectivity
of $M$ by induction on the length of the first argument.
The dimensions of the socles give the indicated multiplicities.

In the projection triangle $I\to\cO_x\to W$, the object $I$ is
locally free near $x$. Thus $\RHom_X(\cO_x,I)$ is concentrated
in degree two by local duality. Applying $\RHom_X(\cO_x,-)$,
using left adjunction on $W$, gives
\[
 \End_{\cA}(W)=\C,\qquad
 T_xX\lhook\joinrel\longrightarrow \Ext^1_{\cA}(W,W).
\]
Thus $W=N[s]$ with $N$ an indecomposable nonprojective module of the
other hereditary heart. For each vertex projective $P$ in this heart
one has $\Hom(N,P)=0$. The image of any nonzero map is a nonzero
projective submodule of $P$, hence would split off from $N$.
Consequently $\RHom_X(\cO_x,P)=\RHom_{\cA}(W,P)$ is concentrated
in degree $s+1$. Serre duality, using $S_X(\cO_x)\simeq\cO_x[2]$,
implies that $\RHom_X(P,\cO_x)$ is concentrated in degree $1-s$.
Also $\RHom_X(P,W)$ is concentrated in degree $-s$ by projectivity.
The triangle now gives
\[
   \Ext_X^n(P,I)=0\quad(n\ne 1-s).
\]
This holds for every vertex projective $P$ of $\cA$. Since the
multiplicities $r_j$ are positive, it holds with $I$ replaced
by any $J_j$.

Let $P_0$ and $I_0$ be the sums of the vertex projectives of $\cA$
and the vertex injectives of $\cB$, respectively. By semiorthogonality,
$T=P_0[s-1]\oplus I_0$ is a tilting object generating $\Db(X)$.
Its ordinary endomorphism algebra $D$ is triangular, with hereditary
diagonal blocks. Its global dimension is at most three. To see this,
choose the convention for left modules in which a module over $D$ is a
triple $(U,V,M\otimes_H U\to V)$, for the two hereditary algebras
$H,K$ and the $(K,H)$-bimodule $M$. The modules $(0,V)$ have projective
dimension at most one. For projective $P$ over $H$, the natural
surjection $(P,M\otimes_H P)\to(P,0)$ has kernel $(0,M\otimes_H P)$,
so $(P,0)$ has projective dimension at most two. A projective
resolution of length one over $H$ now gives projective dimension at
most three for $(U,0)$. Every triple is an extension of $(U,0)$ by
$(0,V)$, proving the bound.

The tilting equivalence sends $T$ to the regular module and $S_XT$ to
its vector-space dual $D^*$. Hence
\[
 0=\Ext_D^4(D^*,D)
 =\Hom_X(S_XT,T[4])
 =\Hom_X(T,T\otimes\omega_X^2)^*.
\]
A nonzero bicanonical section induces a nonzero morphism in the last
Hom space: at the generic point it acts invertibly on $T$, whose generic
fiber is nonzero since $T$ generates. This gives $P_2(X)=0$. The two
quiver components give a full exceptional collection on $X$, so
Hochschild homology and HKR imply $q(X)=0$. Castelnuovo's criterion
proves the assertion.
\end{proof}
\begin{corollary}\label{cor:two-hereditary-exceptional-new}
Suppose a full exceptional collection consists of two consecutive
blocks. Suppose the first generates the derived category of a finite
acyclic quiver. In the second block let every member be generically
a nonzero vector bundle in a common degree and let the forward
morphism spaces be concentrated in degree one. Then $X$ is rational.
In particular the second block may consist of vector bundles of
arbitrary positive ranks.
\end{corollary}
\begin{proof}
For the second block take the quiver with one vertex for each member
and with arrows chosen so that the Ext-quiver of its simple modules is
the prescribed degree-one quiver. It is acyclic by semiorthogonality.
In the strictly unital minimal $A_\infty$ model, each operation on $m\ge2$
nonidentity inputs would land in degree $m+(2-m)=2$, which is zero.
The simple modules of the path algebra have the identical minimal
model by heredity. Equivalence on the generators therefore induces an
equivalence of the two generated perfect categories, identifying the
second block with the simple modules. Apply Theorem~\ref{thm:two-hereditary-new}.
\end{proof}

\subsection{Unirationality from the two projections of skyscraper sheaves}
\begin{theorem}\label{thm:discrete-projections}
Let $X$ be a smooth connected projective variety over $\C$ with
$\Db(X)=\langle\cA,\cB\rangle$. Suppose $\cA$ has at most countably
many isomorphism classes of objects. Suppose the indecomposables of
$\cB$ are exhausted by countably many isomorphism classes and
countably many perfect algebraic families parametrized by open subsets
of $\mathbb P^1$. Then $X$ is unirational. In dimension two $X$ is rational.
The roles of the two components may be interchanged.
\end{theorem}
\begin{proof}
For each closed $x$ take its projection triangle
$B_x\to\cO_x\to A_x\to B_x[1]$. The possible $A_x$ belong to a
countable list. Write $B_x$ as a finite direct sum of indecomposable
objects, using Krull--Schmidt. There are countably many possible lists
of families containing these summands, where a fixed object is regarded
as a family over a point.

Fix $A$ and a list of $m$ families $B_{t_i}$, parametrized by curves
or points $T_i$. The complexes on $T_i$ computing
$\RHom_X(A,B_{t_i}[1])$ are perfect by properness and smoothness.
On each curve, generic base change and further restriction give an
open subset where their degree-zero cohomology is a trivial bundle
and commutes with passage to the fiber. The remaining points of the
curve may be treated as separate strata. We may take the open strata
affine. A basis of the degree-zero cohomology lifts to actual relative
morphisms there, since higher cohomology of the cohomology bundles
vanishes on an affine curve. Thus, for each product of strata, all
morphisms
\[
        A\longrightarrow \bigoplus_{i=1}^m B_{t_i}[1]
\]
form a trivial vector bundle over the product, with a universal
morphism on the parameter space times $X$. The parameter space is
rational. Taking the shifted cone gives a perfect family containing
the middle terms of all triangles with these data. If the variable
summands occur on the source side of the connecting morphism, use
$\RHom_X(B_{t_i},A[1])$ instead, by the same argument.

The locus in any such family where the fibers are coherent sheaves
in degree zero and flat over the parameter is open. On any chart
containing a skyscraper, the Hilbert polynomial is identically one
since it is constant for a connected perfect family. Thus this open
locus parametrizes length-one skyscrapers and carries a morphism to
$\operatorname{Hilb}^1(X)=X$. Explicitly the pushforward of the flat
length-one family is a line bundle; tensoring with its inverse gives a
quotient of the structure sheaf that defines the corresponding family
of length-one subschemes. This gives countably many constructible images
of rational parameter spaces covering the closed points of $X$. Over
$\C$ a variety cannot be a countable union of proper closed subvarieties,
so one of the images is dense. Therefore $X$ is unirational, and a
smooth projective unirational complex surface is rational.
\end{proof}
\begin{remark}\label{rem:discrete-projections}
The proof uses the countability and rational-family hypotheses only for
the objects occurring in the two projections of skyscraper sheaves.
If both components have only countably many isomorphism classes, the
countable Rouquier dimension of the whole category is at most one, so in
dimension at least two such a decomposition is excluded by
\cite[Proposition~9]{Ola23}. A component equivalent to the derived
category of an extended Dynkin quiver supplies the required families on
the variable side, through its universal quiver representations and the
classification of its indecomposable objects.
\end{remark}

\subsection{Exceptional generators with square-zero radical}
\label{sec:radical-square}

\begin{theorem}\label{thm:radical-square}
Let $X$ be a smooth connected projective complex surface. Let
$\mathcal A=\mathrm{Perf}(X)$, or let $\mathcal A=\langle G\rangle^\perp$
where $G$ is exceptional and its restriction to a dense open
subset of $X$ is a shifted line bundle. Suppose that $\mathcal A$
is equivalent, with its enhancement, to $\mathrm{Perf}(A)$ for a
finite-dimensional basic complex algebra of finite global dimension
whose radical $J$ satisfies $J^2=0$. Then $X$ is rational.
\end{theorem}

\begin{proof}
First we construct a family $W_x\in\mathcal A$ on a dense open
$U\subset X$ whose members at distinct closed points are nonisomorphic
and satisfy $\dim\operatorname{Ext}^1(W_x,W_x)=2$. For
$\mathcal A=\mathrm{Perf}(X)$ use $W_x=k_x$.
In the second case shift $G$ so that $G|_U$ is a line bundle and let
$W_x$ be the cone of the evaluation
$Q_x=\operatorname{RHom}(G,k_x)\otimes G\to k_x$. This is an
algebraic family, obtained by relative evaluation against
$\mathcal O_\Delta$ on $U\times X$. For $x,y\in U$,
$\operatorname{RHom}(k_x,Q_y)$ is one-dimensional in degree two by
local duality. Adjunction gives the triangle
\[
 \operatorname{RHom}(k_x,Q_y)\longrightarrow
 \operatorname{RHom}(k_x,k_y)\longrightarrow
 \operatorname{RHom}_{\mathcal A}(W_x,W_y).
\]
It gives $\operatorname{Hom}(W_x,W_y)=0$ for $x\ne y$ and
$\operatorname{End}(W_x)=\mathbb C$. For $x=y$ the first arrow
on degree two is an isomorphism. Its Serre-dual map sends the
identity, after trivializing $\omega_X|_x$, to the evaluation of
$G$ at $x$, which is nonzero. Hence
$\operatorname{Ext}^1(W_x,W_x)=T_xX$, as asserted.

We describe the family by complexes of projective right $A$-modules.
If $T$ is the image of the regular module under
the given equivalence and $\mathcal W$ the family, take the relative
complex
$\mathbf R p_{1*}\mathbf R\mathcal H\!om(p_2^*T,\mathcal W)$.
It is a complex of right modules over the dg endomorphism algebra
of $T$, which is quasi-isomorphic to $A$ by truncation since its
cohomology is concentrated in degree zero. It is perfect over $U$,
and derived base change gives the modules of the objects $W_x$.
Over $K=\mathbb C(U)$ choose a bounded minimal complex of projective
$A_K$-modules for its generic fiber. Such a complex is obtained
from a bounded projective representative by splitting off
contractible summands until each differential has image in the
radical. The algebra $A$ is split basic, so its terms are
$P^i\otimes_{\mathbb C}K$ for fixed finite sums of vertex
projectives $P^i$. Shrink to an affine open so that the matrix entries
of the differentials are regular. The generic isomorphism with the
given family also spreads on shrinking: the source is a bounded
complex of direct summands of finite free $A$-modules, so morphisms
out of it commute with localization. Its cone has bounded coherent
cohomology over the base and vanishes generically, hence vanishes
on a further open. We obtain a morphism from $U$ to the affine space
\[
 Z=\bigoplus_i\operatorname{Hom}_A(P^i,P^{i+1}J).
\]
Every point of $Z$ is a differential on the specified graded
projective module, since $J^2=0$.

The group $H=\prod_i\operatorname{Aut}_A(P^i)$ acts on $Z$ by change
of bases. It is irreducible as a product of open subsets of vector
spaces. Two minimal complexes of these terms are isomorphic in the
derived category exactly if they belong to the same orbit. Indeed, a
homotopy equivalence between bounded projective complexes reduces
modulo $J$ to an isomorphism on each term, since the differentials
then vanish. Nakayama's lemma makes each component of the
homotopy equivalence an isomorphism. Thus different $x\in U$
give different orbits.

At $d\in Z$ the tangent space is $Z$ itself. For each $b\in Z$ the
identity $db+bd=0$ shows that $b$ is a degree-one endomorphism
cocycle of $(P^\bullet,d)$. The coboundaries $da-ad$ for all
degree-zero $a$ form exactly the orbit tangent space by
differentiating the change-of-bases action. Over $\mathbb C$ group
orbits are smooth, giving an injection
\[
 T_d Z/T_d(Hd)\lhook\joinrel\longrightarrow
 \operatorname{Ext}^1_A((P^\bullet,d),(P^\bullet,d)).
\]
Every orbit belonging to our family therefore has codimension at
most two in $Z$. Let $Y$ be the closure in $Z\times U$ of the image
of $(h,x)\mapsto(hd_x,x)$. Its constructible image has dimension
at least $\dim Z$, since the fiber over each $x$ is its orbit.
Projection of the constructible image to $Z$ is injective on
closed points by disjointness of the orbits. The dimension theorem
therefore gives equality and density in $Z$. The irreducible
variety $Y$ contains a dense open lying in that image. Its
complement has smaller dimension and cannot dominate $Z$, so
$Y\to Z$ is generically one-to-one. In characteristic zero it is
birational. Composing its rational inverse with $Y\to U$ gives a
dominant rational map from affine space to $X$. A smooth projective
unirational complex surface is rational.
\end{proof}

\begin{remark}\label{rem:radical-square}
The hypothesis of Theorem~\ref{thm:radical-square} holds for the
endomorphism algebra of a full strong exceptional collection in which the
composite of any two composable off-diagonal morphisms vanishes, with no
bound on the number of objects, and also when such a collection forms the
complement of $G$. The algebra may have wild representation type.
\end{remark}

\begin{corollary}\label{cor:square-zero-impossible}
If $X$ is a smooth connected projective complex variety of dimension
$n\ge2$, then $\mathrm{Perf}(X)$ itself admits no enhanced equivalence
to $\mathrm{Perf}(A)$ for an algebra $A$ as in
Theorem~\ref{thm:radical-square}. In particular the case
$\mathcal A=\mathrm{Perf}(X)$ in the theorem is impossible.
\end{corollary}

\begin{proof}
Use the same construction with the $n$-dimensional family $W_x=k_x$.
The orbit codimension in $Z$ is at most
$\dim\operatorname{Ext}^1(k_x,k_x)=n$. On shrinking $U$ the orbit
dimension is constant. Since the orbits of distinct points are
disjoint, their union in $Z$ has dimension $\dim U+\dim Hd_x$.
Hence the codimension is at least $n$, and equality makes the
injection from the orbit tangent quotient to
$\operatorname{Ext}^1(k_x,k_x)$ an isomorphism. All its vectors
are represented by degree-one endomorphisms with image in the
radical. The composite of any two such representatives is zero,
since $J^2=0$. This contradicts the nonzero product in the
Koszul self-Ext algebra $\bigwedge^*T_xX$ in degree two.
\end{proof}

\section{Applications to products with the projective line}\label{sec:products}

Let $S$ be a smooth connected projective complex surface. Every line bundle on $S\times\PP^1$ is uniquely of the form $M\boxtimes\cO_{\PP^1}(a)$. We first describe how the coefficient comparison behaves under external products, and then use mutations to reduce collections with certain repetition patterns to Theorem~\ref{thm:rational}.

\subsection{Transport of the simultaneous comparison}

\begin{proposition}[Product transport of coefficient matrices]\label{prop:transport}
Let $S,T$ be smooth projective complex varieties, let $\cB_S\subset\Db(S)$ be admissible, and put $\cB=\cB_S\boxtimes\Db(T)\subset\Db(S\times T)$, the thick subcategory generated by the external products. For probes $A_a\in\Perf(S)$ and $U_a\in\Perf(T)$, write $Z_a=A_a\boxtimes U_a$ and $D_{ab}=\RHom_T(U_a,U_b)$. The comparison triangle of the coefficient $Z_a^\vee\boxtimes S_{S\times T}(Z_b)$ is canonically
\begin{equation}\label{eq:transport}
  N^S_{ab}\otimes D_{ab}\longrightarrow C^S_{ab}\otimes D_{ab}\longrightarrow R^S_{ab}\otimes D_{ab},
\end{equation}
where the superscript $S$ denotes the comparison for $A_a^\vee\boxtimes S_S(A_b)$ and $\cB_S\subset\Db(S)$, compatibly with the comparison maps and with composition. If $D_{ab}\ne0$, then $h^{S\times T}_{ab}=h^S_{ab}+\inf D_{ab}$. In particular, for a fixed $t\in T$ and $Z_a=\cO_{s_a}\boxtimes\cO_t$,
\begin{equation}\label{eq:transport-points}
  h^{S\times T}_{ab}=h^S_{ab},
\end{equation}
including infinite heights.
\end{proposition}

\begin{proof}
If $Q_S$ is the projector onto $\cB_S$, the kernel $Q_S\boxtimes\cO_{\Delta_T}$ on $S^2\times T^2\simeq(S\times T)^2$ gives the projector onto $\cB$: one has $Q_\cB(A\boxtimes U)\simeq Q_S(A)\boxtimes U$, with counit the external product of the counit on $S$ with the identity, and the projection property follows on external products and hence on their thick span. Theorem~\ref{thm:external} and the K\"unneth formula give
\[
  \RHom_{S\times T}(Z_a,Q_\cB Z_b)\simeq\RHom_S(A_a,Q_SA_b)\otimes D_{ab},
\]
and the same formula without the projectors; naturality identifies the maps and their cones. Over a field the initial degrees of nonzero bounded complexes add. For point probes, $\Ext^q_T(\cO_t,\cO_t)=\bigwedge^qT_tT$ has initial degree zero.
\end{proof}

\begin{corollary}\label{cor:product-nonfull}
Let $\bE_S$ be a nonempty exceptional collection of line bundles on a smooth projective complex surface satisfying \eqref{eq:nonneg}. Its external product with the Beilinson collection on $\PP^{r-1}$, $r\ge2$, is not full on $S\times\PP^{r-1}$: for every $s\in S^\circ_{\bE_S}$ and every $t\in\PP^{r-1}$,
\[
  h_{P_{(s,t)}}\bigl(\langle\bE_S\rangle\boxtimes\Db(\PP^{r-1})\bigr)=2>0=\alpha_{P_{(s,t)}}.
\]
\end{corollary}

\begin{proof}
The product collection generates $\langle\bE_S\rangle\boxtimes\Db(\PP^{r-1})$ by Beilinson's theorem \cite{Bei78}. Apply Theorem~\ref{thm:point-height}, Proposition~\ref{prop:transport} and Theorem~\ref{thm:intro-localization}.
\end{proof}

\subsection{Two consecutive relative degrees}

\begin{proposition}[A relative-degree window]\label{prop:window}
Let $X=S\times\PP^1$ and let
\[
  \bE=\bigl(M_i\boxtimes\cO_{\PP^1}(a_i)\bigr)_{i=1}^N,\qquad N>0,\qquad a\le a_i\le a+1,
\]
be exceptional, with $a=\min_ia_i$. Let $\bE_0$ be the collection of the $M_i$ with $a_i=a$, in the induced order, and put $\cB=\langle\bE\rangle$, $\cB_0=\langle\bE_0\rangle$. For $Z_s=\cO_s\boxtimes\cO_{\PP^1}(a)$, the comparison of each ordered coefficient $Z_s^\vee\boxtimes S_X(Z_{s'})$ is naturally isomorphic to the comparison on $S$ for $\cO_s^\vee\boxtimes S_S(\cO_{s'})$ and $\cB_0$. In particular, if $K_S$ is nef, then for $s$ in a nonempty open subset of $S$
\[
  h_{Z_s^\vee\boxtimes S_X(Z_s)}(\cB,X)=2,\qquad\alpha_{Z_s^\vee\boxtimes S_X(Z_s)}=0 .
\]
\end{proposition}

\begin{proof}
If two adjacent members have degrees $a+1$ and $a$ in this order, the Hom complex in the reverse direction vanishes by exceptionality and the one in the forward direction by $R\Gamma(\PP^1,\cO(-1))=0$. Exchanging such orthogonal pairs orders the collection by degree. Full faithfulness of $F_a(A)=A\boxtimes\cO(a)$ makes $\bE_0$ exceptional. For the projectors $Q_0$ onto $\cB_0$ and $Q_\cB$ onto $\cB$ we have, compatibly with the counits,
\begin{equation}\label{eq:window-projector}
  Q_\cB F_a(A)\simeq F_aQ_0(A).
\end{equation}
Indeed, apply $F_a$ to the projection triangle of $A$: its first term belongs to $\cB$, and its cone is right orthogonal to the members of degree $a$ by full faithfulness and to those of degree $a+1$ by $R\Gamma(\cO(-1))=0$. Consequently $\RHom_X(Z_s,Q_\cB Z_{s'})\simeq\RHom_S(\cO_s,Q_0\cO_{s'})$, and full faithfulness also identifies the ambient complexes. Naturality in both probes identifies the comparison matrices, and Theorem~\ref{thm:point-height} gives the diagonal height.
\end{proof}

\subsection{Descent from a product with the projective line}\label{sec:parity-descent}

\begin{theorem}\label{thm:parity-descent}
Let $S$ be a smooth projective complex variety. Suppose
$(M_i\boxtimes\mathcal O_{\mathbb P^1}(a_i))_{i=1}^N$ is a full exceptional collection on $S\times\mathbb P^1$, where $M_i\in\mathrm{Perf}(S)$.
For each $e\in\mathbb Z/2$, the subsequence $(M_i)_{a_i\equiv e\bmod 2}$ in its original order is a full exceptional collection on $S$.
In particular, if $S$ is a surface and the collection on the product consists of line bundles, then $S$ is rational.
\end{theorem}

The distinction between the two parities can be realized over $\mathbb C$ using a nonsemisimple Hopf algebra. We give the equivariant argument in the Kronecker quiver presentation of $\mathrm{Perf}(\mathbb P^1)$.

\begin{proof}
Write $k=\mathbb C$. A common twist of even degree on $\mathbb P^1$ allows all $a_i$ to be positive without changing either parity. The K\"unneth formula gives exceptionality of each $M_i$. If $i<j$ and $a_j\equiv a_i\bmod 2$, then $a_i-a_j\ne-1$, so $R\Gamma(\mathcal O(a_i-a_j))\ne0$. Backward semiorthogonality on the product thus gives $\operatorname{RHom}(M_j,M_i)=0$. It remains to prove fullness for each parity.

Let $H$ be the four-dimensional Sweedler Hopf algebra, generated by $g,u$ with $g^2=1$, $u^2=0$, $gu=-ug$, with coproducts $\Delta g=g\otimes g$ and $\Delta u=u\otimes 1+g\otimes u$. Following Walton and Witherspoon \cite[Example~4.16]{WalWit14}, make the polynomial ring $B=k[x,y]$ an $H$-module algebra by putting
\[
 g(x)=x,\quad g(y)=-y,\quad u(x)=0,\quad u(y)=x,
 \qquad u(fh)=u(f)h+g(f)u(h).
\]
Indeed the rule preserves the relation $xy=yx$ and
\[
 u(x^{n-j}y^j)=
 \begin{cases}x^{n-j+1}y^{j-1},&j\text{ odd},\\0,&j\text{ even}.\end{cases}
\]
Each graded part $B_n$ is an $H$-module. For odd $n\ge1$, successive pairs of monomials indexed by $j=2l+1,2l$ span copies of the projective left module $H(1-g)/2$. In particular $B_n$ is projective. In contrast the trivial module $k$ has infinite projective dimension. Indeed $H(1\pm g)/2$ has the two one-dimensional characters as top and socle in opposite order. Their nonsplit length-two projective covers give a periodic resolution of $k$. A projective $H$-module tensored on the left of any finite dimensional $H$-module remains projective for the diagonal action, since tensoring with a finite dimensional dual is an exact right adjoint. Thus perfect $H$-complexes remain perfect when tensored on the right with bounded complexes of finite dimensional $H$-modules.

Here is an $H$-enrichment of the relevant derived Hom complexes. Take the Kronecker path algebra $A$ with arrow space $B_1$ and trivial action on its two primitive idempotents. There are no paths of length two, so this defines an $H$-module algebra. Use right $A$-modules, with an arrow action $U_1\otimes B_1\to U_0$. For $a\ge1$ let $D(a)$ denote the representation
\[
 (D(a)_1,D(a)_0)=(B_{a-1},B_a),\qquad
 B_{a-1}\otimes B_1\longrightarrow B_a ,
\]
whose arrow action is multiplication. It is an equivariant right $A$-module. Under Beilinson's equivalence \cite{Bei78} defined by $\mathcal O\oplus\mathcal O(1)$ on $\mathbb P^1$, it is the image of $\mathcal O(a)$. For $t\ge a$, left polynomial multiplication gives an equivariant isomorphism
\begin{equation}\label{eq:odd-projectives}
 B_{t-a}\simeq \operatorname{RHom}_A(D(a),D(t))
          \quad\text{in }D(H\text{-}\mathrm{Mod}).
\end{equation}
For the equivariance, multiplication $B_{t-a}\otimes D(a)\to D(t)$ is $H$-linear and right $A$-linear. After forgetting $H$, it identifies $B_{t-a}$ with the usual Hom between the line bundles. The higher Ext groups vanish.

The derived internal Homs and compositions just used can be computed by replacing $D(a)$ with its ordinary two-term quiver projective resolution. If $J=B_1\subset A$ and $I=ke_0\oplus ke_1$, the resolution of a representation $D$ has terms
\[
 D\otimes_I J\otimes_I A\longrightarrow D\otimes_I A,
 \quad m\otimes z\otimes b\longmapsto mz\otimes b-m\otimes zb .
\]
Both terms carry the diagonal $H$-action. The map is equivariant, and the resolution is projective over $A$. For equivariant right $A$-modules $P,Q$, the action on $\operatorname{Hom}_A(P,Q)$ is $(h f)(m)=\sum h_{(1)}f(\mathsf S(h_{(2)})m)$, with $\mathsf S$ the antipode. Right $A$-linearity is preserved by the antipode identity. Evaluation and composition in their indicated tensor order are equivariant. This construction gives a dg model using the Hom complexes between the resolutions. Tensor it with a dg enhancement of $\mathrm{Perf}(S)$ carrying trivial $H$-actions on morphism complexes. It computes the ordinary Hom complexes on the product for the external products above, with compatible $H$-actions. In particular for an $H$-complex $V$, taking the tensor $V\otimes E$ on the left gives
\[
 \operatorname{RHom}(E',V\otimes E)
       \simeq V\otimes\operatorname{RHom}(E',E)
\]
equivariantly. This is the usual map for finite dimensional $V$ and projective models, and extends to bounded complexes and their cones.

Fix $e$ and take any $C\in\mathrm{Perf}(S)$ right orthogonal to the subsequence of parity $e$. Choose $t\ge\max_i a_i$ of parity $e$, and put $F=C\boxtimes\mathcal O(t)$ and $E_i=M_i\boxtimes\mathcal O(a_i)$. Formula~\eqref{eq:odd-projectives} identifies
\[
 \operatorname{RHom}(E_i,F)
 = \operatorname{RHom}_S(M_i,C)\otimes B_{t-a_i}.
\]
For matching parity the first factor vanishes. Otherwise the second factor is projective. Thus every displayed complex is perfect over $H$.

Set $F_N=F$, and for $i=N,\ldots,1$ define
\[
 V_i=\operatorname{RHom}(E_i,F_i),\qquad
 F_{i-1}=\operatorname{Cone}(V_i\otimes E_i\longrightarrow F_i).
\]
At each step where $V_i$ is perfect over $H$, choose a bounded complex of finite projective $H$-modules mapping quasi-isomorphically to the indicated Hom complex, and use that model in the evaluation. The evaluation triangles are then equivariant in the dg model. Tensoring such a model with an $H$-complex having bounded finite-dimensional cohomology preserves perfection: this follows by the cohomological truncation triangles and the projective tensor assertion above. Inductively all $\operatorname{RHom}(E_j,F_i)$ are perfect over $H$. Indeed the induction step follows by applying $\operatorname{RHom}(E_j,-)$ to the triangle, since its first term is $V_i\otimes\operatorname{RHom}(E_j,E_i)$ and the latter tensor is perfect. The underlying Hom complexes have bounded finite dimensional cohomology by properness. After forgetting the action, $F_0$ is right orthogonal to every $E_i$, by semiorthogonality and the evaluation triangles. Therefore $F_0=0$ by fullness. Now apply $\operatorname{RHom}(F,-)$ to the same triangles. The complex $\operatorname{RHom}(F,F_0)$ is acyclic, also as an $H$-complex. Each $V_i\otimes\operatorname{RHom}(F,E_i)$ is perfect. Induction in the reverse direction shows that $\operatorname{RHom}(F,F)$ is perfect over $H$.

Finally $\operatorname{RHom}(F,F)\simeq\operatorname{RHom}_S(C,C)\otimes B_0$ by~\eqref{eq:odd-projectives}, and $B_0=k$ carries the trivial $H$-action. The complex on the right is a direct sum of shifts of trivial modules. For $C\ne0$ its degree-zero cohomology contains the identity, so one of these summands is $k$. This contradicts the infinite projective dimension of $k$. The right orthogonal on $S$ is zero, proving fullness. For surfaces the line-bundle case gives rationality by \cite[Theorem~1.2]{LZ26}.
\end{proof}

\begin{corollary}\label{cor:iterated-parity}
A full exceptional collection of line bundles on $S\times(\mathbb P^1)^m$ induces a full exceptional collection of line bundles on $S$ for each of the $2^m$ parity vectors of the relative degrees. If $S$ is a complex surface, it is rational.
\end{corollary}
\begin{proof}
Apply Theorem~\ref{thm:parity-descent} successively, regarding one projective line at a time as the second factor. The Picard formula for a product with projective space gives the required external-product form at each step.
\end{proof}

\subsection{Mutation through a product subcategory}

For a collection $(M_i\boxtimes\cO_{\PP^{r-1}}(a_i))$ on $S\times\PP^{r-1}$ we call an isomorphism class of objects of $\Db(S)$ a \emph{base class} if it occurs among the $M_i$.

\begin{lemma}\label{lem:mutation-product}
Let $\cB\subset\Db(S)$ be generated by an exceptional collection, put $\cA=\{A\in\Db(S):\RHom(A,B)=0\text{ for all }B\in\cB\}$, so that $\Db(S)=\langle\cB,\cA\rangle$, and let $\cC=\cB\boxtimes\Db(\PP^{r-1})$. Write $R_\cB$ for the mutation through $\cB$, given by the triangle $R_\cB(A)\to A\to jj^*A$. Then, for $M\in\Db(S)$ and $a\in\Z$,
\begin{equation}\label{eq:mutation-product}
  R_\cC\bigl(M\boxtimes\cO(a)\bigr)\simeq(R_\cB M)\boxtimes\cO(a).
\end{equation}
In particular, if a full exceptional collection of external products on $S\times\PP^{r-1}$ contains generators of $\cC$, moving them to the beginning of the collection by right mutations leaves them unchanged and replaces every other member $M\boxtimes\cO(a)$ by $(R_\cB M)\boxtimes\cO(a)$, an object of $\cA\boxtimes\Db(\PP^{r-1})$.
\end{lemma}

\begin{proof}
As in the proof of Proposition~\ref{prop:transport}, the kernel $P_\cB\boxtimes\cO_{\Delta_{\PP^{r-1}}}$, where $P_\cB$ is the kernel of $jj^*$, represents the projection onto $\cC$; the triangle $R_\cB(M)\to M\to jj^*M$, tensored externally with $\cO(a)$, is therefore the defining triangle of $R_\cC(M\boxtimes\cO(a))$. Successive right mutations through the generators of $\cC$ compose to $R_\cC$ on the remaining members.
\end{proof}

\subsection{The product of a Barlow surface with the projective line}

\begin{example}\label{ex:barlow}
Let $S$ be a generic determinantal Barlow surface in the family of \cite[Theorem~1.1]{BGKS15}, and let $(M_1,\dots,M_{11})$ be its exceptional collection of line bundles, whose orthogonal is a phantom; here $K_0(S)=\Z^{11}$, and $S$ is minimal of general type with $p_g=q=0$, so $K_S$ is nef. On $X=S\times\PP^1$ consider
\[
  \bigl(M_1\boxtimes\cO,\dots,M_{11}\boxtimes\cO,\ M_1\boxtimes\cO(1),\dots,M_{11}\boxtimes\cO(1)\bigr).
\]
This collection is exceptional, its twenty-two classes form a basis of $K_0(X)$, and its orthogonal has zero Hochschild homology, by the projective bundle formula and additivity. Thus neither of these two invariants detects its nonfullness. However, for every $s\in S^\circ_{\bE_S}$ and $t\in\PP^1$,
\[
  h_{P_{(s,t)}}=2>0=\alpha_{P_{(s,t)}},
\]
by Corollary~\ref{cor:product-nonfull}. The collection is known not to be full, since its orthogonal contains the phantom of $S$; the point coefficient certifies this by a strict height gap, computed from the collection alone.
\end{example}

\subsection{Higher projective spaces}

\begin{corollary}\label{cor:higher-projective}
Let $S$ be a smooth connected projective complex surface, $r\ge2$, and suppose that $(M_i\boxtimes\cO_{\PP^{r-1}}(a_i))_{i=1}^N$ is a full exceptional collection of line bundles on $S\times\PP^{r-1}$. Delete all members whose base class occurs exactly $r$ times. If at most $r$ members remain, then $S$ and $S\times\PP^{r-1}$ are rational.
\end{corollary}

\begin{proof}
The line bundles on $\PP^{r-1}$ with vanishing cohomology are those of degrees $-r+1,\dots,-1$. So the occurrences of a base class have strictly increasing degrees with pairwise differences at most $r-1$; there are at most $r$ of them, with consecutive degrees when there are $r$.

Let $q$ be the number of base classes occurring $r$ times, and $m\le r$ the number of remaining members. Order the $q$ complete classes by first occurrence. They are exceptional on $S$, since $r$ consecutive differences of relative degrees cannot all lie in a range of length $r-1$. Their span $\cB$ gives the subcategory $\cB\boxtimes\Db(\PP^{r-1})$, and moving its generators to the front by Lemma~\ref{lem:mutation-product} leaves a full exceptional collection
\[
  \bigl(F_1\boxtimes\cO(a_1),\dots,F_m\boxtimes\cO(a_m)\bigr)\quad\text{in }\cA\boxtimes\Db(\PP^{r-1}),\qquad\Db(S)=\langle\cB,\cA\rangle .
\]
Each $F_i$ is exceptional. Fullness and the projective bundle formula show that $K_0(S)$ is free of finite rank and $m=r\cdot\rk K_0(\cA)$, so $m=0$ or $m=r$. For $m=0$, $\cB=\Db(S)$ as before. Otherwise $K_0(\cA)=\Z$, the rank-one Euler form gives $\chi(F_j,F_i)=\pm1$, and semiorthogonality forces $a_i-a_j\in\{-r+1,\dots,-1\}$ for $i<j$, hence $a_i=a_1+i-1$. If $H\in\cA$ is right orthogonal to $F_1$, then $H\boxtimes\cO(a_1)$ is right orthogonal to every generator, since $R\Gamma(\cO(1-i))=0$ for $1<i\le r$. Hence $H=0$ and $\cA=\langle F_1\rangle$. Apply Theorem~\ref{thm:rational} to the resulting collection on $S$.
\end{proof}

\section{The general case and two gluing maps}\label{sec:remaining-glue}
Let $(E_1,\ldots,E_N)$ be any full exceptional collection on a smooth
connected complex projective surface $X$. Set $\cA=\langle
E_1,\ldots,E_{N-1}\rangle$ and let $M$ be the module on $\cA$ defined
by $A\mapsto\RHom_X(A,E_N)$. Denote the action maps in the Hochschild
gluing sequence by
\[
 \theta_t:\HH^t(\cA)\longrightarrow
             \Ext^t_{\mathrm{Mod}(\cA)}(M,M)\qquad(t>0).
\]
Here the module and its endomorphisms are derived. They are determined
by the restricted exceptional category, the forward morphisms into
$E_N$, and the higher compositions with those morphisms.

\begin{proposition}\label{prop:last-generator}
The surface $X$ is rational if and only if $\theta_3$ is surjective
and $\theta_4$ is injective. More precisely there is an exact sequence
\[
 0\longrightarrow\coker(\theta_3)\longrightarrow
 H^0(X,2K_X)^\vee\longrightarrow\ker(\theta_4)\longrightarrow0 .
\]
\end{proposition}
\begin{proof}
Apply the gluing sequence in the proof of
Theorem~\ref{thm:gluing-hereditary} to
$\Db(X)=\langle\cA,E_N\rangle$ and use
$\HH^{>0}(\Perf(\C))=0$. Its five terms in degrees three and four
give the displayed short exact sequence by top-degree HKR and Serre
duality. A full collection gives $q(X)=0$, so the equivalence follows
from Castelnuovo's criterion.
\end{proof}

For a nonrational surface carrying a full collection, the degree-zero
map from the two centres to $\End(M)$ is surjective, since $\HH^1(X)=0$.
In the hereditary case of Section~\ref{sec:hereditary-new} an
indecomposable summand of $M$ is then a shifted module, which gives the
vanishing of $\Ext^3(M,M)$ used there. For an arbitrary directed dg
category the degree-zero surjection does not determine $\Ext^3(M,M)$,
and higher products enter both $\theta_3$ and $\theta_4$; the ranks of
the exceptional objects and the dimensions of their forward Ext spaces
do not determine these two maps. Proving the surjectivity of $\theta_3$
and the injectivity of $\theta_4$ for one full presentation of $\Db(X)$
would prove Orlov's conjecture for surfaces.

\appendix
\section{Ordinary height}\label{app:ordinary}

We record two consequences of Kuznetsov's localization triangle \cite[Theorems~3.3 and~6.3]{Kuz15}, which describe what the ordinary height does and does not see. Let $V$ be smooth connected projective over an algebraically closed field $k$, and let $\Perf(V)=\langle\cA,\cB\rangle$. Write $\NHH^t(\cB,V)=H^tN_{\id}(\cB,V)$ and $h_{\Kuz}(\cB)=\inf\{t:\NHH^t(\cB,V)\ne0\}$.

\begin{proposition}\label{prop:ordinary}
\begin{enumerate}[(1)]
  \item If $\cA\ne0$, then $\NHH^t(\cB,V)\simeq\HH^{t-1}(\cA)$ for $t\le0$. Consequently $h_{\Kuz}(\cB)>0$ if and only if $\HH^{<0}(\cA)=0$, and $h_{\Kuz}(\cB)=0$ if and only if $\HH^{-1}(\cA)\ne0$ and $\HH^{\le-2}(\cA)=0$.
  \item If $\cB=\Perf(V)$, then $h_{\Kuz}(\cB)=0$.
  \item The comparison map $\gamma_\cB\colon\NHH^0(\cB,V)\to\HH^0(V)=k$ is nonzero if and only if $\cB=\Perf(V)$.
\end{enumerate}
\end{proposition}

\begin{proof}
In the triangle $N_{\id}(\cB,V)\to\HH(V)\to\HH(\cA)$ one has $\HH^j(V)=0$ for $j<0$ and $\HH^0(V)=k$, since $\cO_\Delta$ is a sheaf. The long exact sequence gives (1) for $t<0$, and in degree zero it reads
\[
  0\longrightarrow\HH^{-1}(\cA)\longrightarrow\NHH^0(\cB,V)\xrightarrow{\ \gamma_\cB\ }k\longrightarrow\HH^0(\cA).
\]
The last map sends $1$ to $\id_\cA$. If $\cA\ne0$ it is injective, so $\gamma_\cB=0$ and (1) holds in degree zero; if $\cA=0$ the normal complex is $\HH(V)$, which gives (2) and the remaining case of (3).
\end{proof}

So positive and negative ordinary height both exclude fullness, whereas ordinary height zero is compatible with fullness and with nonfullness; in the latter case nonfullness is visible in the comparison map but not in the initial degree. The example of Appendix~\ref{app:fivefold} shows that a geometric coefficient can then restore a strict gap of initial degrees.

\section{A rational fivefold with ordinary height zero}\label{app:fivefold}

Recall that a \emph{phantom} is a nonzero admissible category $\cP$ with $K_0(\cP)=0$ and $\HH_*(\cP)=0$; it is \emph{universal} if $K_0(\cP\otimes\Perf(W))=0$ for every smooth projective $W$ \cite[Definition~1.9]{GO13}. We combine Krah's universal phantom on a rational surface \cite{Kra24} with a graded Kronecker category from \cite[Example~4.8]{RS22}.

\subsection{A two-object category with negative Hochschild cohomology}

Let
\[
  p\colon Y=\PP_{\PP^2}(\cO\oplus\cO(-3))\to\PP^2,\qquad L=\cO_Y(1),\qquad H=p^*\cO_{\PP^2}(1),
\]
with the quotient convention, so that $p_*L=\cO\oplus\cO(-3)$. The variety $Y$ is rational.

\begin{lemma}\label{lem:Y-collection}
The collection
\begin{equation}\label{eq:Y-collection}
  (\cO_Y,\ L,\ LH,\ LH^2,\ L^2H^4,\ L^2H^5)
\end{equation}
is full and exceptional. With $\cQ=\langle\cO_Y,L\rangle$, $(C_1,C_2,C_3,C_4)=(LH,LH^2,L^2H^4,L^2H^5)$ and $\cC=\langle C_1,\dots,C_4\rangle$, one has $\Db(Y)=\langle\cQ,\cC\rangle$ and
\[
  \RHom_Y(L,\cO_Y)=0,\qquad\RHom_Y(\cO_Y,L)=\C\oplus\C[-2].
\]
\end{lemma}

\begin{proof}
The fibres of $p$ are projective lines and $Rp_*\cO_Y=\cO_{\PP^2}$, so every line bundle is exceptional. The restriction of $L^{-1}$ to each fibre is $\cO_{\PP^1}(-1)$, whence $Rp_*L^{-1}=0$, and $R\Gamma(Y,L)=R\Gamma(\PP^2,\cO\oplus\cO(-3))=\C\oplus\C[-2]$. Orlov's projective bundle formula \cite{Orl93} and Beilinson's collection \cite{Bei78} give the full exceptional collection $(H^{-2},H^{-1},\cO_Y,L,LH,LH^2)$. Since $\omega_Y=L^{-2}H^{-6}$, Serre duality allows the first two objects to be moved to the end as $L^2H^4$ and $L^2H^5$. This proves \eqref{eq:Y-collection}.
\end{proof}

\begin{proposition}\label{prop:Q}
The smooth proper category $\cQ$ has
\[
  \HH^t(\cQ)=\begin{cases}\C,&t=-1,0,1,3,\\0,&\text{otherwise.}\end{cases}
\]
\end{proposition}

\begin{proof}
Write the space of arrows from $\cO_Y$ to $L$ as $V=\C a\oplus\C b$ with $|a|=0$ and $|b|=2$. There are no composable nonidentity arrows, so in a strictly unital minimal model all higher products vanish, and $\cQ$ is the graded Kronecker category of \cite[Example~4.8]{RS22}. Over the vertex algebra $S=\C e_0\oplus\C e_1$ one has $V\otimes_SV=0$, so the relative Hochschild cochain complex is the fibre of
\begin{equation}\label{eq:Q-cochains}
  \C\oplus\C\xrightarrow{\ \delta\ }\RHom_\C(V,V),\qquad\delta(\lambda_0,\lambda_1)=(\lambda_1-\lambda_0)\,\id_V .
\end{equation}
The endomorphism space of $V$ has dimensions one, two, one in internal degrees $-2,0,2$; in the fibre the target is shifted by $[-1]$, so these degrees become $-1,1,3$. The map $\delta$ has rank one, with one-dimensional kernel and cokernel in internal degree zero. This gives the stated cohomology. The class of degree $-1$ is represented by the map $b\mapsto a$, $a\mapsto0$.
\end{proof}

\subsection{The surface phantom and the collection on the product}

Let $X=\Bl_{p_1,\dots,p_{10}}\PP^2$ at ten general points. Write $h$ for the pullback of a line and $E_i$ for the exceptional curves, and define
\begin{equation}\label{eq:krah-divisors}
  D_i=-6h+2\sum_{j=1}^{10}E_j-E_i\quad(1\le i\le10),\qquad D=-19h+6\sum_{j=1}^{10}E_j .
\end{equation}
Krah \cite[Theorem~1.1]{Kra24} constructs the exceptional collection
\begin{equation}\label{eq:krah-collection}
  (B_1,\dots,B_{13})=\bigl(\cO_X,\cO_X(D_1),\dots,\cO_X(D_{10}),\cO_X(D),\cO_X(2D)\bigr),
\end{equation}
whose forward Hom spaces vanish \cite[Section~5]{Kra24} and whose right orthogonal $\cP$ is a universal phantom \cite[Corollary~5.4]{Kra24}. Its ordinary height is four \cite[Remark~5.5]{Kra24}, so Kuznetsov's localization theorem gives
\begin{equation}\label{eq:krah-HH}
  \HH^{<0}(\cP)=0,\qquad\HH^0(\cP)=\C .
\end{equation}
We also use the full exceptional collection
\begin{equation}\label{eq:standard-collection}
  (G_1,\dots,G_{13})=\bigl(\cO_X,\cO_X(E_1),\dots,\cO_X(E_{10}),\cO_X(h),\cO_X(2h)\bigr)
\end{equation}
of \cite[Section~5]{Kra24}; its fullness follows from Orlov's blowup formula \cite{Orl93}.

\begin{proposition}\label{prop:product-collection}
On $Z=X\times Y$, consider the collection, with $i$ increasing within each block,
\begin{equation}\label{eq:Z-collection}
  \bE_Z=\bigl((B_i\boxtimes\cO_Y)_{i},\ (B_i\boxtimes L)_{i},\ (G_i\boxtimes C_1)_{i},\ (G_i\boxtimes C_2)_{i},\ (G_i\boxtimes C_3)_{i},\ (G_i\boxtimes C_4)_{i}\bigr),\qquad1\le i\le13 .
\end{equation}
It is exceptional of maximal length $78$, and its right orthogonal is $\cA\simeq\cP\otimes\cQ$, a universal phantom.
\end{proposition}

\begin{proof}
The decompositions of the factors give
\[
  \Db(Z)=\bigl\langle\Db(X)\otimes\cQ,\ \Db(X)\otimes\cC\bigr\rangle=\bigl\langle\cP\otimes\cQ,\ \langle B_1,\dots,B_{13}\rangle\otimes\cQ,\ \Db(X)\otimes\cC\bigr\rangle,
\]
and the last two components are generated by \eqref{eq:Z-collection}. Derived Hom of external products factors, so the collection is exceptional. Replacing the $B_i$ by the $G_i$ in the first two blocks gives a full exceptional collection of length $13\cdot6=78$, so $K_0(Z)$ is free of rank $78$, and the Euler matrix of \eqref{eq:Z-collection} is upper triangular with diagonal entries one; hence its length is maximal.

The two exceptional objects generating $\cQ$ give a semiorthogonal decomposition of $\cP\otimes\cQ$ into two copies of $\cP$. Additivity yields $K_0(\cA)=0$ and $\HH_*(\cA)=0$, and $\cA\ne0$ since it contains a copy of $\cP$. After tensoring with $\Perf(W)$ the same decomposition gives $K_0(\cA\otimes\Perf(W))=K_0(\cP\otimes\Perf(W))^{\oplus2}=0$, so $\cA$ is universal.
\end{proof}

\subsection{The two height calculations}

\begin{proposition}\label{prop:two-heights}
The collection \eqref{eq:Z-collection} has ordinary height $h_{\Kuz}(\bE_Z)=0$, and its ordinary comparison map in degree zero is the zero map between one-dimensional spaces. For $x\in X$, let $T=\cO_x\boxtimes\cO_Y$ and $K=T^\vee\boxtimes S_Z(T)$. Then
\[
  \alpha_{F_K}=0,\qquad\ph_{F_K}(\bE_Z)=h_{F_K}(\bE_Z)=2 .
\]
\end{proposition}

\begin{proof}
For smooth proper categories, the K\"unneth formula for Hochschild cohomology gives $\HH(\cP\otimes\cQ)\simeq\HH(\cP)\otimes_\C\HH(\cQ)$. By \eqref{eq:krah-HH} and Proposition~\ref{prop:Q}, the only term of total degree $-1$ is $\HH^0(\cP)\otimes\HH^{-1}(\cQ)=\C$, and there are no terms of lower degree. Proposition~\ref{prop:ordinary}(1) gives $h_{\Kuz}(\bE_Z)=0$ and $\NHH^0(\langle\bE_Z\rangle,Z)=\C$, and the comparison map in degree zero vanishes by Proposition~\ref{prop:ordinary}(3), since $\cA\ne0$.

For the coefficient $K$, Theorem~\ref{thm:external} gives
\[
  C_K(Z)=\RHom_Z(T,T)\simeq\RHom_X(\cO_x,\cO_x)\otimes R\Gamma(Y,\cO_Y)\simeq\RHom_X(\cO_x,\cO_x),
\]
so $\alpha_{F_K}=0$. The terminal factor \eqref{eq:external-terminal} vanishes unless the path ends in the first block of \eqref{eq:Z-collection}: for the second block this follows from $\RHom_Y(L,\cO_Y)=0$, and for the last four blocks from $\RHom_Y(C_j,\cO_Y)=0$ in the exceptional collection \eqref{eq:Y-collection}. Since indices increase along a path, every contributing path lies entirely in the first block. In that block the factor on $Y$ is $\cO_Y$, with scalar endomorphisms, so the internal morphisms, both endpoint actions and the differentials reduce to those of the point coefficient for the collection \eqref{eq:krah-collection} on $X$. This collection has vanishing forward Hom spaces, and Lemma~\ref{lem:no-hom} gives pseudoheight and height two.
\end{proof}

Both factors are rational, so $Z$ is a rational fivefold. Here the orthogonal is a phantom by construction; the point of the example is that the coefficient $K$, supported on $(\{x\}\times Y)\times(\{x\}\times Y)$, exhibits this by an excess height of two computed from the collection, while the ordinary height of the same collection is zero.


\end{document}